\documentclass[]{dukedissertation}

\usepackage{amsmath, amssymb, amsfonts, amsthm, mathrsfs}
\usepackage{algpseudocode}
\usepackage{algorithm}
\usepackage{graphicx}
\usepackage{color}
\usepackage[usenames,dvipsnames]{xcolor}
\usepackage{bm}
\usepackage{subfigure}
\usepackage{graphicx}
\usepackage{caption}
\usepackage{mathabx}
\usepackage{multirow}
\usepackage{setspace}
\usepackage[english]{babel}
\usepackage[makeroom]{cancel}
\usepackage{wrapfig}

\newtheorem{theorem}{Theorem}[chapter]

\newtheorem{problem}{Problem}[chapter]
\newtheorem{conjecture}{Conjecture}[chapter]
\newtheorem{lemma}[theorem]{Lemma}
\newtheorem{proposition}[theorem]{Proposition}
\newtheorem{coro}[theorem]{Corollary}
\newtheorem{mydef}{Definition}[section]
\numberwithin{equation}{section}

\newtheorem*{remark}{Remark}
\newtheorem{example}{Example}[section]

\newcommand\RR{\mathbb{R}}

\newcommand\<{\langle}
\newcommand\p{\partial}
\newcommand\na{\nabla}
\newcommand\Ric{\text{Ric}}
\newcommand\R{\text{R}}
\newcommand\n{\textbf{n}}
\newcommand\gbar{\overline{g}}
\newcommand\Rbar{\overline{R}}
\newcommand\nabar{\overline{\na}}
\newcommand\naperp{\na^{\perp}}
\newcommand\Gammabar{\overline{\Gamma}}
\newcommand\ppxi{\frac{\p}{\p x^{i}}}
\newcommand\ppxj{\frac{\p}{\p x^{j}}}
\newcommand\ppxk{\frac{\p}{\p x^{k}}}
\newcommand\ppxl{\frac{\p}{\p x^{l}}}
\newcommand\ppt{\frac{\p}{\p t}}
\newcommand\ppr{\frac{\p}{\p r}}
\newcommand\pptheta{\frac{\p}{\p \theta}}
\newcommand\ppphi{\frac{\p}{\p \phi}}

\newcommand\overlinenaxi{\overline{\na}_{\frac{\p}{\p x^{i}}}}
\newcommand\overlinenaxj{\overline{\na}_{\frac{\p}{\p x^{j}}}}

\newcommand\overlinenat{\overline{\na}_{\frac{\p}{\p t}}}
\newcommand\II{\text{II}}
\newcommand\trace{\text{trace}}

\newcommand\Div{\text{div}}
\newcommand\divg{\text{div}}
\newcommand\adj{\text{adj}}
\newcommand\Red{\color{Black}}
\newcommand\Black{\color{Black}}

\newcommand\Green{\color{Black}}
\newcommand\YellowOrange{\color{Black}}
\newcommand\Peach{\color{Black}}
\newcommand\BlueViolet{\color{Black}}
\newcommand\Blue{\color{Black}}

\newcommand\CornflowerBlue{\color{Black}}

\newcommand\RedViolet{\color{Black}}
\newcommand\RoyalPurple{\color{Black}}
\renewcommand\>{\rangle}

\graphicspath{{./Pictures/}}

\author{Hangjun Xu}
\title{Uniformly Area Expanding Flows in Spacetimes}
\supervisor{Hubert L. Bray}
\department{Mathematics} 

\date{2014} 

\copyrighttext{ All rights reserved except the rights granted by the\\
   \href{http://creativecommons.org/licenses/by-nc/3.0/us/}
        {Creative Commons Attribution-Noncommercial Licence}
}

\member{Mark A. Stern}
\member{Leslie D. Saper}
\member{Lenhard L. Ng}

\usepackage[pdftex, pdfusetitle, plainpages=false, 
				letterpaper, bookmarks, bookmarksnumbered,
				colorlinks, linkcolor=black, citecolor=black,
	         filecolor=black, urlcolor=black]
				{hyperref}

\begin{document}

\maketitle

\abstract

The central object of study of this thesis is inverse mean curvature vector flow of two-dimensional surfaces in four-dimensional spacetimes. Being a system of forward-backward parabolic PDEs, inverse mean curvature vector flow equation lacks a general existence theory. Our main contribution is proving that there exist infinitely many spacetimes, not necessarily spherically symmetric or static, that admit smooth global solutions to inverse mean curvature vector flow. Prior to our work, such solutions were only known in spherically symmetric and static spacetimes. The technique used in this thesis might be important to prove the Spacetime Penrose Conjecture, which remains open today. 

Given a spacetime $(N^{4}, \gbar)$ and a spacelike hypersurface $M$. For any closed surface $\Sigma$ embedded in $M$ satisfying some natural conditions, one can ``steer'' the spacetime metric $\gbar$ such that the mean curvature vector field of $\Sigma$ becomes tangential to $M$ while keeping the induced metric on $M$. This can be used to construct more examples of smooth solutions to inverse mean curvature vector flow from smooth solutions to inverse mean curvature flow in a spacelike hypersurface.



\dedication{Dedicated to my parents: Sihong and Juhua.}

\tableofcontents 
\listoffigures	

\acknowledgements

First of all, I would like to express my deepest gratitude to my advisor, Professor Hubert Bray, for introducing me to the field of mathematical relativity, for suggesting this problem, and for countlessly many inspiring discussions on geometric analysis. I would like to thank him for his guidance and patience every step along the way during my research. Without his constant support, this project wouldn't be possible.

I'm also in debt to my committee member and mentor Professor Mark Stern. I took four courses from him, and I want to thank him for teaching me differential geometry and PDE.

I would also like to thank my other committee members, Professor Leslie Saper, and Lenny Ng, as well as my mentor Professor Richard Hain, for various help and support during the last five years.
 
I want to thank Andrew Goetz for checking certain parts of this thesis, and for teaching me how to draw pictures with mathematica. I would also like to thank Tatsunari Watanabe, Kevin Kordek, Tingran Gao and many others for making my time at Duke so memorable. 

Last but not the least, I would like to thank PCMI for their summer school funding in 2012; American Mathematical Society for their travel support during Fall 2011; Duke University Graduate School for the summer research fellowship in 2012; and MSRI for their funding of the mathematical relativity summer school in 2012.

%
%
%
\chapter{Introduction}
\label{chap:introduction}

\section{Motivation: Mass in General Relativity}
General relativity is the study of large scale structures of the universe. One fundamental object in general relativity is the notion of mass. Pointwise energy density and total mass of a spacetime are both well-defined in general relativity. However, the local mass of a given region in a spacetime (called quasi-local mass), as well as the relationship between local mass and pointwise energy density and total mass of the spacetime are still not very well understood. 

Despite of many attempts in defining the quasi-local mass (e.g. \cite{Bartnik1989, Bartnik2002, Bray2001, BrownYork1992, BrownYork1993, Hawking1968, ShiTam2007}), none of the proposed functionals satisfy all the desired properties. One such natural property is that the total mass of the spacetime should be bounded from below by the mass of a region in it, assuming some positivity condition on the pointwise energy density (e.g. dominant energy condition). 

Given a spacetime $(N^{4}, \gbar)$ and a complete asymptotically flat spacelike hypersurface $M^{3}$ (also called a \emph{slice}) with the induced Riemannian metric $g$. Let $k$ be the second fundamental form of $M$. The triple $(M^{3}, g, k)$ is called a \emph{Cauchy data} of this hypersurface (see Figure \ref{fig:slice in spacetime} below).

\begin{figure}[!htb]
\centering
\includegraphics[scale = .4]{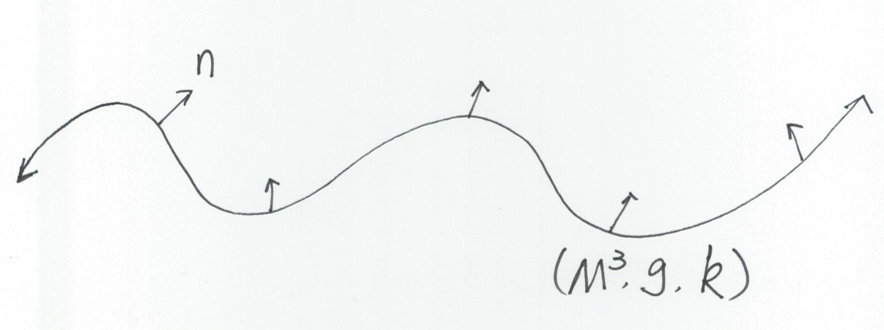}
\caption{Slice in a spacetime with Cauchy data.}
\label{fig:slice in spacetime}
\end{figure}

There is a well-defined quantity called the \emph{ADM mass} (defined by R. Arnowitt, S. Deser and C. Misner in \cite{ADM1956}) that measures the total mass of this hypersurface. Suppose $M$ has a compact outermost minimizing surface $\Sigma$. Physically, $\Sigma$ can be viewed as the apparent horizon of blackholes. 

\begin{figure}[!htb]
\centering
\includegraphics[scale = .4]{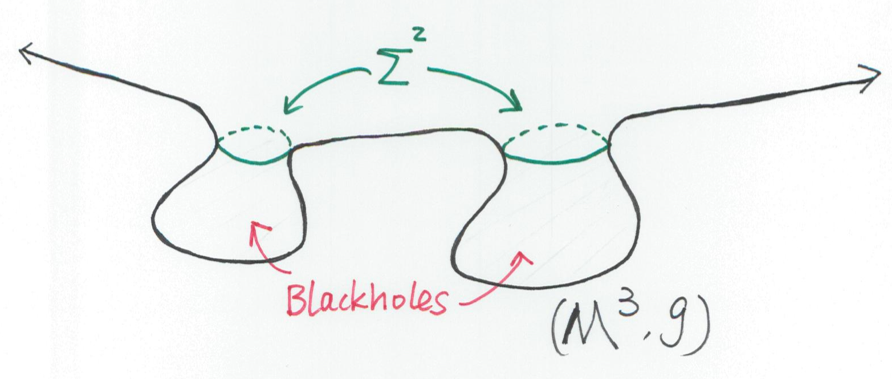}
\caption{Totally geodesic spacelike hypersurface with compact outermost minimal surface $\Sigma$.}
\label{fig:penrose inequality}
\end{figure}

In the case that $M$ is \emph{totally geodesic}, i.e. $k = 0$, then the pointwise energy density equals the scalar curvature of $M$. In this case, the Riemannian Penrose Inequality states that:
\begin{theorem}[Riemannian Penrose Inequality]
Let $M$ and $\Sigma$ be given as above. If the scalar curvature of $(M, g)$ is non-negative, then its ADM mass is greater than or equal to $\sqrt{|\Sigma|/16\pi}$, where $|\Sigma|$ is the total area of $\Sigma$ (see Figure \ref{fig:penrose inequality}).
\end{theorem}

Penrose \cite{Penrose1973} first conjectured this inequality in 1973, and he gave a heuristic proof based on physical considerations, explained as follows. It turns out that the lower bound $\sqrt{|\Sigma|/16\pi}$ in the Riemannian Penrose Inequality equals the \emph{Hawking mass}, which is a quasi-local mass functional proposed by Hawking \cite{Hawking1968}, of the minimal surface $\Sigma$. This can be viewed as the mass of the blackholes inside $\Sigma$. Thus, the Riemannian Penrose Inequality states that the total mass of $M$ should be at least the mass contributed by the blackholes, assuming that the energy density (which equals the scalar curvature in the case of a totally geodesic hypersurface) is non-negative everywhere.

Geroch \cite{Geroch1973}, Jang and Wald \cite{JangWald1977} first discovered a monotone property of Hawking mass of surfaces under smooth \emph{inverse mean curvature flow}. Based on this, Huisken and Ilmanen \cite{HuiskenIlmanen2001} gave a proof of this inequality in the case of a single blackhole (i.e. $\Sigma$ is connected). In the same year, Bray \cite{Bray2001} proved the full Riemannian Penrose Inequality using a different technique. 

In the case of no blackholes, the Riemannian Penrose Inequality is also known as the Riemannian Positive Mass Theorem:

\begin{theorem}[Riemannian Positive Mass Theorem]
Given a complete asymptotically flat Riemannian manifold $(M^{3}, g)$ with non-negative scalar curvature. The ADM mass of $M$ is non-negative. 
\end{theorem}
In 1979, Schoen and Yau \cite{SchoenYau1979} proved this result using a variational method. In the same year, they \cite{SchoenYau19792} generalized this result to Riemannian manifolds of dimension less than eight. In 1981, they \cite{SchoenYau1981} removed the assumption that $M$ is totally geodesic and proved the Riemannian Positive Mass Theorem for an arbitrary spacelike hypersurface in a spacetime that satisfies the dominant energy condition.

\section{Inverse Mean Curvature Vector Flow}
So far all the discussions assume that our spacetime has a hypersurface with zero second fundamental form, in which case inverse mean curvature flow naturally bridges the Hawking mass of an apparent horizon of blackholes and the ADM mass of the hypersurface. However, a spacetime in general does not admit such a totally geodesic hypersurface. This is because that the second fundamental form of a hypersurface in a spacetime has six components, but the hypersurface only has one degree of freedom. Thus it is not generic to have all six components vanish. Therefore it is desirable to obtain a similar bound on the total mass of the spacetime by the mass of the blackholes without this assumption. This leads to the general Spacetime Penrose Conjecture, which is still open today.

A viable candidate for proving this conjecture is the codimension-two analogue of inverse mean curvature flow, called the \emph{inverse mean curvature vector flow}. However, there are two major problems with this flow. First, unlike the inverse mean curvature flow, which is a forward parabolic PDE, inverse mean curvature vector flow is a system of forward-backward parabolic PDEs: forward parabolic in spacelike directions, and backward parabolic in timelike directions (see \cite{HuiskenIlmanen2001}). Backward parabolic equations lack a general existence theory. For instance, the reverse heat flow is backward-parabolic. Given many initial conditions, the reverse heat flow would develop singularities instantaneously. However, the reverse heat flow would exist for time $t > 0$ if we first perform the heat flow for time $t$ and then start flowing backwards.

Second, for some initial surfaces, even the inverse mean curvature vector flow exist, their Hawking mass still won't give us a lower bound on the total mass of the spacetime as in the inverse mean curvature flow case simply because the former is too large. To illustrate this, take a $t = \text{constant}$ slice in the Minkowski spacetime. The round sphere in that slice has zero Hawking mass. Spacial perturbations will decrease the Hawking mass making it negative, whereas timelike perturbations will increase the Hawking mass making it positive. During inverse mean curvature vector flow, the spacial ``wiggles'' will smooth out due to the parabolic nature of the flow. However, timelike ``wiggles'' will get amplified since the flow is reverse parabolic in the timelike directions. With these surfaces with positive Hawking mass, inverse mean curvature vector flow will not provide a lower bound for the ADM mass of Minkowski space, which is zero.

However, these two problems seem to solve each other because they are both suggesting that solutions to inverse mean curvature vector flow exist only when the ``right'' initial surface is given. The important question is then: Given a spacetime. Do such ``right'' initial surfaces always exist? The answer is affirmative if the spacetime is \emph{spherically symmetric} or \emph{static}.

Inverse mean curvature vector flow of surfaces in spherically symmetric spacetimes was first studied by E. Malec, and N. \'{O}Murchadha \cite{MalecOMurchadha1994}. They showed that inverse mean curvature vector flow of spherically symmetric spheres exist for all time. Intuitively, the spherical symmetries prevent timelike ``wiggles'' to occur. Moreover, the Hawking mass is monotonically non-decreasing under inverse mean curvature vector flow of spacelike surfaces with spacelike mean curvature vectors, assuming the spacetime satisfies the dominant energy condition. 

Later Frauendiener \cite{Frauendiener2001} showed that, in an arbitrary spacetime that satisfies the dominant energy condition, if smooth inverse mean curvature vector flow exists, then the Hawking mass is monotonically non-decreasing. In 2004, H. Bray, S. Hayward, M. Mars and W. Simom \cite{BrayHaywardMarsSimon2007} showed that we can in fact flow along a one-parameter family of directions and the Hawking mass is still monotone. 

In spherically symmetric spacetimes, the ``right'' initial surfaces for inverse mean curvature vector flow are spherically symmetric spheres. What about spacetimes that are not necessarily spherically symmetric?

Bray and Ye Li were trying to develop a general existence theory for inverse mean curvature vector flow back in 2009, and one of their intuitions was that if one can somehow control the flow of the surfaces so that they stays tangential to a spacelike slice, then the flow might not develop singularities. In fact it has been shown that:

\begin{proposition}[\cite{BrayHaywardMarsSimon2007}]
\label{prop:imcvf and imvf}
The family of closed embedded spacelike surfaces $\{\Sigma_{s}\}$ is a solution to the smooth inverse mean curvature vector flow with spacelike inverse mean curvature vector $\vec{I}_{\Sigma_{s}}$ everywhere on the surfaces if and only if there exists a spacelike hypersurface $M^{3}\subset N$, such that the mean curvature vector $\vec{H}_{\Sigma_{s}}$ is tangential to $M$ at all $(x, s)$, and $\{\Sigma_{s}\}$ is a solution to the smooth inverse mean curvature flow in $M$.
\end{proposition}

Following the intuition, we prove the main theorem in Chapter \ref{chap:imcvf}:
\begin{theorem}[\textbf{Main Theorem}]\label{thm:main theorem}
\label{thm:existence of spacetimes that admit IMCVF coordinates}
There exist infinitely many non-spherically symmetric, non-static spacetimes that admit inverse mean curvature vector flow coordinate charts. Given such a spacetime $U$ with an inverse mean curvature vector flow coordinate chart $(t, r, \theta, \phi)$ and the constructed spacetime metric $\gbar$. The coordinate spheres $S_{t,r}$ contained in each $t = \text{constant}$ slice, when reparameterized by $r^{2} = e^{s}$, are smooth global solutions to the inverse mean curvature vector flow equation. 
\end{theorem}
This theorem is restated and proved in Theorem \ref{thm:existence of spacetimes that admit IMCVF coordinates}. The proof is based on explicit constructions of spacetime metrics that admit inverse mean curvature vector flow coordinate charts, defined in Chapter \ref{chap:imcvf}. Theorem \ref{thm:main theorem} seems to suggest that spacetimes that admit smooth solutions to inverse mean curvature vector flow exist generically. However, this general problem of find solutions to inverse mean curvature vector flow (i.e. the ``right'' initial surface) in arbitrary spacetimes is still open. 

There also exists a coordinate-free analogue of Theorem \ref{thm:main theorem}:

\begin{theorem}
\label{thm:steering theorem}
Given a spacetime $(N^{4}, \gbar)$,  a spacelike hypersurface $M^{3}$ and a closed embedded surface $(\Sigma^{2}, g_{\Sigma}) \subset M$. Suppose $\Sigma$ is area expanding (defined in (\ref{equ:area expanding condition})), then there exists a unique smooth steering parameter $Q = Q_{\Sigma} \in C^{\infty}(N)$, such that in the steered spacetime metric $\gbar_{Q}$ (defined in (\ref{def:steering of spacetime metric})), $\vec{H}_{\Sigma}$ is tangential to $M$ everywhere on $\Sigma$.
\end{theorem}

This can be used to generate more examples of solutions to inverse mean curvature vector flow. Consider a smooth solution to inverse mean curvature flow in a spacelike hypersurface $M$. One can then smoothly adjust the spacetime metric along the flow such that the mean curvature vector of each flow surface becomes tangential to $M$. These steered surfaces are then solutions to inverse mean curvature vector flow, since the area expanding condition is already satisfied, and now the mean curvature vectors are tangential to a spacelike hypersurface (see Proposition \ref{prop:imcvf and imvf}).

\section{Uniformly Area Expanding Straight Out Flows and Time Flat Surfaces}

Inverse mean curvature vector flow is a type of flow that has bad existence theory, but very good properties: the Hawking mass is monotone under smooth inverse mean curvature vector flow. There is another flow studied in Chapter \ref{chap:straightoutflow}, called \emph{uniformly area expanding straight out flow} (or simply \emph{straight out flow}), that has solutions with a wide class of initial surfaces. In that chapter, we try to construct spacetimes that admit a coordinate chart in which straight out flow of coordinate spheres exists for all time. Partial results have been obtained while complete understanding of this problem is still work in progress:

\begin{proposition}\label{prop:partial results on straight out flows}
Suppose a spacetime $(N^{4}, \gbar)$ admits a coordinate chart $\{t, r, \theta, \phi\}$ such that the coordinate representation of $\gbar$ is
\begin{equation}
\overline{g} = \bordermatrix{~ & t & r & \theta & \phi \cr
                  t & -v^{2} & d & e & f\cr
                  r & d & u^{2} & 0 & 0\cr
                  \theta & e & 0 & a & c\cr
                  \phi & f & 0 & c & b\cr}
\end{equation}
with $ab - c^{2} = r^{4}\sin^{2}{\theta}$. Then $e_{r}: = \frac{1}{u}\ppr$ is straight out if and only if $d$ satisfies a second order elliptic PDE in $d$: $\Delta_{g_{S}}d + G(d, d') = 0$, where $G$ is given by (\ref{equ:scaled straight out equation simplified}) and (\ref{equ:straight out flow equals second order elliptic in d}).
\end{proposition}

We conjecture that solutions to the above elliptic PDE always exist. 

In addition to general existence, another reason for studying the straight out flow is that the Hawking mass is also monotonically non-decreasing under such flow (e.g. see \cite{BrayJauregui2013}, \cite{BrayJauregui2014}) if the spacetime also satisfies the dominant energy condition.

The general existence of straight out flows can serve as a disadvantage since we can even flow surfaces with positive Hawking mass, too large to be used as a lower bound of the total mass of some spacetimes, in straight out directions. To see this, again consider the Minkowski spacetime. All surfaces that are contained in a spacelike plane have non-positive Hawking mass. Thus, non-planer surfaces have positive Hawking mass. For such surfaces, inverse mean curvature vector flow would not work since there are time ``wiggles''. However, those surfaces can still flow in straight out directions. Since the total mass of Minkowski space is zero, having a surface with positive Hawking mass is not going to give us a lower bound for the total mass since the Hawking mass is monotone. 

Inverse mean curvature vector flow and uniformly area expanding straight out flow are two special cases of \emph{uniformly area expanding flows}: orthogonal flows such that the rate of change of the area form of each flow surface equals the area form itself. The Hawking mass is not necessarily monotone under general uniformly area expanding flows. H. Bray, J. Jauregui and M. Mars very recently (\cite{BrayJauregui2013}, \cite{BrayJauregui2014}) obtained a variational formula of the Hawking mass under general uniformly area expanding flows, which consists of four major terms (see \cite{BrayJauregui2014}). The first three terms are non-negative if the spacetime satisfies the dominant energy condition. The fourth term is an integral term with integrant a function of the spacetime multiplied by the divergence of the connection one-form associated with their mean curvature vector of the flow surfaces. Thus, if the connection one-form is divergence free, then the fourth term vanishes and the Hawking mass is monotone. Surfaces with divergence free connection one-form associated with the mean curvature vectors are called \emph{time-flat} (defined in \cite{BrayJauregui2013}, \cite{BrayJauregui2014}). While the conditions on inverse mean curvature vector flow coordinate chart can be viewed as a global ``flatness'' condition on the surface, the time-flat condition is a local ``flatness'' condition.
 
The organization of this thesis is given as follows. In Chapter \ref{chap:hiimcf}, we study the monotonicity of Hawking mass under smooth inverse mean curvature flow. In Chapter \ref{chap:background}, we study inverse mean curvature vector flow in spherically symmetric spacetimes. The notations used in this thesis are also introduced in that chapter. In Chapter \ref{chap:imcvf}, we prove the main theorems \ref{thm:main theorem} and \ref{thm:steering theorem}. In Chapter \ref{chap:straightoutflow}, we study uniformly area expanding straight out flows, and prove Proposition \ref{prop:partial results on straight out flows}. Finally in Chapter \ref{chap:futurework}, some open problems and future works are discussed.
\chapter{Huisken-Ilmanen Inverse Mean Curvature Flow and Monotonicity of Hawking Mass}
\label{chap:hiimcf}

In this chapter we study \emph{inverse mean curvature flow} of a closed embedded surface $\Sigma^{2}$ in an \emph{asymptotically flat} Riemannian manifold $(M^{3}, g)$. The motivation is the Riemannian Penrose Inequality:

\begin{theorem}[Riemannian Penrose Inequality]\label{thm:Riemannian Penrose conjecture}
Let $(M^{3}, g)$ be a complete, asymptotically flat Riemannian manifold with non-negative scalar curvature and a compact outermost minimal surface $\Sigma$ of total area $|\Sigma|$, then
\begin{equation}\label{equ:Riemannian Penrose inequality}
m_{ADM} \geq  \sqrt{\frac{|\Sigma|}{16\pi}}
\end{equation}
with equality if and only if $(M^{3}, g)$ is isometric to the Schwarzschild metric with mass $m > 0$:
\begin{equation}\label{equ:case of equality in Riemannian Penrose inequality}
\left(\RR^{3}\backslash\{0\}, \left(1 + \frac{m}{2|x|}\right)^{4}\delta_{ij}\right)
\end{equation}
outside their respective outermost minimal surfaces. 
\end{theorem}

$\Sigma$ here can be viewed at the apparent horizon of blackholes. The lower bound for the ADM mass, $\sqrt{\frac{|\Sigma|}{16\pi}}$ has the physical interpretation as the mass of the blackholes. Penrose \cite{Penrose1973} first conjectured the Riemannian Penrose Inequality in 1973, and he gave a heuristic proof based on the physical considerations. In 2001, Huisken and Ilmanen \cite{HuiskenIlmanen2001} proved this inequality using inverse mean curvature flow in the case of a single blackhole. In the same year, H. Bray \cite{Bray2001} proved this inequality using conformal flow of metrics that works for any number of blackholes. In 2009, H. Bray and D. Lee \cite{BrayLee2009} generalized the inequality to all dimensions less than eight. In 2010, Lam \cite{Lam2010} proved the Riemannian Penrose Inequality for graphs in all dimensions. In 2011, Schwartz \cite{Schwartz2011} proved a volumetric version of the Penrose inequality for conformally flat manifolds. The general Spacetime Penrose Conjecture is still open today (see \cite{Mars2009, BrayKhuri2010, BrayKhuri2011, Mars2009} for more discussions of this conjecture).

In Section \ref{sec:AF manifolds, ADM mass and Hawking mass}, we define asymptotically flat manifolds, ADM mass and Hawking mass of closed surfaces. In Section \ref{sec:GJW approach monotonicity formula}, we study the Geroch, Jang-Wald monotonicity formula of Hawking mass under smooth inverse mean curvature flow. In Section \ref{sec:HI IMCF with jumps}, we briefly discuss Huisken and Ilmanen's proof of the Riemannian Penrose Inequality using such flows.

\section{Asymptotically Flat Manifolds, ADM Mass and Hawking Mass}
\label{sec:AF manifolds, ADM mass and Hawking mass}
\begin{mydef}\label{def:asymptotically flat manifold}
An $n$-dimensional Riemannian manifold $(M^{n}, g)$ is called \emph{asymptotically flat} if it satisfies the following two conditions:
\begin{enumerate}
\item[(1)] There exists a compact set $K\subset M$ and a diffeomorphism
\[\Phi: E: = M\backslash K \longrightarrow \RR^{n}\backslash \bar{B_{1}},\]
where $B_{1}$ is the unit open ball in $\RR^{n}$; and

\item[(2)] In the coordinate chart $(x^{1}, x^{2}, \cdots, x^{n})$ on $E$ induced by the above diffeomorphism $\Phi$, called an asymptotically flat coordinate chart, the metric components $g_{ij}$ and the scalar curvature $R$ satisfy the following decay conditions at any point $x\in E$, $i, j, k, l = 1, 2, \cdots, n$:
\begin{enumerate}
\item[(1)] $g_{ij}(x) = \delta_{ij}(x) + O(|x|^{-p})$;
\item[(2)] $|x||g_{ij, k}(x)| + |x|^{2}|g_{ij, kl}(x)| = O(|x|^{-p})$;
\item[(3)] $|R(x)| = O(|x|^{-q})$,
\end{enumerate}

with some constants $p > \frac{n-2}{2}$ and $q > n$. Here $g_{ij, k}$ and $g_{ij, kl}$ are coordinate derivatives.
\end{enumerate}
$E$ is called an \emph{asymptotically flat end} of $M$. An asymptotically flat manifold can have multiple asymptotically flat ends.
\end{mydef}

\begin{mydef}\label{def:ADM mass of asymptotically flat manifold}
Given an asymptotically flat Riemannian manifold $(M^{n}, g)$ and asymptotically flat coordinate chart. The ADM mass of of $M$ is:
\begin{equation}\label{equ:AD mass of asymptotically flat manifold}
m_{ADM}(M, g): = \lim_{r\to\infty}\frac{1}{2(n-1)\omega_{n-1}}\int_{S_{r}}\sum_{i, j = 1}^{n}(g_{ij, i} - g_{ii, j})\nu^{j}dS_{r}
\end{equation}
where $\omega_{n-1}$ is the volume of the $(n-1)$-dimensional round sphere; $S_{r}$ is the coordinate sphere of radius $r$; $\nu$ is the outward unit normal along $S_{r}$; and $dS_{r}$ is the volume form of $S_{r}$. 
\end{mydef}
The ADM mass of an asymptotically flat manifold was defined by Richard Arnowitt, Stanley Deser and Charles W. Misner \cite{ADM1956}. They proved that the above definition is independent of the choice of asymptotically flat coordinate charts. Thus, the notion of the ADM mass is well-defined. We sometimes simply write $m_{ADM}(g)$ instead of $m_{ADM}(M, g)$ if the underlying manifold is clear.

In dimension $3$, we have
\begin{equation}
m_{ADM}(M^{3}, g) = \lim_{r\to \infty}\frac{1}{16\pi}\int_{S_{r}}\sum_{i,j = 1}^{3}(g_{ij, i} - g_{ii, j})\nu^{j}dS_{r}.
\end{equation}

\begin{mydef}\label{def:Hawking mass of surface with codim 1}
Given a Riemannian manifold $(M^{3}, g)$ and a closed embedded surface $(\Sigma^{2}, g_{\Sigma})$ with the induced metric. The Hawking mass of $\Sigma$ is defined to be:
\begin{equation}\label{equ:Hawking mass of surface with codim 1}
m_{H}(\Sigma) := \sqrt{\frac{|\Sigma|}{16\pi}}\left(1 - \frac{1}{16\pi}\int_{\Sigma}|H_{\Sigma}|^{2}\,\,dA_{\Sigma} \right),
\end{equation}
where $H_{\Sigma}$ is the scalar mean curvature of $\Sigma$ in $M$. 
\end{mydef}

\begin{example}[Euclidean Space]\label{example:Rn as asymptotically flat manifold}
$\RR^{n}$ with the standard Euclidean metric is an asymptotically flat manifold with zero ADM mass.
\end{example}

\begin{example}[Conformal Transformation of Metric]\label{example:conformal transformation of asymptotically flat manifold}
Given an asymptotically flat manifold $(M^{n}, g)$. Consider a conformal transformation $\gbar = u^{\frac{4}{n-2}}g$ of the metric $g$, with $u \in C^{\infty}(M)$, $u > 0$. By Equation (\ref{conformal transformation of scalar curvature in dimension three and higher}) the scalar curvatures $\Rbar$ and $R$ of $\gbar$ and $g$ respectively, are related by:
\begin{equation}\label{equ:change of scalar curvature under conformal transformation}
\Rbar = u^{-\frac{n+2}{n-2}}\left(Ru  - \frac{4(n-1)}{n-2}\Delta_{g}u\right).
\end{equation}
If $u$ and its coordinate derivatives satisfy the following decay conditions, $i, j, k = 1, 2, \cdots, n$: 
\begin{enumerate}
\item[(1)] $u$ tends to $1$ at $\infty$;
\item[(2)] $u_{,i} = O(|x|^{-p-1})$;
\item[(3)] $u_{,jk} = O(|x|^{-p-2})$;
\item[(4)] $\Delta_{g}u = O(|x|^{-q})$
\end{enumerate}
in an asymptotically flat coordinate chart of $(M^{n}, g)$ for some constants $p > \frac{n-2}{2}$ and $q > n$, then $(M^{n}, \gbar)$ is also asymptotically flat in that coordinate chart. Moreover,

\begin{equation}\label{equ:change of ADM mass under conformal transformation}
m_{ADM}(\gbar) = m_{ADM}(g) - \lim_{r\to\infty}\frac{2}{(n-2)\omega_{n-1}}\int_{S_{r}}\frac{\p u}{\p \nu}dS_{r},
\end{equation}
where $\frac{\p u}{\p \nu}$ is the outward normal derivative of $u$ along $S_{r}$.
\end{example}

\begin{example}[Schwarzschild Manifold]\label{example:Riemannian Schwarzschild metric}
Combining Example \ref{example:Rn as asymptotically flat manifold} and \ref{example:conformal transformation of asymptotically flat manifold}, consider the following one-parameter family of conformal transformations of $(\RR^{n}\backslash\{0\}, \delta_{ij})$, parameterized by a constant $m > 0$:
\begin{equation}\label{def:Riemannian Schwarzschild metric}
\left(\RR^{n}\backslash\{0\}, \left(1 + \frac{m}{2|x|^{n-2}} \right)^{\frac{4}{n-2}}\delta_{ij} \right).
\end{equation}
This is called the \emph{Schwarzschild manifold} of dimension $n$ and \emph{mass} $m$. It is easy to verify that $u = 1 + \frac{m}{2|x|^{n-2}}$ satisfies the desired decay condition to make the resulting metric asymptotically flat. $m$ is called the mass because the ADM mass of this metric is exactly $m$. This can be seen quite easily for the three-dimensional Schwarzschild manifold: $\left(\RR^{3}\backslash\{0\},\left(1 + \frac{m}{2|x|}\right)^{4}\delta_{ij}\right)$. By Equation (\ref{equ:change of ADM mass under conformal transformation}):

\begin{align}
m_{ADM}\left(\left(1 + \frac{m}{2|x|}\right)^{4}\delta_{ij}\right) & = m_{ADM}(\delta) - \lim_{r\to\infty}\frac{1}{2\pi}\int_{S_{r}}\frac{\p}{\p r}\left(1 + \frac{m}{2r} \right)\, dS_{r}\notag\\
& = -\lim_{r\to\infty}\frac{1}{2\pi}\int_{S_{r}}-\frac{m}{2}r^{-2}\,dS_{r}\notag\\
& = \lim_{r\to\infty}\frac{1}{2\pi}\frac{m}{2} r^{-2}4\pi r^{2}\notag\\
& = m \label{equ:ADM mass of Schwarzschild is m}
\end{align}

We now study further the geometry of the three-dimensional Schwarzschild manifold. Let $r: = |x|$. First, note that the Schwarzschild metric is symmetric under the mapping $r\mapsto \frac{m^{2}}{4r}$. Thus the Schwarzschild manifold has two ends, with the center of symmetry being $\frac{m}{2r} = \frac{2r}{m}$, that is $r = \frac{m}{2}$, which is a two-sphere. Recall that if $g = u^{4}\delta$ for some positive function $u$, then the mean curvature of a sphere  of radius $r$ with respect to $g$ is given by:
\begin{equation}\label{equ:mean curvature of sphere under conformal transformation}
H = \frac{1}{u^{2}}\left(\frac{2}{r} + \frac{4}{u}\frac{d u}{d r} \right).
\end{equation}
Therefore at $r = \frac{m}{2}$, the mean curvature is zero. Hence, the sphere $r = \frac{m}{2}$ is a minimal surface, which can be viewed as the apparent horizon of a blackhole. The region outside of the blackhole is called the \emph{exterior region} of the Schwarzschild manifold:
\[\left(\RR^{3}\backslash B_{\frac{m}{2}}, \left(1 + \frac{m}{r}\right)^{4}\delta\right).\]
The exterior region is an asymptotically flat end (see Figure \ref{fig:exterior region of three dimensional schwarzschild})
 
\begin{figure}[!htb]
\centering
\includegraphics[scale = .55]{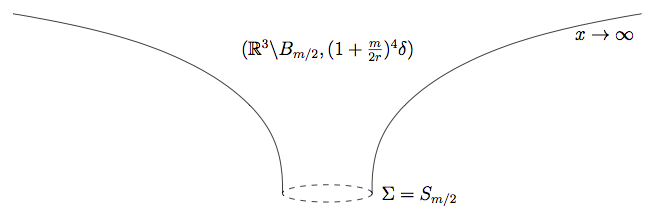}
\caption{Exterior region of three-dimension Schwarzschild manifold with boundary the minimal sphere $S_{\frac{m}{2}}$. Figure courtesy of Mau-Kwong G. Lam.}
\label{fig:exterior region of three dimensional schwarzschild}
\end{figure}

There exists an isometric embedding of the three-dimensional Schwarzschild manifold into $\RR^{4}$ such that 
\[r = \frac{w^{2}}{8m} + 2m.\]
The image of this embedding is a parabola (see Figure \ref{fig:embedding of three dimensional Schwarzschild as parabola}), and the minimal sphere $S_{\frac{m}{2}}$ gets mapped to the sphere $S_{2m}\subset \RR^{4}$:

\begin{figure}[!htb]
\centering
\includegraphics[scale = .55]{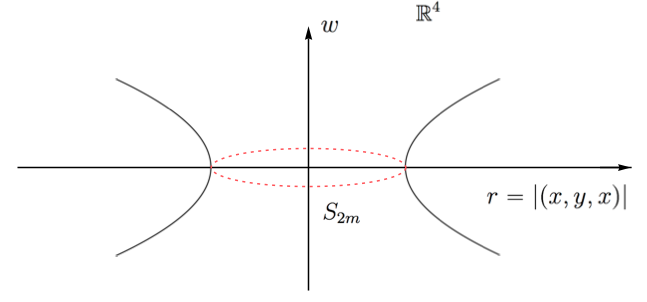}
\caption{Isometric embedding of three-dimensional Schwarzschild manifold into $\RR^{4}$.}
\label{fig:embedding of three dimensional Schwarzschild as parabola}
\end{figure}

The area of $S_{\frac{m}{2}}$ is then given by the Euclidean area of $S_{2m}$: $|S_{\frac{m}{2}}| = 4\pi (2m)^{2} = 16\pi m^{2}$. Therefore
\begin{equation}\label{equ:Hawking mass of minimal sphere in Schwarzschild}
m_{H}(S_{\frac{m}{2}}) = \sqrt{\frac{|S_{\frac{m}{2}}|}{16\pi}}\left(1 - \int_{S_{\frac{m}{2}}}\cancelto{0}{H_{S_{\frac{m}{2}}}^{2}}dA_{S_{\frac{m}{2}}} \right) = m.
\end{equation}

Combing this with the ADM mass (\ref{equ:ADM mass of Schwarzschild is m}), we see that
\begin{proposition}\label{prop:ADM mass of Schwarzschild equals Hawking mass of minimal separating sphere}
The ADM mass of the three-dimensional Schwarzschild manifold $\left(\RR^{3}\backslash\{0\},\left(1 + \frac{m}{2|x|}\right)^{4}\delta_{ij}\right)$ equals the Hawking mass of the minimal sphere $S_{\frac{m}{2}}$, which is exactly $m$.
\end{proposition}
\end{example}

More generally, Huisken and Ilmanen \cite{HuiskenIlmanen2001} proved (see also \cite{Schoen2005}):
\begin{theorem}\label{thm:limit of Hawking mass equals ADM mass of an asymptotically flat manifold}
Given an asymptotically flat Riemannian manifold $(M^{n}, g)$ and an asymptotically flat coordinate chart. Then
\begin{equation}\label{equ:limit of Hawking mass equals ADM mass of an asymptotically flat manifold}
\lim_{r\to\infty}m_{H}(S_{r}) = m_{ADM}(M),
\end{equation}
where $S_{r}$ is the coordinate sphere of radius $r$.
\end{theorem}

\section{Geroch, Jang-Wald's Approach and Their Monotonicity Formula}
\label{sec:GJW approach monotonicity formula}

Huisken and Ilmanen's proof of the Riemannian Penrose Inequality is based on the monotonicity property of the Hawking mass under smooth inverse mean curvature flow, first discovered by Geroch and Jang-Wald.  

\begin{mydef}[Inverse Mean Curvature Flow]\label{def:smooth IMCF}
Given a Riemannian manifold $(M^{3}, g)$ and a closed embedded surface $\Sigma^{2}$ in $M$. A smooth inverse mean curvature flow of $\Sigma$ in $M$ is a smooth family of surfaces
 $F: \Sigma\times[0, T]\longrightarrow M$ of $\Sigma$ such that the following parabolic evolution equation is satisfied:
\begin{equation}\label{equ:smooth IMCF}
\frac{\p F}{\p t} = \frac{\nu_{t}}{H_{t}}, \quad t\in [0, T], \,\,
\end{equation}
where $\nu_{t}$ and $H_{t}$ is the unit outward normal vector field and scalar mean curvature of $\Sigma_{t}: = F(\Sigma, t)$, respectively. 
\end{mydef}
A family of closed embedded surfaces $\{\Sigma_{t}\}$ in $M$ is called a smooth solution to inverse mean curvature flow if they satisfy (\ref{equ:smooth IMCF}). Given such a family of surfaces $\{\Sigma_{t}\}$. The first variation of area formula (\ref{equ:first variation of area}) implies that

\begin{equation}\label{equ:first variation of area in IMCF}
\frac{d}{dt}|\Sigma_{t}| = \int_{\Sigma_{t}}H_{t}\frac{1}{H_{t}}dA_{\Sigma} = |\Sigma_{t}|.
\end{equation}

Therefore the area of $\Sigma_{t}$ grows exponentially under inverse mean curvature flow. 

\begin{example}[Inverse Mean Curvature Flow of Spheres]\label{example:IMCF of sphere}
Consider a round sphere $S_{r_{0}}$ in $\RR^{3}$ with radius $r_{0} > 0$, and flow this sphere out by inverse mean curvature flow. By (\ref{equ:first variation of area in IMCF}), the flow surfaces are still round spheres, and the area grows exponentially. Thus $\{S_{e^{t/2}r_{0}}\}$ is a solution to this flow for all time. 
\end{example}

Geroch \cite{Geroch1973}, Jang and Wald \cite{JangWald1977} discovered the following nice connection between solutions to inverse mean curvature flow and the Hawking mass:

\begin{theorem}[Geroch, Jang-Wald]\label{thm:GJW monotonicity of Hawking mass under IMCF}
Given $(M^{3}, g)$ with non-negative scalar curvature. If a family of closed embedded surfaces $\{\Sigma_{t}\}$ is a smooth solution to inverse mean curvature flow in $M$, then for all $t > 0$,
\begin{equation}
\frac{d}{dt}m_{H}(\Sigma_{t}) \geq 0,
\end{equation}
i.e. the Hawking mass is monotonically non-decreasing.
\end{theorem}

\begin{proof}
Let $H_{t}$ and $dA_{t}$ be the scalar mean curvature and the volume form of $\Sigma_{t}$ in $M$, respectively. 

\begin{align}
& \frac{d}{dt}m_{H}(\Sigma_{t}) = \frac{d}{dt}\left[\sqrt{\frac{|\Sigma_{t}|}{16\pi}}\left(1 - \frac{1}{16\pi}\int_{\Sigma_{t}}H_{t}^{2}\,\,dA_{t} \right)\right]\notag\\
& = \frac{d}{dt}\left( \sqrt{\frac{|\Sigma_{t}|}{16\pi}}\right)\left(1 - \frac{1}{16\pi}\int_{\Sigma_{t}}H_{t}^{2}\,\,dA_{t} \right) + \sqrt{\frac{|\Sigma_{t}|}{16\pi}}\left( - \frac{1}{16\pi}\frac{d}{dt}\int_{\Sigma_{t}}H_{t}^{2}\,\,dA_{t} \right)\notag\\
& = \frac{1}{2}\left(\frac{|\Sigma_{t}|}{16\pi} \right)^{-1/2}\frac{1}{16\pi}|\Sigma_{t}|\left(1 - \frac{1}{16\pi}\int_{\Sigma_{t}}H_{t}^{2}\,\,dA_{t} \right)\tag{By Equation (\ref{equ:first variation of area in IMCF})}\notag\\
& + \sqrt{\frac{|\Sigma_{t}|}{16\pi}}\left(-\frac{1}{16\pi} \int_{\Sigma_{t}}\left[2H_{t}\frac{d}{dt}(H_{t})dA_{t} +  H_{t}^{2}\frac{d}{dt}(dA_{t})\right]\right)\notag\\
& = \frac{1}{2}\sqrt{\frac{|\Sigma_{t}|}{16\pi}}\left(1 - \frac{1}{16\pi}\int_{\Sigma_{t}}H_{t}^{2}\,\,dA_{t} \right) + \sqrt{\frac{|\Sigma_{t}|}{16\pi}}\left(-\frac{1}{16\pi} \int_{\Sigma_{t}}\left[2H_{t}\frac{d}{dt}(H_{t})dA_{t} +  H_{t}^{2}dA_{t}\right]\right)\notag\\
& = \sqrt{\frac{|\Sigma_{t}|}{16\pi}}\left\{\frac{1}{2}\left(1 - \frac{1}{16\pi}\int_{\Sigma_{t}}H_{t}^{2}\,\,dA_{t} \right) - \frac{1}{16\pi} \int_{\Sigma_{t}}\left[2H_{t}\frac{d}{dt}(H_{t}) +  H_{t}^{2}\right]dA_{t}\right\}\label{equ:variation of Hawking mass along IMCF step 1}
\end{align}
By the first variation of mean curvature formula (\ref{equ:first variation of mean curvature}):
\begin{equation}\label{equ:variation of mean curvature along IMCF}
\frac{d}{dt}H_{t} = -\Delta_{\Sigma_{t}}\left(\frac{1}{H_{t}}\right) - \frac{1}{H_{t}}\Ric^{M}(\nu, \nu) -\frac{1}{H_{t}} ||\II_{t}||^{2},
\end{equation}
where $\Ric^{M}$ is the Ricci curvature of $M$. Plug (\ref{equ:variation of mean curvature along IMCF}) into (\ref{equ:variation of Hawking mass along IMCF step 1}) we get:

\begin{align}
& \frac{d}{dt}m_{H}(\Sigma_{t}) = \sqrt{\frac{|\Sigma_{t}|}{16\pi}}\bigg\{\frac{1}{2} - \frac{1}{2}\frac{1}{16\pi}\int_{\Sigma_{t}}H_{t}^{2}\,\,dA_{t}\notag\\
& - \frac{1}{16\pi} \int_{\Sigma_{t}}2H_{t}\left[-\Delta_{\Sigma_{t}}\left(\frac{1}{H_{t}}\right) - \frac{1}{H_{t}}\Ric^{M}(\nu, \nu) -\frac{1}{H_{t}} ||\II_{t}||^{2} \right]+  H_{t}^{2}dA_{t}\bigg\}\notag\\
& = \sqrt{\frac{|\Sigma_{t}|}{16\pi}}\bigg\{\frac{1}{2} + \frac{1}{16\pi}\int_{\Sigma_{t}} \left[ 2H_{t}\Delta_{\Sigma_{t}}\left(\frac{1}{H_{t}} \right) + 2\Ric^{M}(\nu, \nu) + 2||\II_{t}||^{2} - \frac{3}{2}H_{t}^{2}\right] dA_{t}\bigg\}\label{equ:variation of Hawking mass along IMCF step 2}
\end{align}

We now compute the first three integral terms in the above. By integration by parts, we get:
\begin{equation}
\int_{\Sigma_{t}}2H_{t}\Delta_{\Sigma_{t}}\left(\frac{1}{H_{t}} \right)dA_{t} = \int_{\Sigma_{t}}-2\<\na_{\Sigma_{t}}H, \na_{\Sigma_{t}}\frac{1}{H_{t}}\>dA_{t} = \int_{\Sigma_{t}}\frac{2||\na_{\Sigma_{t}}H_{t}||^{2}}{H_{t}^{2}}dA_{t}.
\end{equation}

Now by the Gauss equation (see e.g. \cite{Jost2008}), we have
\begin{equation}
\Ric^{M}(\nu, \nu) = \frac{1}{2}R^{M} - K^{\Sigma_{t}} + \frac{1}{2}H_{t}^{2} - \frac{1}{2}||\II_{t}||^{2},
\end{equation}
where $K^{\Sigma_{t}}$ is the Gauss curvature of $\Sigma_{t}.$

Next let $\lambda_{1}(t)$ and $\lambda_{2}(t)$ be the principal curvatures of $\Sigma_{t}$, then 
\[H_{t} = \lambda_{1}(t) + \lambda_{2}(t), \quad ||\II_{t}||^{2} = \lambda_{1}(t)^{2} + \lambda_{2}(t)^{2}.\]
Therefore
\begin{equation}
||\II_{t}||^{2} - \frac{1}{2}H_{t}^{2} = \lambda_{1}(t)^{2} + \lambda_{2}(t)^{2} - \frac{[\lambda_{1}(t) + \lambda_{2}(t)]^{2}}{2} = \frac{[\lambda_{1}(t) - \lambda_{2}(t)]^{2}}{2}.
\end{equation}

Plug them back into (\ref{equ:variation of Hawking mass along IMCF step 2}), we get:
\begin{align}
& \frac{d}{dt}m_{H}(\Sigma_{t}) = \sqrt{\frac{|\Sigma_{t}|}{16\pi}}\bigg\{\frac{1}{2} + \frac{1}{16\pi}\int_{\Sigma_{t}} \left[ 2H_{t}\Delta_{\Sigma_{t}}\left(\frac{1}{H_{t}} \right) + 2\Ric^{M}(\nu, \nu) + 2||\II_{t}||^{2} - \frac{3}{2}H_{t}^{2}\right] \bigg\}\notag\\
& = \sqrt{\frac{|\Sigma_{t}|}{16\pi}}\bigg\{\frac{1}{2} + \frac{1}{16\pi}\int_{\Sigma_{t}} \frac{2||\na_{\Sigma_{t}}H_{t}||^{2}}{H_{t}^{2}} + R^{M} - 2K^{\Sigma_{t}} + H_{t}^{2} + ||\II_{t}||^{2} - \frac{3}{2}H_{t}^{2} dA_{t}\bigg\}\notag\\
& \geq \sqrt{\frac{|\Sigma_{t}|}{16\pi}}\bigg\{\frac{1}{2} + \frac{1}{16\pi}\int_{\Sigma_{t}}  - 2K^{\Sigma_{t}} +\frac{(\lambda_{1}(t) - \lambda_{2}(t))^{2}}{2}  dA_{t}\bigg\}\tag{$R^{M} \geq 0$}\notag\\
& \geq \sqrt{\frac{|\Sigma_{t}|}{16\pi}}\bigg\{\frac{1}{2} - \frac{1}{8\pi}\int_{\Sigma_{t}}  K^{\Sigma_{t}} dA_{t}\bigg\}\notag\\
& \geq 0
\end{align}
where the last inequality follows from the Gauss-Bonnet formula:
\[\int_{\Sigma_{t}}  K^{\Sigma_{t}} dA_{t} \leq 2\pi \chi(\Sigma_{t}) = 2\pi(2 - 2\cdot \text{genus}(\Sigma_{t})) \leq 4\pi.\]
\end{proof}

Using this, Geroch, Jang-Wald discovered a possible approach to prove the Riemannian Penrose Inequality via the following steps:

$\bullet$ Let $\Sigma$ be the outermost minimal surface in $M$. Its Hawking mass is $\sqrt{\frac{|\Sigma|}{16\pi}}$ since its mean curvature is zero. Notice that this is the lower bound in the Riemannian Penrose Inequality.

$\bullet$ Flow $\Sigma$ out by inverse mean curvature flow, and \emph{assume} that the flow is smooth and exists for all time. Let $\{\Sigma_{t}\}$ be the flow surfaces. Theorem \ref{thm:GJW monotonicity of Hawking mass under IMCF} implies that the Hawking mass is non-decreasing. 

$\bullet$ Let $S_{r}$ be the coordinate sphere of radius $r$ in an asymptotically flat coordinate chart of $M$. Theorem \ref{thm:limit of Hawking mass equals ADM mass of an asymptotically flat manifold} implies that $\displaystyle\lim_{r\to\infty}m_{H}(S_{r}) = m_{ADM}(M)$.

Here is the upshot: If smooth inverse mean curvature flow of $\Sigma$ in $M$ exists for all time, and the flow surfaces approach large coordinate spheres near infinity sufficiently fast, then:

\begin{equation}
m_{ADM}(M) = \lim_{t\to \infty}m_{H}(\Sigma_{t}) \geq m_{H}(\Sigma_{0}) = m_{H}(\Sigma) = \sqrt{\frac{|\Sigma|}{16\pi}},
\end{equation}
and this would prove the Riemannian Penrose Inequality. However, smooth solutions to inverse mean curvature flow do not always exist. In fact, in the case that the flow surface becomes minimal (i.e. mean curvature is zero), the flow is no longer defined since the flow speed is the reciprocal of the mean curvature (see (\ref{equ:smooth IMCF})). There are other cases where singularities can occur:

\begin{example}[Inverse Mean Curvature Flow of Disjoint Spheres]\label{example:IMCF of disjoint spheres}
Suppose $\Sigma$ is a disjoint union of two spheres. Inverse mean curvature flow of $\Sigma$ will develop self-intersection in finite time. 
\end{example}

\begin{example}[Inverse Mean Curvature Flow of Thin Torus]\label{example:IMCF of thin torus}
Consider a thin torus in $\RR^{3}$, obtained as the boundary of an $\epsilon$-neighborhood of a large round circle. Thus its mean curvature is positive everywhere. Now starting flowing the torus by inverse mean curvature flow (see Figure \ref{fig:IMCF_thin_torus}). By the first variation of the mean curvature (\ref{equ:first variation of mean curvature}) and the parabolic maximum principle, the flow speed has a lower bound. As a consequence the torus will fatten up and eventually the mean curvature will become negative in the hole of the torus. Thus, the mean curvature must be zero at some point. However, the flow is not defined when the mean curvature is zero. 
\begin{figure}[!htb]
\centering
\includegraphics[scale = .55]{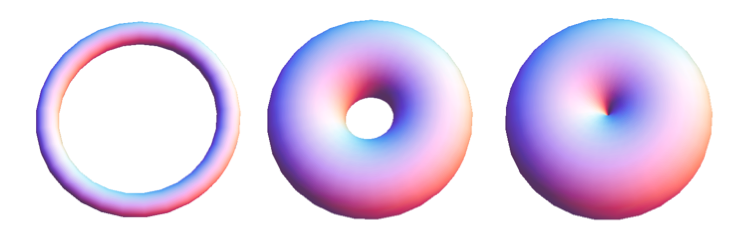}
\caption{Inverse mean curvature flow of a thin torus which develops a singularity in finite time. Picture courtesy of Andrew Goetz.}
\label{fig:IMCF_thin_torus}
\end{figure}
\end{example}


\section{Huisken-Ilmanen's Approach and Level Set Formulation of Inverse Mean Curvature Flow with Jumps}
\label{sec:HI IMCF with jumps}

Because inverse mean curvature flow does not always have solutions, Huisken and Ilmanen defined a generalized inverse mean curvature flow which \emph{always} has solutions. The basic idea is that, in this generalized flow, when a surface is enclosed by another surface of less area, it \emph{jumps} outward to its outermost minimal area enclosure (see \cite{HuiskenIlmanen2001, Bray2001, Bray2002}), and then resume inverse mean curvature flow. Huisken and Ilmanen used a level set formulation to characterize this jumping phenomenon. Consider a scalar-valued function $f$ on $M$, and let $\Sigma_{t}$ be the level set of $f$:
\[\Sigma_{t} = \{x\in M | f(x) = t\}.\]
In this setting, the inverse mean curvature flow equation (\ref{equ:smooth IMCF}) becomes:

\begin{equation}\label{equ:level set formulation of IMCF}
\divg\left(\frac{\na f}{|\na f|} \right) = |\na f|.
\end{equation}
In the above, notice that the left hand side is the mean curvature of $\Sigma_{t}$, and the right hand side is the reciprocal of the flow speed. Thus when $|\na f| \neq 0$, equation (\ref{equ:level set formulation of IMCF}) describes inverse mean curvature flow of the level sets. Note that this formulation allows jumps in a natural way, since if $f$ is constant on some region of $M$, then the level sets of $f$ would just jump over that region during the flow. Huisken and Ilmanen defined a weak solution to (\ref{equ:level set formulation of IMCF}) using an energy minimization principle, and proved existence of such a weak solution by regularizing the degenerate elliptic equation (\ref{equ:level set formulation of IMCF}). They showed that the Hawking mass is still monotone as in the smooth inverse mean curvature flow case. In particular, the Hawking mass is non-decreasing during jumps. In this way, they proved the Riemannian Penrose Inequaltiy (\ref{thm:Riemannian Penrose conjecture}) in the case of a single blackhole (i.e. the outermost minimal surface $\Sigma$ is connected).
\chapter{Inverse Mean Curvature Vector Flow in Spherically Symmetric Spacetimes}
\label{chap:background}

In this chapter, we study inverse mean curvature vector flows in spherically symmetric spacetimes. The motivation comes from the fact that despite lack of general existence theory, inverse mean curvature vector flow always works in spherical symmetry. In Section \ref{sec:technical background} and \ref{sec:hypersurface of spacetime}, terminologies and notations used in later discussions are provided. In Section \ref{sec:codim two surface geometry}, we study closed embedded codimension-two surfaces in a spacetime and the geometry of their normal bundles. We define mean curvature vector fields and Hawking mass. We then define inverse mean curvature vector flow of a closed embedded surface. In Section \ref{sec:spherically symmetric spacetime}, we show that spherically symmetric spheres are smooth global solutions to inverse mean curvature vector flow, and the Hawking mass is monotonically non-decreasing if the spacetime also satisfies the dominant energy condition. 

\section{Spacetime, Einstein Equation, Dominant Energy Condition}\label{sec:spacetime geometry and dominant energy condition}
\label{sec:technical background}
\begin{figure}[!htb]
\centering
\includegraphics[scale = .3]{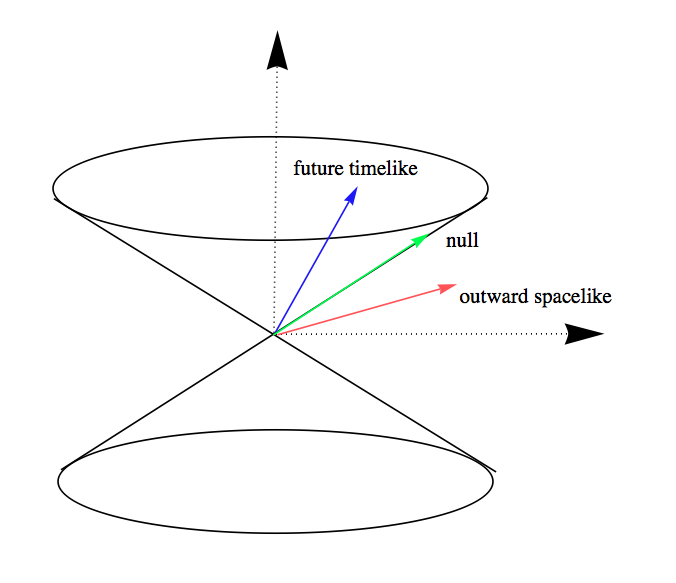}
\caption{Timelike, null and spacelike vectors in spacetime.}
\label{fig:spacetime geometry}
\end{figure}

A \emph{spacetime} $(N^{4}, \gbar, \nabar)$ considered in this thesis is a connected, smooth, time-oriented, four-dimensional manifold with Lorentian metric $\gbar$ of signature $(-, +, +, +)$. $\nabar$ is the associated Levi-Civita connection. A tangent vector of $N$ is called (see Figure \ref{fig:spacetime geometry})
\begin{enumerate}
\item \emph{timelike} if $\gbar(v, v) < 0$;
\item \emph{null} if $\gbar(v, ,v) = 0$;
\item \emph{spacelike} if $\gbar(v, ,v) > 0$.
\end{enumerate}

$N$ is called \emph{time-orientable} if it admits a smooth timelike vector field $\vec{T}$. $N$ is \emph{time-oriented} if such a vector field $\vec{T}$ is chosen. A vector field $X$ is called \emph{future-pointing} if $\gbar(X, \vec{T}) > 0$, or \emph{past-pointing} if $\gbar(X, \vec{T}) < 0$.

$\<\cdot, \cdot\> := \<\cdot, \cdot\>_{\gbar}$ is used to denote the inner product with respect to $\gbar$, unless otherwise specified.

A vector field $X\in \Gamma(TN)$ is said to be \emph{timelike} (resp. \emph{null} or \emph{spacelike}) if at every point $p\in N$, $X(p)$ is timelike (resp. null or spacelike). A submanifold $M$ of $N$ is said to be \emph{timelike} (resp. \emph{null} or \emph{spacelike}) if every tangent vector of $M$ is timelike (resp. null or spacelike). 

Let $\Ric^{N}$ and $\R^{N}$ be the Ricci curvature and scalar curvature of the spacetime respectively. The \emph{Einstein curvature tensor} $G$ is defined as:
\begin{equation}\label{def:Einstein curvature}
G: = \Ric^{N} - \frac{1}{2}\R^{N}\cdot \gbar.
\end{equation}
We assume that the spacetime satisfies the \emph{Einstein Equation}:
\begin{equation}\label{equ:Einstein equation}
G = 8\pi T,
\end{equation}
where $T$ is the \emph{stress energy tensor}. For any tangent vectors $u, v$ of $N$, $T(u, v)$ has the physical meaning as the energy density going in the direction of $u$ as observed by someone going in the direction of $v$. The \emph{dominant energy condition} is:
\begin{equation}\label{equ:dominant energy condition}
T(u, v) \geq 0, \quad \forall\,\, u, v\,\,\text{future-pointing, timelike}.
\end{equation}

\section{Spacelike Hypersurfaces}\label{sec:hypersurface of spacetime}
Given a spacetime. Consider a spacelike hypersurface $M$ with globally defined, future-timelike unit normal vector field $\n$ on $M$, and induced Riemannian metric $g = g^{M}$ by restricting the spacetime metric onto $M$. $M$ is also called a \emph{slice} in a spacetime. Let $k$ be the second fundamental form of $M$, then the triple $(M^{3}, g, k)$ is called a \emph{Cauchy data}. Given a Cauchy data and the unit normal vector field $\n$, we define the \emph{energy density} of $M$ as $\mu: = T(\n, \n)$.
We can compute that
\begin{align}
\mu & = \frac{1}{8\pi}G(\n, \n) = \frac{1}{8\pi}\left(\Ric^{N}(\n, \n) - \frac{1}{2}\R^{N} \cdot\gbar(\n,\n) \right)\tag{By the Einstein equation}\notag\\
& = \frac{1}{8\pi}\left(\Ric^{N}(\n, \n) + \frac{1}{2}\R^{N} \right) \tag{$\n$ is unit time like}\notag\\
& = \frac{1}{16\pi}\left(R^{M} + (\trace_{g}k)^{2} - ||k||^{2}_{g} \right).\label{energy density equation}
\end{align}
where $R^{M}$ is the scalar curvature of $(M, g)$. 
We define the \emph{momentum density} of $M$ to be a one-form $J(\cdot)$ on $M$, such that $J(X) := T(\n, X)$, for any tangent vector field $X$ on $M$.  Then:
\begin{align}
J(\cdot) & = T(\n, \cdot) = \frac{1}{8\pi}\left(\Ric^{N}(\n, \cdot) - \frac{1}{2}R^{N}\gbar(\n, \cdot) \right)\notag\\
& = \frac{1}{8\pi}\Ric^{N}(\n, \cdot)\tag{$\n$ is normal to $M$}\notag\\
& = \frac{1}{8\pi}\divg_{g}((k - \trace_{g}k)\cdot g). \label{momentum density equation}
\end{align}
Equation $(\ref{energy density equation})$ and $(\ref{momentum density equation})$ follow from the Gauss and Codazzi equations respectively, and they are called the \emph{constraint equations}. The dominant energy condition (\ref{equ:dominant energy condition}) implies that
\begin{equation}
\mu \geq ||J||_{g}.
\end{equation}
In the time-symmetric case where the second fundamental form $k = 0$, we see that $\mu = \frac{R^{M}}{16\pi}$, and $J = 0$. The dominant energy condition thus reduces to $R\geq 0$.

\section{Geometry of Codimension Two Surfaces in Spacetime}
\label{sec:codim two surface geometry}
Let $\Sigma$ be an \emph{closed, embedded, spacelike} surface in $N$ with codimension two. We assume that $\Sigma$ is an oriented surface such that at each point the notion of ``outward'' and ``inward'' is well-defined. Let $g_{\Sigma}$ be the induced metric on $\Sigma$ from the spacetime metric $\gbar$.

\subsection{Rank-two Normal Bundle Geometry: Normal Connection and Connection One-form}\label{subsec:normal bundle geometry}
Let $N\Sigma$ be the rank-two normal bundle of $\Sigma$. Define a connection on $N\Sigma$, denoted as $\na^{\perp}$, to be the projection of $\nabar$ onto $N\Sigma$. Notice that since $\Sigma$ is spacelike, $N\Sigma$ has an induced metric with signature $(-, +)$. Therefore, at each point $p\in \Sigma$, $N_{p}\Sigma$ has four quadrants: future-timelike, past-timelike, outward-spacelike and inward-spacelike. Given a local orthonormal frame $\{e_{1}, e_{2}\}$ of $N\Sigma$. Suppose $e_{1}$ is outward-spacelike and $e_{2}$ is future-time like, then the geometry of $N_{p}\Sigma$ is depicted by Figure \ref{fig:codim2 surface}:

\begin{figure}[!htb]
\centering
\includegraphics[scale = .45]{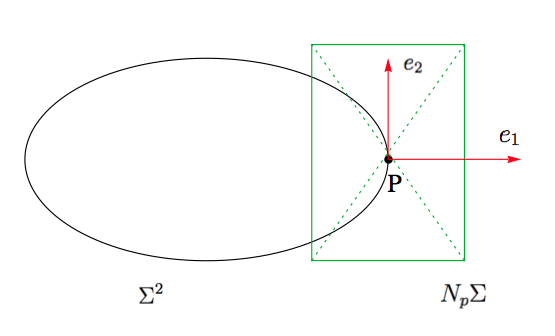}
\caption{Rank two normal bundle at a point $p$ on a surface with orthonormal basis $e_{1}(p)$ and $e_{2}(p)$.}
\label{fig:codim2 surface}
\end{figure}

Define a linear isomorphism, denoted as ``$\perp$'', on each fiber $N_{p}\Sigma$ of the normal bundle as follows: for any orthonormal basis $\{u, v\}$ of $N_{p}\Sigma$ such that $u$ is outward-spacelike and $v$ is future-timelike, set $u^{\perp}: = v$, and then extend linearly to the entire fiber. This definition is independent of the choice of basis, and is an involutive isomorphism. Notice that
\[\<u^{\perp}, u^{\perp}\> = \<v, v\> = -1 = -\<u, u\>.\]
This isomorphism can be viewed as analogy of the $90^{\circ}$ rotation in Euclidean space (see \cite{BrayJauregui2013}). 

Given any tangent vector field $X$ of $\Sigma$. If $v$ is an outward-spacelike unit normal vector field along $\Sigma$, then:
\[0 = \naperp_{X}\<v, v\> = 2\<\naperp_{X}v, v\>.\]
Therefore $\naperp_{X}v$ is perpendicular to $v$, and hence is parallel to $v^{\perp}$. From this one can define a one-form $\alpha_{v}$ on $\Sigma$, uniquely depends on $v$, such that:
\begin{equation}\label{def:connection one form}
\alpha_{v}(X): = \<\naperp_{X}v, v^{\perp}\>, \quad \forall\, X\in \Gamma(T\Sigma).
\end{equation}

This definition yields the following straightforward corollary (see \cite{BrayJauregui2013}):
\begin{coro}
For a smooth section $v\in \Gamma(N\Sigma)$, the associated connection one-form $\alpha_{v}$ vanishes if and only if $v$ is parallel with respect to the normal connection $\naperp$, i.e,
\[\naperp_{X}v = 0,  \quad \forall\, X\in \Gamma(T\Sigma).\]
\end{coro}

\subsection{Mean Curvature Vector Field, Hawking Mass and Inverse Mean Curvature Vector Flow}\label{subsec:mean curvature vector, hawking mass, normal variation}
We define
\begin{equation}\label{def:second fundamental form}
\vec{\II}: T\Sigma \times T\Sigma \longrightarrow N\Sigma, \quad (X, Y) \mapsto (\nabar_{X}Y)|_{N\Sigma}
\end{equation}
to be the \emph{second fundamental form} of $\Sigma$, where $X, Y$ are tangent vector fields along $\Sigma$, and $(\nabar_{X}Y)|_{N\Sigma}$ is the projection of $\na_{X}Y$ onto the normal bundle of $\Sigma$. Define the \emph{mean curvature vector field} of $\Sigma$ to be the trace of the second fundamental form with respect to $g_{\Sigma}$:

\begin{equation}\label{def:mean curvature vector}
\vec{H}_{\Sigma}: = \trace_{g_{\Sigma}}\vec{\II}.
\end{equation}

Therefore, $\vec{H}_{\Sigma}$ is an \emph{inward-pointing} normal vector field along $\Sigma$. Given any normal vector field $\vec{n}$, define the \emph{mean curvature scalar in the direction of $\vec{n}$} to be:
\begin{equation}\label{def:mean curvature scalar}
H_{\vec{n}}: = -\<\vec{H}_{\Sigma}, \vec{n}\>.
\end{equation}

Given $(\Sigma, g_{\Sigma})$ and the mean curvature vector $\vec{H}_{\Sigma}$, the \emph{Hawking mass} of $\Sigma$ is defined to be:
\begin{equation}\label{def:hawking mass}
m_{H}(\Sigma): = \sqrt{\frac{|\Sigma|}{16\pi}}\left(1 - \frac{1}{16\pi}\int_{\Sigma}\gbar(\vec{H}_{\Sigma}, \vec{H}_{\Sigma})\,dA_{\Sigma} \right),
\end{equation}
where $|\Sigma|$ is the area of the surface $\Sigma$. 

\begin{mydef}\label{def:inverse mean curvature vector and smooth IMCVF}
Given a spacetime $(N, \gbar)$ and a surface $(\Sigma^{2}, g_{\Sigma})$. Define the \emph{inverse mean curvature vector field} of $\Sigma$ to be:
\begin{equation}
\vec{I}_{\Sigma}: = - \frac{\vec{H}_{\Sigma}}{\<\vec{H}_{\Sigma}, \vec{H}_{\Sigma}\>}.
\end{equation}
\end{mydef}

According to our sign convention for $\vec{H}_{\Sigma}$, $\vec{I}_{\Sigma}$ thus defined is \emph{outward-pointing}. 

\begin{mydef}[Smooth Inverse Mean Curvature Vector Flow]\label{def:smooth IMCVF}
Given a closed embedded surface $\Sigma^{2}$ in a spacetime $(N^{4}, \gbar)$. A \emph{smooth inverse mean curvature vector flow} of $\Sigma$ is a smooth family of surfaces $F: \Sigma\times [0, T]\longrightarrow N$ of $\Sigma$ satisfying the following evolution equation:

\begin{equation}\label{equ:IMCVF equation}
\frac{\p}{\p s}F(x, s) = \vec{I}_{\Sigma_{s}}(x, s), \quad s\in [0, T]\,\, \text{and}\,\,(x, s)\in \Sigma_{s}: = F(\Sigma, s).
\end{equation}
\end{mydef}

$T > 0$ could also be $\infty$. By the first variation of area formula (see Equation \ref{equ:first variation of area} in appendix \ref{append:first variation of area}), the rate of change of area form of the flow surfaces under smooth inverse mean curvature vector flow is given by:

\begin{equation}\label{equ:first variation of area in IMCVF}
\frac{d}{ds}dA_{\Sigma_{s}} = -\<\vec{H}_{\Sigma_{s}}, \vec{I}_{\Sigma_{s}}\>dA_{\Sigma_{s}} = dA_{\Sigma_{s}}.
\end{equation}

That is, the rate of the area form of each surface is the area form itself, everywhere on each surface. This is a special case of a \emph{uniformly area expanding flow} first defined by H. Bray, J. Jauregui and M. Mars in \cite{BrayJauregui2014}.

\section{Model Spacetime: Spherically Symmetric Spacetime}\label{sec:spherically symmetric spacetime}
In this section, we study spherically symmetric spacetimes. The main motivation comes from the fact that, even though inverse mean curvature vector flow lacks a general existence theory, smooth solutions still exist in many spherically symmetric spacetimes. Moreover, the Hawking mass is monotonically non-decreasing if the spacetime satisfies the dominant energy condition. Thus, it is critical to understand the geometry of spherical symmetry. 

\begin{mydef}\label{def:spherically symmetric spacetime}
A spacetime $(N^{4}, \gbar)$ is said to be \emph{spherically symmetric} if its isometry group $\text{Isom}(N^{4})$ contains a subgroup $G$ that is isomorphic to the rotation group $SO(3)$; moreover for any point $p\in N$, the orbit of $p$ under the action of $G$ is a two-sphere with metric a multiple of the standard round metric. 
\end{mydef}
From the above definition, $N^{4}$ and standard sphere $S^{2}$ share $SO(3)$ as a subgroup in their isometry groups, thus $N^{4}$ share some of the symmetries with $S^{2}$, and hence the term \emph{spherically symmetric}. 

\begin{proposition}\label{prop:spherically symmetric metric that admit IMCVF coordinate chart and monotonicity}
If $(N^{4}, \gbar)$ is a spherically symmetric spacetime that admits a coordinate chart $\{t, r, \theta, \phi\}$, such that $\gbar$ takes the form:
\begin{equation}\label{equ:coordinate repren of spherically symmetric metric}
g = \bordermatrix{~ & t & r & \theta & \phi \cr
                  t & -v^{2}(t, r) & 0 & 0 & 0\cr
                  r & 0 & u^{2}(t, r) & 0 & 0\cr
                  \theta & 0 & 0 &  r^{2} & 0\cr
                  \phi & 0 & 0 & 0 & r^{2}\sin^{2}{\theta}\cr}
\end{equation}
where $u$ and $v$ are smooth functions of $t$ and $r$ only. Then within each $t = \text{constant}$ slice, smooth inverse mean curvature vector flow of coordinate sphere $S_{t, r}$ exists for all time with monotonically non-decreasing Hawking mass. 
\end{proposition}

\begin{remark}
Roughly speaking, all spherically symmetric spacetimes outside blackholes admit such metrics as in $(\ref{equ:coordinate repren of spherically symmetric metric})$.
\end{remark}

\begin{proof}
Let $g_{t, r}$ be the metric on $S_{t,r}$, then
\[g_{t, r} = r^{2}d\theta^{2} + r^{2}\sin^{2}{\theta}d\phi^{2},\]
and its inverse is given by:
\[g^{-1}_{t,r} = \frac{1}{r^{2}}d\theta^{2} + \frac{1}{r^{2}\sin^{2}{\theta}}d\phi^{2}.\]

Let $\vec{H}_{t,r}$ be the mean curvature vector of $S_{t, r}$. Notice that $\{\ppt, \ppr\}$ forms a frame for the normal bundle of $S_{t, r}$, and thus $H_{t,r}$ can be computed as follows:

\begin{align}
\vec{H}_{t, r} & = g_{t, r}^{ij}(\nabar_{\p_{i}}\p_{j})^{\text{nor}} \notag \\
& = g_{t, r}^{ij} \left(\frac{\<\nabar_{\p_{i}} \p_{j}, \p_{t}\>}{\<\p_{t}, \p_{t}\>} \p_{t} + \frac{\<\nabar_{\p_{i}} \p_{j}, \p_{r}\>}{\<\p_{r}, \p_{r}\>} \p_{r}\right)\notag \\
& = g_{t, r}^{\theta \theta}  \left(\frac{\<\nabar_{\p_{\theta}} \p_{\theta}, \p_{t}\>}{\<\p_{t}, \p_{t}\>} \p_{t} + \frac{\<\nabar_{\p_{\theta}} \p_{\theta}, \p_{r}\>}{\<\p_{r}, \p_{r}\>} \p_{r}\right) +  2 g_{t, r}^{\theta \phi}  \left(\frac{\<\nabar_{\p_{\theta}} \p_{\phi}, \p_{t}\>}{\<\p_{t}, \p_{t}\>} \p_{t} + \frac{\<\nabar{\p_{\theta}} \p_{\phi}, \p_{r}\>}{\<\p_{r}, \p_{r}\>} \p_{r}\right) \notag \\
& +  g_{t, r}^{\phi \phi}  \left(\frac{\<\nabar_{\p_{\phi}} \p_{\phi}, \p_{t}\>}{\<\p_{t}, \p_{t}\>} \p_{t} + \frac{\<\nabar_{\p_{\phi}} \p_{\phi}, \p_{r}\>}{\<\p_{r}, \p_{r}\>} \p_{r}\right) \notag \\
& = g_{t, r}^{\theta \theta}  (\Gammabar_{\theta \theta}^{t}\p_{t} + \Gammabar_{\theta \theta}^{r}\p_{r}) + 0 + g_{t, r}^{\phi \phi}  (\Gammabar_{\phi \phi}^{t}\p_{t} + \Gammabar_{\phi \phi}^{r}\p_{r}) \notag\tag{$g^{-1}_{t, r}$ is diagonal}\\
& = \frac{1}{r^{2}} \left(-\frac{r}{u^{2}}\right)\p_{r} + \frac{1}{r^{2}\sin^{2}{\theta}} \left(-\frac{r\sin^{2}{\theta}}{u^{2}}\right)\p_{r}\tag{See Section \ref{appen:curvatures in spherically symmetric spacetime}}\notag\\
& = -\frac{2}{r}\frac{1}{u^{2}} \p_{r}. 
\end{align}
That is 
\begin{equation}\label{equ:mean curvature of sphere in spherically symmetric spacetime}
\vec{H}_{t, r} = -\frac{2}{r}\frac{1}{u^{2}} \p_{r}.
\end{equation}

The inverse mean curvature vector is:
\begin{equation}\label{equ:inverse mean curvature vector in spherically symmetric spacetime}
\vec{I}_{t, r} = -\frac{\vec{H}_{t, r}}{\<\vec{H}_{t, r}, \vec{H}_{t, r}\>} = \frac{\frac{2}{r}\frac{1}{u^{2}} \p_{r}}{\<-\frac{2}{r}\frac{1}{u^{2}} \p_{r}, -\frac{2}{r}\frac{1}{u^{2}} \p_{r}\>} = \frac{\frac{2}{r}\frac{1}{u^{2}} \p_{r}}{\frac{4}{r^{2}}\frac{1}{u^{4}}u^{2}} = \frac{r}{2}\p_{r}.
\end{equation}

Therefore, inverse mean curvature vector flow of $S_{t, r}$ is a reparametrization of radial flow, and hence is smooth and exits for all time. 

To prove monotonicity of Hawking mass, recall that:
\begin{equation}\label{equ:Hawking mass of sphere in spherically symmetric spacetime}
m_{H}(S_{t, r}) =  \sqrt{\frac{|S_{t, r}|}{16 \pi}} \left(1 - \frac{1}{16\pi} \int_{S_{t, r}} \gbar(\vec{H}_{t, r}, \vec{H}_{t, r})\, dA_{t, r} \right),
\end{equation}
where $dA_{t, r}$ is the area form of $S_{t, r}$, and $|S_{t, r}|$ is the area of $S_{t, r}$. Note that 
\begin{align*}
|S_{t, r}| & = \int_{S_{t, r}}\, dA_{t, r} = \int_{0}^{2\pi}\int_{0}^{\pi} r^{2}\sin{\theta}\, d\theta d\phi\\
& = r^{2}\int_{0}^{2\pi}  2 \, d \phi \\
& = 4\pi r^{2}.
\end{align*}
Plug the mean curvature vector of $S_{t, r}$ (\ref{equ:mean curvature of sphere in spherically symmetric spacetime}) into the Hawking mass equation $(\ref{equ:Hawking mass of sphere in spherically symmetric spacetime})$, we get:
\begin{align}
m_{H}(S_{t,r}) & = \sqrt{\frac{4 \pi r^{2}}{16 \pi}} \left(1 - \frac{1}{16\pi} \int_{S_{t, r}} \gbar\left(-\frac{2}{r}\frac{1}{u^{2}} \p_{r}, -\frac{2}{r}\frac{1}{u^{2}} \p_{r}\right)\, dA_{t, r} \right) \notag \\
& = \frac{r}{2}\left(1 - \frac{1}{16\pi} \int_{S_{t, r}} \frac{4}{r^{2}} \frac{1}{u^{4}} \gbar(\p_{r}, \p_{r})\, dA_{t, r} \right) \notag \\
& = \frac{r}{2}\left(1 - \frac{1}{16\pi} \int_{0}^{2\pi}\int_{0}^{\pi} \frac{4}{r^{2}} \frac{1}{u^{2}}\, r^{2} sin{\theta} \, d\theta \, d\phi \right) \notag \\
& = \frac{r}{2}\left(1 - \frac{1}{4\pi}\frac{1}{u^{2}} 4 \pi \right) \notag \\
& = \frac{r}{2} \left(1 - \frac{1}{u^{2}}\right).
\end{align}

By Equation (\ref{equ:inverse mean curvature vector in spherically symmetric spacetime}), the variation of the Hawking mass of $S_{t, r}$ along inverse mean curvature vector flow is given by:
\begin{align*}\label{variation of Hawking mass of sphere in inverse mean curvature vector flow}
\vec{I}_{t, r} (m_{H}(S_{t, r})) & = \frac{r}{2}\frac{\p m_{H}(S_{t, r})}{\p_{r}} = \frac{r}{2}\left[\frac{1}{2}\left(1 - \frac{1}{u^{2}} \right) + r\frac{u_{, r}}{u^{3}}\right].
\end{align*}
Let $G$ be the Einstein curvature tensor of $(N^{4}, \gbar)$. $G(\p_{t}, \p_{t})$ can be computed as (see Equation $(\ref{Gtt})$ in Appendix \ref{appen:curvatures in spherically symmetric spacetime}):
\begin{equation}
G(\p_{t}, \p_{t}) = \frac{2}{r}\frac{u_{, r}}{u^{3}} v^{2} + \frac{1}{r^{2}}v^{2}\left(1 - \frac{1}{u^{2}} \right)  = \frac{2v^{2}}{r^{2}}\left[\frac{1}{2}\left(1 - \frac{1}{u^{2}} \right) + r\frac{u_{, r}}{u^{3}} \right].
\end{equation}
By the dominant energy condition (\ref{equ:dominant energy condition}), $G(\p_{t}, \p_{t}) \geq 0$. Therefore
\begin{equation}
\frac{1}{2}\left(1 - \frac{1}{u^{2}} \right) + r\frac{u_{, r}}{u^{3}} \geq 0,
\end{equation}
which implies that $\vec{I}_{t, r}(m_{H}(S_{t, r})) \geq 0$, as desired. 
\end{proof}

\begin{figure}[!htb]
\centering
\includegraphics[scale = .4]{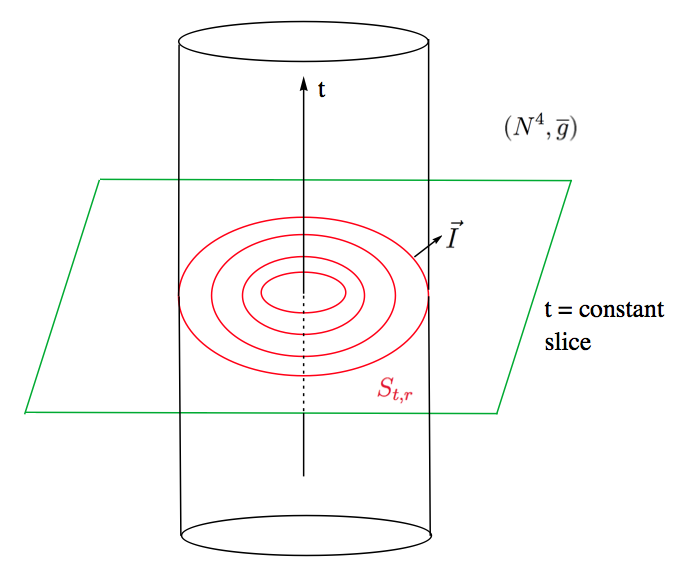}
\caption{Inverse mean curvature vector flow of coordinate spheres $S_{t, r}$ in spherically symmetric spacetime (\ref{equ:coordinate repren of spherically symmetric metric}).}
\label{fig:foliation of spherically symmetric spacetime with spheres in IMCVF}
\end{figure}

For any spherically symmetric spacetime in Proposition \ref{prop:spherically symmetric metric that admit IMCVF coordinate chart and monotonicity}, it can then be foliated by $t = \text{constant}$ spacelike hyperplanes, and each hyperplane can be foliated by smooth inverse mean curvature vector flow spheres (see Figure \ref{fig:foliation of spherically symmetric spacetime with spheres in IMCVF}). 

\chapter{Spacetimes that Admit Inverse Mean Curvature Vector Flow Solutions}
\label{chap:imcvf}

In this chapter we construct non-spherically symmetric, non-static spacetimes that admit smooth global solutions to inverse mean curvature vector flows. Surfaces we study here satisfy the topological and geometric setups defined in Section \ref{sec:codim two surface geometry}. In Section \ref{sec:motivation for imcvf}, we recall the definition of inverse mean curvature vector flow and the bad existence theory of such flows. Motivated by the spherically symmetry case,  in Section \ref{sec:construction of spacetime with imcvf coordinates}, we construct spacetimes that admit a special coordinate chart (called an inverse mean curvature vector flow coordinate chart) in which smooth inverse mean curvature vector flow of coordinate spheres exists for all time. We show that there in fact exist infinitely many spacetimes that admit such coordinate charts, and hence admit smooth solutions to inverse mean curvature vector flow. 

In Section \ref{sec:coordinate free analogue of imcvf and steering parameters}, we give a coordinate-free analogue of our construction, and show that we can actually ``steer'' a spacetime metric in a certain direction to make the mean curvature vector of a surface embedded in a spacelike hypersurface $M$ tangential to $M$. Finally in Section \ref{sec:generalizations of imcvf}, we discuss some generalizations of the technique we use in constructing inverse mean curvature vector flow coordinates. 

\section{Motivations from Spherically Symmetric Spacetimes and Main Results}
\label{sec:motivation for imcvf}
Given a spacetime $(N^{4}, \gbar)$ and a closed codimension two surface $\Sigma$ with induced metric $g_{\Sigma}$ and mean curvature vector $\vec{H}_{\Sigma}$, recall from Section \ref{subsec:mean curvature vector, hawking mass, normal variation} that a \emph{smooth inverse mean curvature vector flow} of $\Sigma$ is a normal variation 
\[F: \Sigma \times [0, T]\longrightarrow N, \quad (\Sigma, s)\mapsto F(\Sigma, s) =: \Sigma_{s},\] 
such that 
\begin{equation}\label{equ:IMCVF equation revisit}
\frac{\p}{\p s}F(x, s) = \vec{I}_{\Sigma_{s}}(x, s),\quad \forall\, (x, s)\in \Sigma_{s},
\end{equation}
where $\vec{I}_{\Sigma_{s}}$ is the inverse mean curvature vector of $\Sigma_{s}$ defined as
\[\vec{I}_{\Sigma_{s}}: = -\frac{\vec{H}_{\Sigma_{s}}}{\<\vec{H}_{\Sigma_{s}}, \vec{H}_{\Sigma_{s}}\>_{\gbar}}.\]
The inverse mean curvature vector flow equation \ref{equ:IMCVF equation revisit} (same as Equation \ref{equ:IMCVF equation}) is a forward-backward parabolic PDE, forward-parabolic in spacelike directions and backward-parabolic in timelike directions. Such a PDE lacks a general existence theory. However, such PDEs can still have solutions if we start with the ``right'' initial conditions. In spherically symmetric spacetimes, the ``right'' initial surfaces are spherically symmetric spheres. Had we chosen some other sphere to start with, inverse mean curvature vector flow is very likely to not exist. This is due to the mean curvature vector computation in Equation \ref{equ:mean curvature of sphere in spherically symmetric spacetime} in spherical symmetry: it is radial and has no components in the timelike direction, therefore the inverse mean curvature vector flow of spheres will be contained inside $t = \text{constant}$ spacelike slices. Geometrically, spherical symmetry eliminates all the timelike ``wiggles'' of the flow surfaces, hence restricting the flow direction to be spacelike. Since the inverse mean curvature vector flow equation is backward-parabolic only in timelike directions, inverse mean curvature vector flow exists without running into singularities in spherical symmetry (see Figure \ref{fig:foliation of spherically symmetric spacetime with spheres in IMCVF}).

How do we generalize the spherically symmetric case to non-symmetric spacetimes? Note that the mean curvature vector $\vec{H}_{\Sigma}$ of $\Sigma$ is a section of the normal bundle $N\Sigma$, and thus has a timelike component and a spacelike component. Motivated by the spherically symmetric case, intuitively if the flow surfaces all have ``purely \emph{spacelike}'' mean curvature vectors, we might hope to get a better existence theory. 

\begin{proposition}[\cite{BrayHaywardMarsSimon2007}]\label{prop:spacelike IMCVF equals tangential mean curvature in IMCF}
The family of closed embedded spacelike surfaces $\{\Sigma_{s}\}$ is a solution to the smooth inverse mean curvature vector flow with spacelike inverse mean curvature vector $\vec{I}_{\Sigma_{s}}$ everywhere on the surfaces if and only if there exists a spacelike hypersurface $M^{3}\subset N$, such that the mean curvature vector $\vec{H}_{\Sigma_{s}}$ is tangential to $M$ at all $(x, s)$, and $\{\Sigma_{s}\}$ is a solution to the smooth inverse mean curvature flow in $M$.
\end{proposition}

\begin{proof}
Given $\{\Sigma_{s}\}$ a solution to the smooth inverse mean curvature vector flow with \emph{spacelike} inverse mean curvature vector $\vec{I}_{\Sigma_{s}}$. Consider the hypersurface $M$ of $N$ defined by the union of all the surfaces $\Sigma_{s}$, i.e. the ``sweep out'' region of the the flow surfaces. Since the flow is smooth and spacelike, $M$ is a smooth manifold which is spacelike as well. $\vec{H}_{\Sigma_{s}}$ is tangential to $M$ by the construction of $M$. Since $\Sigma_{s}$ is of codimension one in $M$, $\vec{H}_{\Sigma_{s}}$ is parallel to the unit outward normal vector field $\nu_{\Sigma_{s}}$ along $\Sigma_{s}$, i.e.,
\begin{equation}
\vec{H}_{\Sigma_{s}} = -\lambda \nu_{\Sigma_{s}},
\end{equation}
at each $(x, s)\in \Sigma_{s}$ for some smooth positive function $\lambda$. $\lambda$ is chosen to positive since $\vec{H}_{\Sigma_{s}}$ points inward. Therefore
\begin{equation}\label{spacelike IMCVF implies IMCF}
\vec{I}_{\Sigma_{s}} = -\frac{\vec{H}_{\Sigma_{s}}}{\<\vec{H}_{\Sigma_{s}}, \vec{H}_{\Sigma_{s}}\>} = \frac{\lambda \nu_{\Sigma_{s}}}{\<-\lambda \nu_{\Sigma_{s}}, -\lambda \nu_{\Sigma_{s}}\>} = \frac{\lambda \nu_{\Sigma_{s}}}{\lambda^{2} \<\nu_{\Sigma_{s}}, \nu_{\Sigma_{s}}\>} = \frac{\nu_{\Sigma_{s}}}{\lambda} = \frac{\nu_{\Sigma_{s}}}{H_{\Sigma_{s}}},
\end{equation}
where $H_{\Sigma_{s}}$ is the mean curvature scalar of $\Sigma_{s}$ in the direction of $\nu_{\Sigma_{s}}$, defined by Equation \ref{def:mean curvature scalar}:
\begin{equation*}
H_{\Sigma_{s}} : = -\<\vec{H}_{\Sigma_{s}}, \nu_{\Sigma_{s}}\> = -\<-\lambda \nu_{\Sigma_{s}}, \nu_{\Sigma_{s}}\> = \lambda.
\end{equation*}
From equation (\ref{spacelike IMCVF implies IMCF}), we see that $\{\Sigma_{s}\}$ indeed is a solution to the smooth inverse mean curvature flow in $M$.
\\
Conversely, suppose $M$ is a spacelike hypersurface and $\{\Sigma_{s}\}$ is a solution to the smooth inverse mean curvature flow in $M$. Assuming that $\vec{H}_{\Sigma_{s}}$ is tangential to $M$ at each $(x, s)$, we know that $\vec{H}_{\Sigma_{s}}$ is spacelike as well. Moreover, by reversing the computations in equation ($\ref{spacelike IMCVF implies IMCF}$), we have 
\[\vec{I}_{\Sigma_{s}} = -\frac{\vec{H}_{\Sigma_{s}}}{\<\vec{H}_{\Sigma_{s}}, \vec{H}_{\Sigma_{s}}\>} .\]
Thus $\{\Sigma_{s}\}$ is a solution to the smooth inverse mean curvature vector flow equation in $N$, with spacelike inverse mean curvature vectors.
\end{proof}

Therefore spacelike smooth inverse mean curvature vector flow solutions in $N$ correspond to smooth inverse mean curvature flow solutions in a spacelike hypersurface of $N$ with tangential mean curvature vector fields. Huisken and Ilmanen defined a weak notion of inverse mean curvature flow in which jumps are allowed. This suggests the following definition of a weak solution of inverse mean curvature vector flow:

\begin{mydef}[Weak Solution to Inverse Mean Curvature Vector Flow, \cite{BrayHaywardMarsSimon2007}]\label{def:weak solution to IMCVF}
A family of spacelike surfaces $\{\Sigma_{s}\}$ is said to be a weak solution to the inverse mean curvature vector flow equation if there exists a spacelike hypersurface $M$ in $N$ containing $\Sigma_{s}$ such that $\vec{H}_{\Sigma_{s}}$ is tangential to $M$ everywhere for all $s\in [0, T]$, and $\{\Sigma_{s}\}$ is a solution to Huisken-Ilmanen inverse mean curvature flow in $M$ (i.e. with jumps).
\end{mydef}

Now we focus on the problem of finding spacetimes with inverse mean curvature vector flow solutions. 

Suppose a spacetime $(N^{4}, \gbar)$ admits the following special foliation: $N$ is foliated by spacelike hyperplanes, and then within each hyperplane, smooth inverse mean curvature vector flow of spheres exists and foliates the entire hyperplane. Consequently, the mean curvature vector of the flow spheres are tangential to the hyperplane. If $N$ admits such a special ``double'' foliation (e.g. spherically symmetric spacetimes), then $N$ has to be topologically equivalent to $(\RR^{3}\backslash B_{1}) \times \RR$, where $B_{1}$ is the closed unit ball in $\RR^{3}$.

\begin{figure}[!htb]
\centering
\includegraphics[scale = .5]{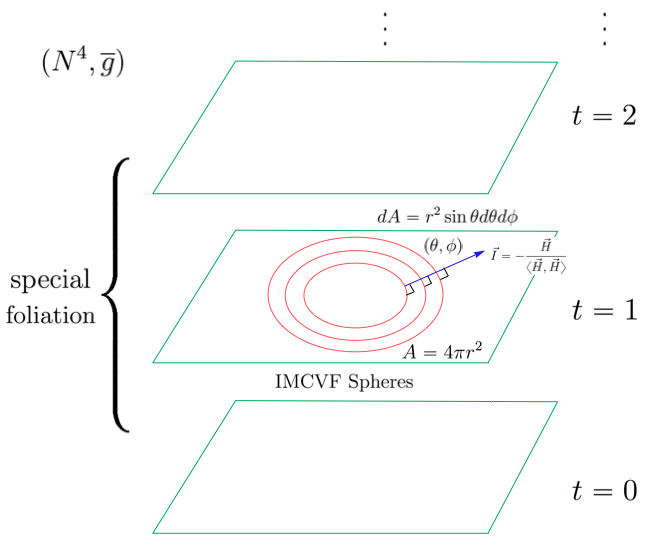}
\caption{Special foliation of a spacetime $(N^{4}, \gbar)$: first foliated by hyperplanes, then each hyperplane is foliated by inverse mean curvature vector flow of spheres. This generalizes the spherically symmetric case in Figure \ref{fig:foliation of spherically symmetric spacetime with spheres in IMCVF}.}
\label{fig:special foliation}
\end{figure}

Suppose $(N^{4}, \gbar)$ admits such a special foliation. We can use this to define coordinates that generalize the spherically symmetric coordinates. We define the $t$-coordinate by setting each hyperplane as $t = \text{constant}$, thus the $t$-coordinate tells us which hyperplane we are on. For each inverse mean curvature vector flow sphere, define $A = 4\pi r^{2}$, where $A$ is the area of that sphere. This defines a very natural $r$-coordinate. Since inverse mean curvature vector flow is area expanding, the $r$-coordinate is well-defined. For simplicity we assume that $r \geq 1$, i.e. the initial spheres on each hyperplane have area $4\pi$. Then define $(\theta, \phi)$-coordinates on an initial sphere, $0 < \theta < \pi$ and $0 < \phi < 2\pi$, such that the area form satisfies 
\begin{equation}\label{equ:area form condition}
dA_{0} = \frac{A(0)}{4\pi}\sin{\theta}d\theta d\phi = \sin{\theta}d\theta d\phi,
\end{equation} 
where $A(0)$ is the area of the initial sphere. Extend $(\theta, \phi)$ by setting them to be constant along perpendicular directions of the initial sphere, such that $(\theta, \phi)$ coordinates are defined for each sphere. By the extension, $\ppr$ will be perpendicular to each sphere. The equation for the area form (\ref{equ:area form condition}) will be preserved: $dA_{r} = \frac{A(r)}{4\pi}\sin{\theta}d\theta d\phi$\footnote{In smooth inverse mean curvature vector flow, $\frac{d}{ds}(dA_{s}) = -\<\vec{I_{s}}, \vec{H_{s}}\>dA_{s} = dA_{s}$, where $dA_{s}$ is the area form of $\Sigma_{s}$. The solution to this equation is $dA_{s} = e^{s}dA_{0}$. The area of $\Sigma_{s}$ is thus given by $A(s) = A(0)e^{s} = 4\pi r^{2}$, by the definition of the $r$-coordinate. Thus $e^{s}\frac{A(0)}{4\pi}= r^{2}$. Therefore the area form in the $r$ parameter is $dA_{r} = e^{s}dA_{0} = e^{s}\frac{A(0)}{4\pi}\sin{\theta}d\theta d\phi = r^{2}\sin{\theta}d\theta d\phi = \frac{A(r)}{4\pi}\sin{\theta}d\theta d\phi$.}.
See Figure \ref{fig:special foliation} for an illustration. 

Therefore we have proved the ``only-if'' direction of the following theorem:
\begin{theorem}
\label{thm:equivalence between existence of special foliation and existence of IMCVF coordinate chart}
A spacetime $(N^{4}, \gbar)$ is foliated by spacelike hyperplanes, and each hyperplane is foliated by smooth inverse mean curvature vector flow of spheres if and only if there exists a coordinate chart $\{t, r, \theta, \phi\}$ of $N$, such that in this coordinate chart the metric has the form:
\begin{equation}\label{equ:local representation of gbar in special foliation coordinates}
\overline{g} = \bordermatrix{~ & t & r & \theta & \phi \cr
                  t & -v^{2} & d & e & f\cr
                  r & d & u^{2} & 0 & 0\cr
                  \theta & e & 0 & a & c\cr
                  \phi & f & 0 & c & b\cr}
\end{equation}
where $u, v, a, b, c, d, e, f$ are smooth functions on $N$, 
and the following four conditions are satisfied:
\begin{align}
& (1) \quad  \<\ppr, \pptheta\> = 0;\label{def:first conditions of IMCVF coordinate chart}\\
& (2) \quad  \<\ppr, \ppphi\> = 0;\label{def:second conditions of IMCVF coordinate chart}\\
& (3) \quad  dA_{t, r} = r^{2}\sin{\theta}d\theta d\phi \,\,(\text{i.e.}\,\, ab - c^{2} = r^{4}\sin^{2}{\theta});\label{def:third conditions of IMCVF coordinate chart}\\
& (4) \quad  \vec{H}_{t, r}\,\,\text{is tangential to the $t = \text{constant}$ hyperplane};\label{def:fourth conditions of IMCVF coordinate chart}
\end{align}
where $dA_{t, r}$ and $\vec{H}_{t, r}$ are the area form and the mean curvature vector of the coordinate sphere $S_{t, r}$, respectively.
\end{theorem}

\begin{proof}[Proof of the ``if'' direction]
Given a coordinate chart $\{t, r, \theta, \phi\}$ of $(N, \gbar)$ such that the $\gbar$ satisfies the four conditions, $N$ is then foliated by $t = \text{constant}$ slices which are spacelike since the metric has the form $(\ref{equ:local representation of gbar in special foliation coordinates})$. For any $t = \text{constant}$ slice, the coordinate spheres $\{S_{t,r}\}$ are solutions of a normal flow since $\<\ppr, \pptheta\> = \<\ppr, \ppphi\> = 0$. We reparametrize the flow by setting $s : = C + 2\ln{r}$, where $C$ is a positive constant. Then 
\begin{equation}
\frac{d}{dr}(dA) = \frac{d}{ds}(dA)\frac{ds}{dr} = \frac{2}{r}\frac{d}{ds}(dA).
\end{equation}
On the other hand by condition (3)
\begin{equation}
\frac{d}{dr}(dA) = \frac{d}{dr}(r^{2}\sin{\theta}d\theta d\phi) = 2r\sin{\theta}d\theta d\phi = \frac{2}{r}dA.
\end{equation}
Combing the two equations above, we have $\frac{d}{ds}(dA) = dA$. Thus, by the first variation of area formula, $\{S_{t, r}\}$ when reparameterized by $r^{2} = Ce^{s}$, are smooth solutions to inverse mean curvature flow. By condition (4), the mean curvature vector of $S_{t, r}$ stays tangential to the slice, therefore $\{S_{t, r}\}$ with the above reparameterization are smooth solutions to inverse mean curvature vector flow. 
\end{proof}

\begin{mydef}[Inverse Mean Curvature Vector Flow Coordinate Chart]\label{def:IMCVF coordinates}
If a spacetime $(N^{4}, \gbar)$ admits a coordinate chart $\{t, r, \theta, \phi\}$ such that the four conditions (\ref{def:first conditions of IMCVF coordinate chart}), (\ref{def:second conditions of IMCVF coordinate chart}), (\ref{def:third conditions of IMCVF coordinate chart}) and (\ref{def:fourth conditions of IMCVF coordinate chart}) are satisfied, then $\{t, r, \theta, \phi\}$ is called an inverse mean curvature vector flow coordinate chart, and $N$ is called a spacetime that admits an inverse mean curvature vector flow coordinate chart. 
\end{mydef}

We sometimes refer the fourth condition (\ref{def:fourth conditions of IMCVF coordinate chart}) as the \emph{steering condition}, as it forces the coordinate spheres to stay inside the spacelike hyperplane during inverse mean curvature vector flow. 

Many spherically symmetric spacetimes admit an inverse mean curvature vector flow coordinate chart (e.g. coordinate chart (\ref{equ:coordinate repren of spherically symmetric metric}) with $d = e = f = c = 0$, and $a = r^{2}$, $b = r^{2}\sin^{2}{\theta}$ and radial mean curvature vector by Equation (\ref{equ:mean curvature of sphere in spherically symmetric spacetime})).  However, given an arbitrary spacetime $(N^{4}, \gbar)$, it is generally impossible to reparametrize it with an inverse mean curvature vector flow coordinate chart (e.g. a spacetime that is not topologically equivalent to ($\RR^{3}\backslash B_{1})\times \RR$). However, is it possible to \emph{construct} a spacetime that admits such a coordinate chart? In the next section we show that we can actually construct many such spacetimes:

\begin{proposition}[Existence of Spacetimes That Admit an Inverse Mean Curvature Vector Flow Coordinate Chart] 
\label{prop:existence of spacetimes that admit IMCVF coordinate chart}

Let $U:= (\RR^{3}\backslash B_{1})\times \RR\subset \RR^{4}$. There exist infinitely many spacetime metrics of the form of (\ref{equ:local representation of gbar in special foliation coordinates}) that admits an inverse mean curvature vector flow coordinate chart. 
\end{proposition}

Combining Proposition \ref{prop:existence of spacetimes that admit IMCVF coordinate chart} and Theorem \ref{thm:equivalence between existence of special foliation and existence of IMCVF coordinate chart}, we have the following main theorem of this thesis:

\begin{theorem}[\textbf{Main Theorem}]
\label{thm:existence of spacetimes that admit IMCVF coordinates}
There exist infinitely many non-spherically symmetric, non-static spacetimes that admit inverse mean curvature vector flow coordinate charts. Given such a spacetime $U$ with an inverse mean curvature vector flow coordinate chart $(t, r, \theta, \phi)$ and the constructed spacetime metric $\gbar$. The coordinate spheres $S_{t,r}$ contained in each $t = \text{constant}$ slice, when reparameterized by $r^{2} = e^{s}$, are smooth global solutions to the inverse mean curvature vector flow equation. 
\end{theorem}

\section{Construction of Spacetimes That Admit Inverse Mean Curvature Vector Flow Coordinate Charts}
\label{sec:construction of spacetime with imcvf coordinates}

In this section we prove Proposition \ref{prop:existence of spacetimes that admit IMCVF coordinate chart}. Let $U= (\RR^{3}\backslash B_{1})\times \RR$. It is easy to construct a spacetime metric $\gbar$ that admits a coordinate chart $\{t, r, \theta, \phi\}$ that satisfies condition (\ref{def:first conditions of IMCVF coordinate chart}) and (\ref{def:second conditions of IMCVF coordinate chart}). Simply define
\begin{equation}\label{equ:local representation of gbar in IMCVF coordinates step 1}
\overline{g}: = \bordermatrix{~ & t & r & \theta & \phi \cr
                  t & -v^{2} & d & e & f\cr
                  r & d & u^{2} & 0 & 0\cr
                  \theta & e & 0 & a & c\cr
                  \phi & f & 0 & c & b\cr}
\end{equation}
where $a,b,c,d,e,f,u,v$ are arbitrary smooth functions on $U$. Choosing two of the three variables $a, b$ and $c$ such that $ab - c^{2} = r^{4}\sin^{2}{\theta}$ satisfies condition (\ref{def:third conditions of IMCVF coordinate chart}).

The fourth condition requires $\vec{H}_{t, r}$, the mean curvature vector field of $S_{t, r}$, to be tangential to the $t = \text{constant}$ slice. This is equivalent to requiring $\vec{H}_{t, r}$ to be parallel to $\ppr$. We compute the conditions on the metric components such that this is true. 

\begin{lemma}\label{lemma:inverse and det of spacetime metric in IMCVF coordinates}
The determinant of the spacetime metric $\gbar$ in (\ref{equ:local representation of gbar in IMCVF coordinates step 1}) is given by:
\begin{align}
|\overline{g}| := \text{det}(\overline{g})&  = (-u^{2} v^{2} - d^{2}) (ab - c^{2}) + eu^{2} (cf - be) + fu^{2} (ce - af)\label{equ:det of spacetime metric in IMCVF coordinates}\\
& =  (-u^{2} v^{2} - d^{2}) (ab - c^{2}) + u^{2}(2cef - be^{2} - af^{2}) 
\end{align}
Moreover, the coordinate representation of the inverse $(\gbar^{-1})$ is given by:
\begin{equation}\label{equ:local representation of gbar inverse}
(\overline{g})^{-1} = \bordermatrix{~ & t & r & \theta & \phi \cr
                  t & \frac{u^{2}(ab - c^{2})}{|\gbar|} & \color{Red}\frac{-d(ab - c^{2})}{|\gbar|}& \color{Blue}\frac{u^{2}(cf - be)}{|\gbar|} & \color{Green}\frac{u^{2}(ce - af)}{|\gbar|} \cr
                  r & \color{Red}\frac{-d(ab - c^{2})}{|\gbar|} & \frac{-v^{2} (ab - c^{2}) + f(ce - af) + e (cf - be)}{|\gbar|} & \color{Orange}\frac{-d (cf - be)}{|\gbar|} & \color{Plum}\frac{-d(ce - af)}{|\gbar|} \cr
                  \theta & \color{Blue}\frac{u^{2} (cf - be)}{|\gbar|} & \color{Orange}\frac{-d (cf - be)}{|\gbar|} & \frac{-u^{2}v^{2} b - u^{2}f^{2} - bd^{2}}{|\gbar|} & \color{MidnightBlue}\frac{u^{2}v^{2}c + u^{2}ef + cd^{2}}{|\gbar|}\cr
                  \phi & \color{Green}\frac{u^{2}(ce - af)}{|\gbar|} & \color{Plum}\frac{-d (ce - af)}{|\gbar|} & \color{MidnightBlue} \frac{u^{2}v^{2}c + u^{2}ef + cd^{2}}{|\gbar|} & \frac{-u^{2}v^{2}a - u^{2}e^{2} - ad^{2}}{|\gbar|}\cr}
\end{equation}
\end{lemma}
\begin{proof}
See Section \ref{subsec:determinant and inverse of spacetime metric} in Appendix \ref{chap:geometric calculations}.
\end{proof}

\subsection{Geometry of $S_{t, r}$ and the Normal Bundle $NS_{t, r}$}
\label{subsec:geometry and normal bundle of coordinate sphere}
Fix a $t = \text{constant}$ slice $M$. Let $S_{t, r}$ be a coordinate sphere in $M$. We endow $S_{t,r}$ with the induced metric from $\gbar$, denoted as $g_{S}$. Then in the $\{\theta, \phi\}$ coordinate system, $g_{S}$ has the following representation:

\begin{equation}\label{equ:local representation of gs}
g_{S} := g|_{S_{t, r}} = \bordermatrix{~ &\theta & \phi \cr
                  \theta & a & c \cr
                  \phi & c & b \cr}
\end{equation}
Thus its inverse metric is:
\begin{equation}\label{equ:local representation of gs}
g^{-1}_{S} = \frac{1}{ab - c^{2}}\bordermatrix{~ &\theta & \phi \cr
                  \theta & b & -c \cr
                  \phi & -c & a \cr}
\end{equation}

The normal bundle $NS_{t, r}$ of $S_{t, r}$ is of rank-two. $\ppr$ is a nonzero section of $NS_{t, r}$, and thus can be used as a basis for the normal bundle. Let $\n$ be a complementary basis vector field of the normal bundle that is orthogonal to $\ppr$. Since $\ppr$ is outward spacelike, we can assume that $\n$ is future timelike. Therefore using the basis $\{\ppt, \ppr, \pptheta, \ppphi\}$ we can write $\n$ as
\[\n = \frac{\p}{\p t} + x\frac{\p}{\p r} + y\frac{\p}{\p \theta} + z\frac{\p}{\p \phi},\]
with $x,y, z$ yet to be determined, such that
\begin{itemize}
\item $\<\textbf{n}, \frac{\p}{\p r}\> = 0$;
\item $\<\textbf{n}, \frac{\p}{\p \theta}\> = 0$;
\item $\<\textbf{n}, \frac{\p}{\p \phi}\> = 0$.
\end{itemize}
We have three equations and three unknowns which give us

\begin{equation}\label{equ:normal vector field n}
\textbf{n} = \frac{\p}{\p t} + \frac{-d}{u^{2}}\frac{\p}{\p r} + \frac{cf - be}{ab - c^{2}}\frac{\p}{\p \theta} + \frac{ce - af}{ab - c^{2}}\frac{\p}{\p \phi}. 
\end{equation} 

\begin{lemma}\label{lemma:norm squred of n}
\begin{equation}\label{equ:norm squred of n}
\<\textbf{n}, \textbf{n}\> = \frac{\text{det}(\overline{g})}{u^{2}(ab - c^{2})} = \frac{\text{det}(\overline{g})}{u^{2}\text{det}(g_{S})} = :\frac{|\gbar|}{u^{2}|g_{S}|},
\end{equation}
where we set $|g_{S}|: = \text{det}(g_{S})$.
\end{lemma}
\begin{proof}
See Section \ref{subsec:inner product of n with n} in Appendix \ref{chap:geometric calculations}.
\end{proof}

\begin{remark}
Caution that since $\n$ is timelike, $\<\n, \n\> < 0$. Thus
\begin{equation}\label{equ:norm of n}
||\n||_{\gbar} = \left( -\<\n, \n\>\right)^{1/2} = \left(\frac{-|\gbar|}{u^{2}|g_{S}|} \right)^{1/2}.
\end{equation}
\end{remark}

Let $\{e_{r}, e_{n}\}$ be the normalized orthonormal frame obtained from $\{\ppr, \n\}$: 
\[e_{r} := \ppr/||\ppr||_{\gbar} = \frac{1}{u}\ppr,\quad e_{n} := \frac{\n}{||\n||_{\gbar}}.\]
Let
\[\vec{I}_{t, r} = -\frac{\vec{H}_{t, r}}{\<\vec{H}_{t, r}, \vec{H}_{t, r}\>}\]
be the outward-pointing inverse mean curvature vector. The geometry of the normal bundle of $S_{t,r}$ is given by Figure \ref{fig:normal bundle of coordinate sphere} below. Recall that $\vec{H}_{t, r}$ points inward by our convention.

\begin{figure}[!htb]
\centering
\includegraphics[scale = .5]{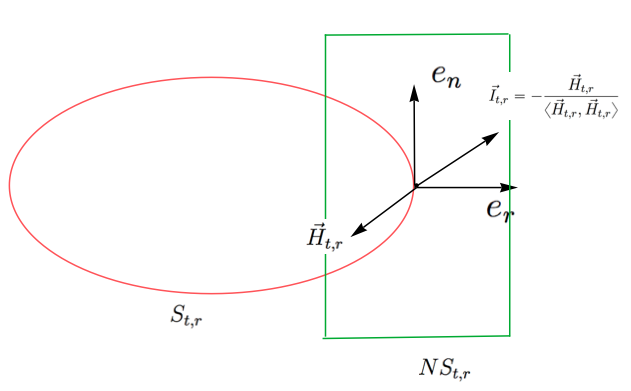}
\caption{Normal bundle of coordinate sphere $S_{t, r}$.}
\label{fig:normal bundle of coordinate sphere}
\end{figure}

\subsection{Mean Curvature Vector of $S_{t,r}$}
\label{subsec:mean curvature of coordinate sphere}
Now we have a coordinate sphere $(S_{t, r}, g_{S})$ with induced metric $g_{S}$ in $(U, \gbar, \nabar)$, where $\nabar$ is the Levi-Civita connection with respect to $\gbar$. We have an orthonormal frame $\{e_{r}, e_{n}\}$ of the normal bundle $NS_{t, r}$ defined in the previous subsection. In this subsection we compute the mean curvature vector $\vec{H}_{t, r}$ of $S_{t, r}$.
Recall the definition of $\vec{H}_{t, r}$ from Equation (\ref{def:mean curvature vector}):
\begin{align}
\vec{H}_{t, r} & = \trace_{g_{S}}\vec{\II}\tag{By definition of $\vec{H}_{t, r}$}\notag\\
& = g_{S}^{ij}\vec{\II}\left(\ppxi, \ppxj\right)\tag{$\ppxi, \ppxj\in \{\pptheta, \ppphi\}$}\notag\\
& = g_{S}^{ij}\left(\nabar_{\ppxi}\ppxj \right)\bigg|_{NS_{t, r}}\tag{By definition of $\vec{\II}(\ppxi, \ppxj)$}\notag\\
& = g_{S}^{ij}\left[\<\nabar_{\ppxi}\ppxj, e_{r}\>e_{r} - \<\nabar_{\ppxi}\ppxj, e_{n}\>e_{n} \right] \tag{$\{e_{r}, e_{n}\}$ is an orthonormal frame of $NS_{t, r}$}\notag\\
& = g_{S}^{ij}\<\nabar_{\ppxi}\ppxj, e_{r}\>e_{r} - g_{S}^{ij}\<\nabar_{\ppxi}\ppxj, e_{n}\>e_{n} \label{equ:mean curvature vector of IMCVF coordinate sphere step 1}
\end{align}

We compute the two parts in $\vec{H}_{t, r}$ separately. First of all
\begin{align}
g_{S}^{ij}\<\nabar_{\ppxi}\ppxj, e_{r}\> & = \frac{1}{u}g_{S}^{ij}\<\nabar_{\ppxi}\ppxj, \ppr\> = \frac{1}{u}g_{S}^{ij}\<\Gammabar^{k}_{ij}\ppxk, \ppr\>\notag\\
& = \frac{1}{u}g_{S}^{ij}(\Gammabar^{t}_{ij}\gbar_{tr} + \Gammabar^{r}_{ij}\gbar_{rr})\tag{$\ppr$ is perpendicular to $\pptheta, \ppphi$}\notag\\
& = \frac{1}{u}\gbar_{tr}(g_{S}^{\theta \theta}\Gammabar^{t}_{\theta \theta} + 2g_{S}^{\theta \phi}\Gammabar^{t}_{\theta \phi} + g_{S}^{\phi \phi}\Gammabar^{t}_{\phi \phi}) \notag \\
& + \frac{1}{u}\gbar_{rr}(g_{S}^{\theta \theta}\Gammabar^{r}_{\theta \theta} + 2g_{S}^{\theta \phi}\Gammabar^{r}_{\theta \phi} + g_{S}^{\phi \phi}\Gammabar^{r}_{\phi \phi})\notag\\
& = \frac{1}{u}\frac{1}{|g_{S}|}\big[\gbar_{tr}(b\Gammabar^{t}_{\theta \theta} -2c\Gammabar^{t}_{\theta \phi} + a\Gammabar^{t}_{\phi \phi}) + \gbar_{rr} (b\Gammabar^{r}_{\theta \theta} - 2c\Gammabar^{r}_{\theta \phi} + a\Gammabar^{r}_{\phi \phi})\big]\notag\\
& = \frac{1}{u}\frac{1}{|g_{S}|}\big[\gbar_{tr} (\star) + \gbar_{rr} (\star \star)\big] \label{equ:first term in mean curvature vector of IMCVF coordinate sphere step 1 step 1}
\end{align}

where
\begin{equation}\label{def:star term}
(\star): = b\Gammabar^{t}_{\theta \theta} -2c\Gammabar^{t}_{\theta \phi} + a\Gammabar^{t}_{\phi \phi},
\end{equation}
and 
\begin{equation}\label{def:double star term}
(\star \star): = b\Gammabar^{r}_{\theta \theta} - 2c\Gammabar^{r}_{\theta \phi} + a\Gammabar^{r}_{\phi \phi}.
\end{equation}

We now compute all the related Christoffel symbols.

\begin{align}
& \Gammabar_{\theta \theta}^{t} = \frac{1}{2} \Big( \overline{g}^{tt}(2 \overline{g}_{\theta t, \theta} - \overline{g}_{\theta \theta, t}) + \overline{g}^{tr} (2\overline{g}_{\theta r, \theta} - \overline{g}_{\theta \theta, r}) + \overline{g}^{t\theta}(2\overline{g}_{\theta \theta, \theta} - \overline{g}_{\theta \theta, \theta}) + \overline{g}^{t \phi} (2 \overline{g}_{\theta \phi, \theta}\notag\\
& - \overline{g}_{\theta \theta, \phi}) \Big) \notag \\
& = \frac{1}{2}\Big( \overline{g}^{tt}(2 \overline{g}_{\theta t, \theta} - \overline{g}_{\theta \theta, t}) - \overline{g}^{tr}  \overline{g}_{\theta \theta, r} + \overline{g}^{t\theta}\overline{g}_{\theta \theta, \theta}  + \overline{g}^{t \phi} (2 \overline{g}_{\theta \phi, \theta} - \overline{g}_{\theta \theta, \phi}) \Big) \notag \\
& = \frac{1}{2|\overline{g}|} \Big(u^{2}(ab - c^{2}) (2e_{,\theta} - a_{, t}) + d(ab - c^{2}) a_{, r} + u^{2}(cf - be)a_{, \theta} + u^{2}(ce - af)(2c_{, \theta} \notag\\
& - a_{,\phi}) \Big)\label{equ:gammarbar t theta theta}
\end{align}

\begin{align}
& \Gammabar_{\theta \phi}^{t} = \frac{1}{2}\Big(\overline{g}^{tt}(\overline{g}_{\phi t, \theta} + \overline{g}_{\theta t, \phi} - \overline{g}_{\theta \phi, t})  + \overline{g}^{tr}(g_{\phi r, \theta} + \overline{g}_{\theta r, \phi} - \overline{g}_{\theta \phi, r}) + \overline{g}^{t \theta}(\overline{g}_{\phi \theta, \theta} + \overline{g}_{\theta \theta, \phi} - \overline{g}_{\theta \phi, \theta})\notag \\ 
& + \overline{g}^{t\phi}(\overline{g}_{\phi \phi, \theta} + \overline{g}_{\theta \phi, \phi} - \overline{g}_{\theta \phi, \phi})\Big)\notag \\
& = \frac{1}{2}\Big(\overline{g}^{tt}(\overline{g}_{\phi t, \theta} + \overline{g}_{\theta t, \phi} - \overline{g}_{\theta \phi, t}) - \overline{g}^{tr} \overline{g}_{\theta \phi, r} + \overline{g}^{t\theta} \overline{g}_{\theta \theta, \phi} + \overline{g}^{t \phi}\overline{g}_{\phi \phi, \theta}\Big) \notag \\
& = \frac{1}{2|\overline{g}|}\Big(u^{2}(ab - c^{2}) (f_{, \theta} + e_{, \phi} - c_{,t})  + d(ab - c^{2})c_{, r} +  u^{2}(cf - be)a_{,\phi} + u^{2}(ce - af)b_{, \theta}\Big)\label{equ:gammabar t theta phi}
\end{align}

\begin{align}
& \Gammabar_{\phi \phi}^{t} = \frac{1}{2}\Big(\overline{g}^{tt}(2\overline{g}_{\phi t, \phi} - \overline{g}_{\phi \phi, t}) + \overline{g}^{tr}(2\overline{g}_{\phi r, \phi} - \overline{g}_{\phi \phi ,r}) + \overline{g}^{t\theta}(2\overline{g}_{\phi \theta, \phi} - \overline{g}_{\phi \phi, \theta}) + \overline{g}^{t\phi}(2\overline{g}_{\phi \phi, \phi}\notag\\
& - \overline{g}_{\phi \phi ,\phi}) \Big) \notag \\
& = \frac{1}{2} \Big(g^{tt}(2\overline{g}_{\phi t, \phi} - \overline{g}_{\phi \phi, t}) - \overline{g}^{tr}\overline{g}_{\phi \phi, r} +  \overline{g}^{t\theta}(2\overline{g}_{\phi \theta, \phi} - \overline{g}_{\phi \phi, \theta}) + \overline{g}^{t\phi}\overline{g}_{\phi \phi, \phi}\Big) \notag \\
& = \frac{1}{2|\overline{g}|}\Big(u^{2}(ab - c^{2})(2f_{,\phi} - b_{,t}) + d(ab - c^{2})b_{, r} + u^{2}(cf - be) (2c_{, \phi} - b_{, \theta}) \notag\\
& + u^{2}(ce - af)b_{,\phi} \Big)\label{equ:gammabar t phi phi}
\end{align}

\begin{align}
& \Gammabar^{r}_{\theta \theta} = \frac{1}{2}\Big( \gbar^{rt}(2\gbar_{\theta t, \theta} - \gbar_{\theta \theta, t}) + \gbar^{rr}(2\gbar_{\theta r, \theta} - g_{\theta \theta ,r}) + \gbar^{r\theta}(2\gbar_{\theta \theta ,\theta} - \gbar_{\theta \theta, \theta}) + \gbar^{r\phi}(2\gbar_{\theta \phi, \theta} \notag\\
& - \gbar_{\theta \theta, \phi})\Big)\notag\\
& = \frac{1}{2}\Big(\gbar^{rt}(2\gbar_{\theta t, \theta} - \gbar_{\theta \theta, t}) - \gbar^{rr}\gbar_{\theta \theta, r} + \gbar^{r\theta}\gbar_{\theta \theta,\theta} + \gbar^{r\phi}(2\gbar_{\theta \phi, \theta} - \gbar_{\theta \theta, \phi}) \Big)\notag\\
& = \frac{1}{2|\gbar|}\Big(-d(ab - c^{2})(2e_{,\theta} - a_{,t}) - [-v^{2}(ab - c^{2}) + f(ce - af) + e(cf - be)] a_{,r}\notag\\
& + (-d)(cf - be)a_{,\theta} + (-d)(ce - af)(2c_{,\theta} - a_{,\phi}) \Big)\label{equ:gammbar r theta theta}
\end{align}

\begin{align}
& \Gammabar^{r}_{\theta \phi} = \frac{1}{2}\Big(\gbar^{rt}(\gbar_{\phi t, \theta} + \gbar_{\theta t, \phi} - \gbar_{\theta \phi, t}) + \gbar^{rr}(\gbar_{\phi r, \theta}+ \gbar_{\theta r, \phi} - \gbar_{\theta \phi, r}) + \gbar^{r\theta}(\gbar_{\phi \theta, \theta}+ \gbar_{\theta \theta, \phi} - \gbar_{\theta \phi, \theta})\notag\\
& + \gbar^{r\phi}(\gbar_{\phi \phi, \theta} + \gbar_{\theta \phi, \phi} - \gbar_{\theta \phi, \phi})\Big)\notag\\
& = \frac{1}{2}\Big(\gbar^{rt}(\gbar_{\phi t, \theta} + \gbar_{\theta t, \phi}- \gbar_{\theta \phi, t}) - \gbar^{rr}\gbar_{\theta \phi, r} + \gbar^{r\theta}\gbar_{\theta \theta, \phi} + \gbar^{r\phi}\gbar_{\phi \phi, \theta}\Big)\notag\\
& = \frac{1}{2|\gbar|}\Big(-d(ab - c^{2})(f_{,\theta} + e_{,\phi} - c_{,t}) - [-v^{2}(ab - c^{2}) + f(ce - af) + e(cf - be)]c_{,r}\notag\\
& + (-d)(cf - be)a_{,\phi} + (-d)(ce - af)b_{,\theta}\Big)\label{equ:gammbar r theta phi}
\end{align}

\begin{align}
& \Gammabar^{r}_{\phi \phi} = \frac{1}{2}\Big(\gbar^{rt}(2\gbar_{\phi t,\phi} - \gbar_{\phi \phi, t}) + \gbar^{rr} (2\gbar_{\phi r,\phi} - \gbar_{\phi\phi, r}) + \gbar^{r\theta}(2\gbar_{\phi \theta, \phi} - \gbar_{\phi \phi,\theta}) + \gbar^{r\phi}(2\gbar_{\phi \phi, \phi} \notag\\
& - \gbar_{\phi \phi,\phi})\Big)\notag\\
& = \frac{1}{2}\Big(\gbar^{rt}(2\gbar_{\phi t,\phi} - \gbar_{\phi \phi, t}) - \gbar^{rr}\gbar_{\phi \phi, r} + \gbar^{r\theta}(2\gbar_{\phi \theta, \phi} - \gbar_{\phi \phi,\theta}) + \gbar^{r\phi}\gbar_{\phi \phi, \phi}\Big)\notag\\
& = \frac{1}{2|\gbar|}\Big(-d(ab - c^{2})(2f_{,\phi} - b_{,t}) - [-v^{2}(ab - c^{2}) + f(ce - af) + e(cf - be)]b_{,r} \notag\\
& + (-d)(cf - be)(2c_{,\phi }- b_{,\theta}) +  (-d)(ce - af)b_{,\phi} \Big)\label{equ:gammbar r phi phi}
\end{align}

With all the Christoffel symbols computed, we have:
\begin{align}
(\star) & = b\Gammabar^{t}_{\theta \theta} -2c\Gammabar^{t}_{\theta \phi} + a\Gammabar^{t}_{\phi \phi} \notag\\
& = \frac{1}{2|\gbar|}\Big\{u^{2}(ab - c^{2})[2be_{,\theta} - \cancel{a_{,t}b} - 2ce_{,\phi} - 2cf_{,\theta} + \cancel{2cc_{,t}} + 2af_{,\phi} - \cancel{ab_{,t}}]\notag\\
& + d(ab - c^{2})[a_{,r}b - 2cc_{,r} + ab_{,r}] + u^{2}(cf - be)[a_{,\theta}b - 2a_{,\phi}c + 2ac_{,\phi} - ab_{,\theta}]\notag\\
& + u^{2}(ce - af)[2bc_{,\theta} - a_{,\phi}b - 2b_{,\theta}c + ab_{,\phi}]\Big\} \label{equ:star term explicit}
\end{align}
The terms cancel in the second line above since $ab - c^{2} = r^{4}\sin^{2}{\theta}$, and thus is not a function of $t$.

Next we compute
\begin{align}
(\star \star) & = b\Gammabar^{r}_{\theta \theta} - 2c\Gammabar^{r}_{\theta \phi} + a\Gammabar^{r}_{\phi \phi} \notag\\
& = \frac{1}{2|\gbar|}\Big\{-d(ab - c^{2})[2be_{,\theta} - \cancel{a_{,t}b} - 2ce_{,\phi} - 2cf_{,\theta} + \cancel{2cc_{,t}} + 2af_{,\phi} - \cancel{ab_{,t}}]\notag\\
& - [-v^{2}(ab - c^{2}) + f(ce - af) + e(cf - be)][a_{,r}b - 2cc_{,r} + ab_{,r}]\notag\\
& -d(cf - be)[a_{,\theta}b - 2a_{,\phi}c + 2ac_{,\phi} - ab_{,\theta}] - d(ce - af)[2bc_{,\theta} - a_{,\phi}b - 2b_{,\theta}c + ab_{,\phi}]\Big\} \label{equ:double star term}
\end{align}

Now plug (\ref{equ:star term explicit}) and (\ref{equ:double star term}) back into (\ref{equ:first term in mean curvature vector of IMCVF coordinate sphere step 1 step 1}):
\begin{align}
& g_{S}^{ij}\<\nabar_{\ppxi}\ppxj, e_{r}\> = \frac{1}{u}\frac{1}{|g_{S}|} \big[ d(\star) + u^{2}(\star \star) \big]\notag\\
& = \frac{1}{u}\frac{1}{|g_{S}|}\frac{1}{2|\gbar|}\Big\{[du^{2}(ab - c^{2}) - u^{2}d(ab - c^{2})][2be_{,\theta} - 2ce_{,\phi} - 2cf_{,\theta} + 2af_{,\phi}]\notag\\
& + [d^{2}(ab - c^{2}) + u^{2}v^{2}(ab - c^{2}) - fu^{2}(ce - af) - eu^{2}(cf - be)][a_{,r}b - 2cc_{,r} + ab_{,r}]\notag\\
& + [du^{2}(cf - be) - u^{2}d(cf - be)][a_{,\theta}b - 2a_{,\phi}c + 2ac_{,\phi} - ab_{,\theta}]\notag\\
& + [du^{2}(ce - af) - u^{2}d(ce - af)][2bc_{,\theta} - a_{,\phi}b - 2b_{,\theta}c + ab_{,\phi}]\Big\}\notag\\
& = \frac{1}{u}\frac{1}{|g_{S}|}\frac{1}{2|\gbar|}[d^{2}(ab - c^{2}) + u^{2}v^{2}(ab - c^{2}) - fu^{2}(ce - af) - eu^{2}(cf - be)][a_{,r}b - 2cc_{,r} \notag\\
& + ab_{,r}] \notag\\
& = \frac{1}{u}\frac{1}{|g_{S}|}\frac{1}{2|\gbar|}(-|\gbar|)(ab - c^{2})_{,r}\tag{By Equation \ref{equ:det of spacetime metric in IMCVF coordinates}}\notag\\
& = -\frac{1}{u}\frac{4r^{3}\cancel{\sin^{2}{\theta}}}{r^{4}\cancel{\sin^{2}{\theta}}}\frac{1}{2} = -\frac{2}{r}\frac{1}{u}. \label{equ:first term in mean curvature vector of IMCVF coordinate sphere step 1 final}
\end{align}

Next we compute the second term in (\ref{equ:mean curvature vector of IMCVF coordinate sphere step 1}):
\begin{align}
g_{S}^{ij}\<\nabar_{\ppxi}\ppxj, e_{n}\> & = \frac{1}{||\n||_{\gbar}}g_{S}^{ij}\<\nabar_{\ppxi}\ppxj, \n\> = \frac{1}{||\n||_{\gbar}}g_{S}^{ij}\<\Gammabar^{t}_{ij}\ppt, \n\>\tag{$\n \perp \ppr, \pptheta, \ppphi$}\notag\\
& = \frac{1}{||\n||_{\gbar}}g_{S}^{ij}\Gammabar^{t}_{ij}\<\n, \n\> \tag{$\n - \ppt$ is spacelike}\notag\\
& = - g_{S}^{ij}\Gammabar^{t}_{ij}||\n||_{\gbar} = - ||\n||_{\gbar}(g_{S}^{\theta \theta}\Gammabar^{t}_{\theta \theta} + 2g_{S}^{\theta \phi}\Gammabar^{t}_{\theta \phi} + g_{S}^{\phi \phi}\Gammabar^{t}_{\phi \phi}) \notag\\
& = - \frac{||\n||_{\gbar}}{|g_{S}|}(\star) \label{equ:second term in mean curvature vector of IMCVF coordinate sphere step 1 step 1}
\end{align}
where $(\star)$ is computed in (\ref{equ:star term explicit}). Now plug (\ref{equ:norm squred of n}) into the above, we get:
\begin{equation}\label{equ:second term in mean curvature vector of IMCVF coordinate sphere step 1 final}
g_{S}^{ij}\<\nabar_{\ppxi}\ppxj, e_{n}\> = -\left(\frac{-|\gbar|}{u^{2}|g_{S}|}\right)^{1/2}\frac{1}{|g_{S}|}(\star) = -\frac{1}{u}\frac{(-|\gbar|)^{1/2}}{|g_{S}|^{3/2}}(\star).
\end{equation}

Combing (\ref{equ:first term in mean curvature vector of IMCVF coordinate sphere step 1 final}) and (\ref{equ:second term in mean curvature vector of IMCVF coordinate sphere step 1 final}), we have:

\begin{proposition}\label{prop:mean curvature vector of IMCVF coordinate sphere}
The mean curvature vector $\vec{H}_{t, r}$ of coordinate sphere $S_{t, r}$ in the spacetime metric $\gbar$ (\ref{equ:local representation of gbar in IMCVF coordinates step 1}) is given by:
\begin{equation}\label{equ:mean curvature vector of IMCVF coordinate sphere}
\vec{H}_{t, r} = -\frac{2}{r}\frac{1}{u}e_{r} + \frac{1}{u}\frac{(-|\gbar|)^{1/2}}{|g_{S}|^{3/2}}(\star)e_{n}.
\end{equation}
where $(\star)$ is given by (\ref{equ:star term explicit}). 
\end{proposition}

\begin{coro}\label{coro:radial mean curvature vector condition in IMCVF coordinates}
The mean curvature vector $\vec{H}_{t, r}$ of coordinate sphere $S_{t, r}$ is parallel to $\ppr$ everywhere if and only if $(\star) = 0$, where $(\star)$ is given by (\ref{equ:star term explicit}). Moreover, in this case, the mean curvature vector equals:

\begin{equation}\label{equ:radial mean curvature vector of IMCVF coordinate sphere}
\vec{H}_{t, r} = -\frac{2}{r}\frac{1}{u}e_{r}.
\end{equation}
\end{coro}

Notice the similarity between the mean curvature vector expression above and the mean curvature vector in the spherically symmetric case in (\ref{equ:mean curvature of sphere in spherically symmetric spacetime}). 

\begin{proof}
This is quite straightforward since $\vec{H}_{t, r}$ is parallel to $\ppr$ everywhere if and only if $\<\vec{H}_{t, r}, e_{n}\> = 0$. Expand this out we get:
\[0 = \<\vec{H}_{t, r}, e_{n}\> = \<\frac{1}{u}\frac{(-|\gbar|)^{1/2}}{|g_{S}|^{3/2}}(\star)e_{n}, e_{n}\> = -\frac{1}{u}\frac{(-|\gbar|)^{1/2}}{|g_{S}|^{3/2}}(\star).\]
Since $u \neq 0$, the above is equivalent to 
\[(\star) = 0,\] as desired.
\end{proof}

Therefore the fourth condition in the definition of inverse mean curvature vector flow coordinate chart is equivalent to $(\star) = 0$, that is:
\begin{align}
0 & = \frac{1}{2|\gbar|}\Big\{u^{2}(ab - c^{2})[2be_{,\theta} - 2ce_{,\phi} - 2cf_{,\theta} + 2af_{,\phi}]\notag\\
& + d(ab - c^{2})(ab - c^{2})_{,r} + u^{2}(cf - be)[a_{,\theta}b - 2a_{,\phi}c + 2ac_{,\phi} - ab_{,\theta}]\notag\\
& + u^{2}(ce - af)[2bc_{,\theta} - a_{,\phi}b - 2b_{,\theta}c + ab_{,\phi}]\Big\} \notag 
\end{align}
Or equivalently, using $ab - c^{2} = r^{4}\sin^{2}{\theta}$,
\begin{align}
0 & = [2be_{,\theta} - 2ce_{,\phi} - 2cf_{,\theta} + 2af_{,\phi}] + d\frac{4r^{3}\sin^{2}{\theta}}{u^{2}}\notag\\
& + \frac{cf - be}{r^{4}\sin^{2}{\theta}}[a_{,\theta}b - 2a_{,\phi}c + 2ac_{,\phi} - ab_{,\theta}] + \frac{ce - af}{r^{4}\sin^{2}{\theta}}[2bc_{,\theta} - a_{,\phi}b - 2b_{,\theta}c + ab_{,\phi}].
\label{equ:key equation form 1}
\end{align}

Here is the upshot: Equation (\ref{equ:key equation form 1}) is \emph{zeroth order} in the metric component $d$, thus we can choose two of the three variables $a, b, c$, and $e, f, u, v$ all together $6$ variables, and solve for $d$ explicitly:
\begin{align}
d & = -\frac{u^{2}}{4r^{3}\sin^{2}{\theta}}\Big\{[2be_{,\theta} - 2ce_{,\phi} - 2cf_{,\theta} + 2af_{,\phi}] + \frac{cf - be}{r^{4}\sin^{2}{\theta}}[a_{,\theta}b - 2a_{,\phi}c + 2ac_{,\phi} - ab_{,\theta}] \notag\\
& + \frac{ce - af}{r^{4}\sin^{2}{\theta}}[2bc_{,\theta} - a_{,\phi}b - 2b_{,\theta}c + ab_{,\phi}]\Big\}.\label{equ:key equation form 2}
\end{align}
Thus all four conditions (\ref{def:first conditions of IMCVF coordinate chart}), (\ref{def:second conditions of IMCVF coordinate chart}) and (\ref{def:third conditions of IMCVF coordinate chart}) and (\ref{def:fourth conditions of IMCVF coordinate chart}) are solvable with infinitely many solutions. Combing the above, we have proved Proposition \ref{prop:existence of spacetimes that admit IMCVF coordinate chart} and hence the main theorem \ref{thm:existence of spacetimes that admit IMCVF coordinates}.

There are six degrees of freedom in constructing our spacetime metric $\gbar$ that admits an inverse mean curvature vector flow coordinate chart. Moreover, the six free variables do not need to be spherically symmetric. The spherically symmetric (\ref{equ:coordinate repren of spherically symmetric metric}) metric is a special case of this large set of spacetime metrics. It is unknown if perturbations of spherically symmetric spacetime with inverse mean curvature vector flow coordinate chart still have inverse mean curvature vector flow solutions. More specifically, we conjecture that:
\begin{conjecture}
\label{conj:perturbation of minkowski spacetime admits IMCVF coordinate chart}
Given Minkowski space with inverse mean curvature vector flow coordinate chart that can be smoothly extended to the boundary, consider a perturbation of the spacetime metric. The resulting spacetime still admits inverse mean curvature vector flow solutions (in a single spacelike hypersurface) that exist for all time. 
\end{conjecture}

Notice that in Equation (\ref{equ:key equation form 2}) if all the variables are smooth then $d$ will be smooth except possibly when $\sin{\theta} = 0$, since $d$ is not defined by our formula there (see (\ref{equ:key equation form 2})). This happens at the north ($\theta = 0$) and south pole ($\theta = \pi$), which are two \emph{coordinate chart singularities}, not metric singularities of the spacetime. 

If $c, e, f$ are chosen to be $0$ at a neighborhood of the north and the south pole, then the right hand side of (\ref{equ:key equation form 2}) will be zero there. In this way $d$ can be extended smoothly across the two coordinate chart singularities, and will be smooth on the entire spacetime. 

Another way to extend $d$ smoothly over the coordinate chart singularities is to choose the metric to be spherically symmetric in a neighborhood of the north and south pole:

\begin{equation}\label{equ:coordinate repren of spherically symmetric metric revisi}
g = \bordermatrix{~ & t & r & \theta & \phi \cr
                  t & -v^{2}(t, r) & 0 & 0 & 0\cr
                  r & 0 & u^{2}(t, r) & 0 & 0\cr
                  \theta & 0 & 0 &  r^{2} & 0\cr
                  \phi & 0 & 0 & 0 & r^{2}\sin^{2}{\theta}\cr}
\end{equation}

This has the advantage that $d = 0$ around the coordinate chart singularities. Then extend $a, b, c, e, f$ smoothly to the entire spacetime while maintaining the condition that $ab - c^{2} = r^{4}\sin^{2}{\theta}$. The resulting metric still satisfies the four conditions since smooth inverse mean curvature vector flow of spheres exists for all time in spherically symmetry. $d$ will be smooth since $d = 0$ identically. 

One could study more general asymptotic conditions for (\ref{equ:key equation form 2}) to be smooth and bounded as $\theta$ approaches $0$ or $\pi$, but we choose not to discuss it further here. 

\begin{remark}
We can actually prove that $\vec{H}_{t, r}$ takes the form (\ref{equ:radial mean curvature vector of IMCVF coordinate sphere}) in inverse mean curvature vector flow coordinates without computing it out explicitly. We now show that as a sanity check of our computation. In inverse mean curvature vector flow coordinates $\vec{H}_{t, r}$ is parallel to $\ppr$, thus let $\lambda = \lambda(t, r, \theta, \phi)$ such that
\[\vec{H}_{t, r} = \lambda\ppr.\]
The inverse mean curvature vector is now
\[\vec{I}_{t, r} =  -\frac{\vec{H}_{t, r}}{\<\vec{H}_{t, r}, \vec{H}_{t, r}\>} = -\frac{\lambda\ppr}{\<\lambda\ppr, \lambda\ppr\>} = -\frac{1}{\lambda u^{2}}\ppr.\]
By the first variation of area formula (\ref{equ:first variation of area}), the rate of change of the area form of $S_{t, r}$ under \emph{outward radial flow} is given by
\begin{align}
\frac{d}{dr}dA_{S_{t, r}} & = -\<\vec{H}_{t, r}, \ppr\>dA_{S_{t, r}} \notag\\
& = -\<\vec{H}_{t, r}, \ppr\> r^{2}\sin{\theta}d\theta d\phi \notag\\
& = -\lambda u^{2}r^{2}\sin{\theta}d\theta d\phi
\end{align}
Notice the left hand side of the above equals to
\[\frac{d}{dr}dA_{S_{t, r}} = \frac{d}{dr}(r^{2}\sin{\theta}d\theta d\phi) = 2r\sin{\theta}d\theta d\phi.\]
Thus matching the two sides we get:
\[2r\sin{\theta}d\theta d\phi =  -\lambda u^{2}r^{2}\sin{\theta}d\theta d\phi,\]
that is $\lambda = -\frac{2}{r}\frac{1}{u^{2}}$. Therefore
\[\vec{H}_{t, r} = \lambda \ppr = -\frac{2}{r}\frac{1}{u^{2}}\ppr = -\frac{2}{r}\frac{1}{u}e_{r},\] which is the same as (\ref{equ:mean curvature vector of IMCVF coordinate sphere}).
\end{remark}

\section{Coordinate Free Analogue and Steering Parameters}
\label{sec:coordinate free analogue of imcvf and steering parameters}

In the previous section we have shown that there exist many spacetimes that admit inverse mean curvature vector flow coordinate chart, in which the coordinates spheres are solutions to the inverse mean curvature vector flow equation. The fourth condition (\ref{def:fourth conditions of IMCVF coordinate chart}) in the definition of inverse mean curvature vector flow coordinates can be viewed as a \emph{steering condition} that keeps the flow direction of coordinate spheres tangential to a spacelike hypersurface. 

\begin{figure}[!htb]
\centering
\includegraphics[scale = .45]{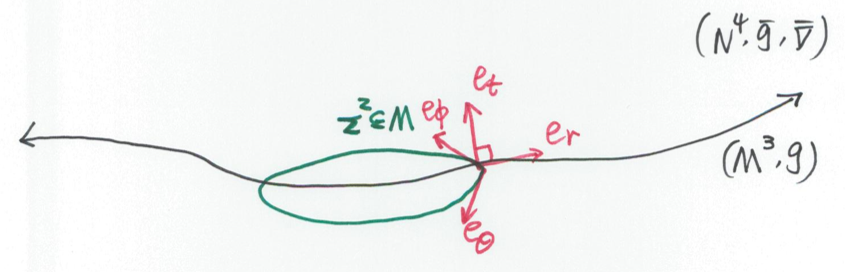}
\caption{Setup of inverse mean curvature vector flow steering.}
\label{fig:IMCVF steering setup}
\end{figure}

Given a spacetime $(N^{4}, \gbar, \nabar)$, a spacelike hypersurface $(M^{3}, g)$ with induced metric $g$, and a closed embedded surface $(\Sigma, g_{\Sigma})$ in $M$ with induced metric $g_{S}$. Assuming the normal bundle of $\Sigma$ is trivial, there is a unique unit outward normal vector field of $\Sigma$ in $M$, denoted as $e_{r}$. Let $e_{t}$ be the unit outward normal vector field of $\Sigma$ in $N$ that is perpendicular to $e_{r}$. Since $M$ is spacelike, $e_{r}$ is spacelike and $e_{t}$ is timelike. Define a local coordinate chart $\{\theta, \phi\}$ on $\Sigma$, and let $\{\pptheta, \ppphi\}$ be the local coordinate frame. Figure \ref{fig:IMCVF steering setup} depicts the above set up, in which the mean curvature vector $\vec{H}_{\Sigma}$ is not necessarily tangential to $M$. We can extend the frame $\{e_{t}, e_{r}, \pptheta, \ppphi\}$ to a frame on a neighborhood of $\Sigma$ in $N$, and we will identify the frame with its extension. 

Since the local frame $\{e_{t}, e_{r}, \pptheta, \ppphi\}$ is not fully a coordinate frame, the commutator coefficients $C^{k}_{ij}$ defined as:
\begin{equation}\label{def:commutator coefficients}
[\alpha_{i}, \alpha_{j}] = C^{k}_{ij}\alpha_{k}, \quad \alpha_{i} \in \{e_{t}, e_{r}, \pptheta, \ppphi\}
\end{equation}
are not necessarily zero. We need these coefficients to compute the connection coefficients with respect this frame later.

With the local frame $\{e_{t}, e_{r}, \pptheta, \ppphi\}$ of the tangent bundle of $N$, let the associated dual frame of the cotangent bundle of $N$ be $\{\beta_{t}, \beta_{r}, d\theta, d\phi\}$, where
\[\beta_{i}(e_{j}) = \delta_{ij}, \quad i, j \in \{t, r\}.\]
With respect to this dual frame, we can write the spacetime metric $\gbar$ as

\begin{equation}\label{equ:local representation of gbar in dual frame}
[\gbar]= \bordermatrix{~ & e_{t} & e_{r} & \pptheta & \ppphi \cr
                  e_{t}        & -1  & 0   & 0   & 0\cr
                  e_{r}         & 0   &1    & 0   & 0\cr
                  \pptheta & 0   & 0   & a   & c\cr
                  \ppphi   & 0   & 0   & c   & b\cr}
\end{equation}
where $a, b, c$ are smooth functions on $N$. 

We want to change the metric so that $\vec{H}_{\Sigma}$ is tangential to $M$. Recall that the fourth condition (\ref{def:fourth conditions of IMCVF coordinate chart}) in the construction of the inverse mean curvature vector flow coordinates is a zeroth order equation for $d$, the $(t, r)$-metric component. This motivates the following definition:

\begin{mydef}[Steering of Spacetime Metric]\label{def:steering of spacetime metric}
Given a spacetime metric $\gbar$. A metric $\gbar_{Q}$ on $N$ is called a \emph{steering} of $\gbar$ if
\begin{equation}\label{def:steering condition}
\gbar_{Q} := \gbar + Q(\beta_{t} \otimes \beta_{r} + \beta_{r} \otimes \beta_{t})
\end{equation}
for some smooth function $Q\in C^{\infty}(N)$. $Q$ is called a \emph{steering parameter}.
\end{mydef}

\begin{figure}[!htb]
\centering
\includegraphics[scale = .45]{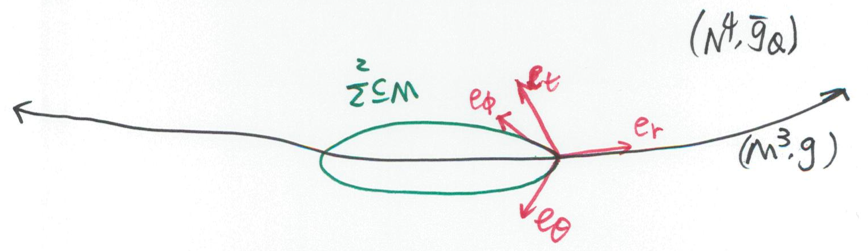}
\caption{After steering the spacetime metric $\gbar$, $e_{t}$ and $e_{r}$ are not necessarily orthogonal to each other. The metric $g$ on $M$ remains the same.}
\label{fig:after steer}
\end{figure}

Note that the coefficient matrix for the steered metric $\gbar_{Q}$ is given by:
\begin{equation}\label{equ:local representation of gbar steered in dual frame}
[\gbar_{Q}]= \bordermatrix{~ & e_{t} & e_{r} & \pptheta & \ppphi \cr
                  e_{t}        & -1  & Q   & 0   & 0\cr
                  e_{r}         & Q   &1    & 0   & 0\cr
                  \pptheta & 0   & 0   & a   & c\cr
                  \ppphi    & 0   & 0   & c   & b\cr}
\end{equation}
Thus geometrically, $e_{t}$ and $e_{r}$ are not necessarily orthogonal to each other in the steered metric (see Figure \ref{fig:after steer}).

\begin{mydef}[Area Expanding Condition]
\label{def:area expanding condition in steering flow}
$\Sigma \subset M$ given as above is said to be area expanding if
\begin{equation}\label{equ:area expanding condition}
e_{r}(ab - c^{2}) > 0.
\end{equation}
\end{mydef}

\begin{theorem}
\label{thm:mean curvature vector is tangential to a spacelike hypersurface in steered spacetime metric}
Let $(N^{4}, \gbar)$,  a spacelike hypersurface $M^{3}$ and a closed embedded surface $(\Sigma^{2}, g_{\Sigma})$ in $M$ be given as above. If $\Sigma$ is area expanding, then there exists a unique smooth steering parameter $Q = Q_{\Sigma} \in C^{\infty}(N)$, such that in the steered spacetime metric $\gbar_{Q}$, $\vec{H}_{\Sigma}$ is tangential to $M$ everywhere on $\Sigma$.
\end{theorem}

\begin{proof}
To show that there exists $Q$ such that $\vec{H}_{\Sigma}$ is tangential to $M$, it suffices to find $Q$ such that $\vec{H}_{\Sigma}$ (with respect to $\gbar_{Q}$) is parallel to $e_{r}$. 

\emph{Notation}: In the following the subscript $Q$ will be dropped for simplicity and all the spacetime metric will be referring to the steered metric $\gbar_{Q}$.

The computations are similar to the inverse mean curvature vector flow coordinates case. We divide the computations into five steps:

\underline{Step 1}: Pick a local normal variation of $\Sigma$ and let $\{\Sigma_{s}\}$ be the variational surfaces, $s\in (0, \epsilon)$. Since $\Sigma_{0} = \Sigma$ is area expanding, we can assume that this normal variation is also area expanding. Extend $\{e_{t}, e_{r}\}$ to a neighborhood of $\Sigma$ in $N$ such that $e_{r}$ remains outward unit normal to each $\Sigma_{s}$ in $M$. Define an $r$-coordinate by requiring the area of $\Sigma_{s}$ to be $A(s) =:4\pi r^{2}$. Since the variation is area expanding, $r$ is well-defined. Extend $\{\pptheta, \ppphi\}$ to a neighborhood of $\Sigma$ in $N$ such that they are constant along normal directions of $\Sigma$ in $M$. By this extension $\ppr$ is perpendicular to each $\Sigma_{s}$. Thus there exists a function $\lambda$ such that $e_{r} = \lambda\ppr$. 
\begin{equation}\label{equ:zero commutator relation}
[e_{r}, \pptheta] = [\lambda \ppr, \pptheta] = -\lambda_{,\theta}\ppt.
\end{equation}
This implies that $C^{\theta}_{r \theta} = 0$. Similarly $C^{\phi}_{r\phi} = 0$. 

\underline{Step 2}: After steering the metric $\gbar$, $e_{r}$ and $e_{t}$ are not necessarily orthogonal to each other. Let $\n$ be the normal vector of $\Sigma$ such that $\<\n, e_{r}\> = 0$. The computation of $\n$ is the same as the inverse mean curvature vector flow coordinate case (\ref{equ:normal vector field n}), and we obtain:
\begin{equation}\label{equ:normal vector field n in steering case}
\n = e_{t} - Qe_{r}.
\end{equation}
Note that $\<\n, \n\> = \<e_{t} - Qe_{r}, e_{t} - Qe_{r}\> = -1 - 2Q^{2} + Q^{2} = -(1 + Q^{2})$, which agrees with (\ref{equ:norm squred of n}) with our metric $\gbar_{Q}$.

Let $e_{n}: = \frac{\n}{||\n||} = \frac{e_{t} - Qe_{r}}{(1 + Q^{2})^{1/2}}$ be the unit timelike normal vector.

\underline{Step 3}: With respect to the orthonormal frame $\{e_{n}, e_{r}\}$ of the normal bundle $N\Sigma$, we can write the mean curvature vector as
\begin{equation}
\label{equ:mean curvature vector decomposition in steered metric}
\vec{H}_{\Sigma} = \<\vec{H}_{\Sigma}, e_{r}\>e_{r} - \<\vec{H}_{\Sigma}, e_{n}\>e_{n}.
\end{equation}

Therefore $\vec{H}_{\Sigma}$ is parallel to $e_{r}$ if and only if $\<\vec{H}_{\Sigma}, e_{n}\> = 0$. Now we compute the condition on $Q$ such that this inner product vanishes.

\begin{align}
\<\vec{H}_{\Sigma}, e_{n}\> & = \<\trace_{g_{\Sigma}}\vec{\II}, e_{n}\> = g_{\Sigma}^{ij}\<\left(\nabar_{\p_{i}}\p_{j}\right)|_{N\Sigma}, e_{n}\>\tag{$\p_{i}, \p_{j} \in \{\pptheta, \ppphi\}$}\notag\\
& = g_{\Sigma}^{ij}\left\<\<\nabar_{\p_{i}}\p_{j}, e_{r}\>e_{r} - \<\nabar_{\p_{i}}\p_{j}, e_{n}\>e_{n}, e_{n}\right\>\notag\\
& = -g_{\Sigma}^{ij}\<\nabar_{\p_{i}}\p_{j}, e_{n}\>\<e_{n}, e_{n}\> = g_{\Sigma}^{ij}\<\nabar_{\p_{i}}\p_{j}, e_{n}\>\tag{$\<e_{n}, e_{n}\> = -1$}\notag\\
& = \frac{1}{||\n||} g_{\Sigma}^{ij}\<\nabar_{\p_{i}}\p_{j}, \n\> = \frac{1}{||\n||} g_{\Sigma}^{ij}\omega^{t}_{ij}\<e_{t}, \n\>\notag\\
& = \frac{1}{||\n||} g_{\Sigma}^{ij}\omega^{t}_{ij} (- (1 + Q^{2})) = -g_{\Sigma}^{ij}\omega^{t}_{ij}\frac{1 + Q^{2}}{(1 + Q^{2})^{1/2}}\notag\\
& = - (1 + Q^{2})^{1/2} g_{\Sigma}^{ij}\omega_{ij}^{t}\label{equ:mean curvature vector in steering metric step 1}
\end{align}
where $\omega$ is the connection coefficients of the connection $\nabar_{Q}$ with respect to the metric $\gbar_{Q}$ and the local frame $\{e_{t}, e_{r}, \pptheta, \ppphi\}$, which is not necessarily a coordinate frame. This is an important difference from the inverse mean curvature vector flow coordinates case.

Recall that
\begin{equation}
\omega_{ij}^{k} = \frac{1}{2}\gbar^{kl}(\gbar_{il, j} + \gbar_{jl, i} - \gbar_{ij, l} + C_{lij} + C_{lji} - C_{ijl}),
\end{equation}
where $C_{ijl}: = \gbar_{lm}C^{m}_{ij}$, and $C_{ij}^{m}$ is the commutator coefficients defined in (\ref{def:commutator coefficients}).

\underline{Notation}: we use $\gbar_{ij, k}$ to denote $e_{k}(\gbar_{ij})$, not necessarily a coordinate derivative. 

\underline{Step 4}:
Now we compute the four connection coefficients as follows:
\begin{align}
\omega_{\theta \theta}^{t} & = \frac{1}{2}\gbar^{ti}(2\gbar_{\theta i, \theta} - \gbar_{\theta \theta, i} + 2C_{i\theta \theta} - C_{\theta \theta i})\notag\\
& = \frac{1}{2}\gbar^{ti}(2\gbar_{\theta i, \theta} - \gbar_{\theta \theta, i} + 2\gbar_{\theta j}C^{j}_{i \theta} - \cancelto{0}{\gbar_{ij}C^{j}_{\theta \theta}})\tag{$\{\pptheta, \ppphi\}$ is a coordinate frame}\notag\\
& = \frac{1}{2}\Big\{\gbar^{tt}(\cancelto{0}{2\gbar_{\theta t, \theta}} - \gbar_{\theta \theta, t} + 2\gbar_{\theta j}C^{j}_{t\theta}) + \gbar^{tr}(\cancelto{0}{2\gbar_{\theta r, \theta}} - \gbar_{\theta \theta, r} + 2\gbar_{\theta j}C^{j}_{r\theta}) + 0 + 0\Big\}\notag\\
& = \frac{1}{2}\Big\{\gbar^{tt}(-\gbar_{\theta \theta, t} + 2\gbar_{\theta \theta}C^{\theta}_{t\theta} + 2\gbar_{\theta \phi}C^{\phi}_{t\theta}) + \gbar^{tr}(-\gbar_{\theta \theta, r} + 2\gbar_{\theta \theta}C^{\theta}_{r\theta} + 2\gbar_{\theta \phi}C^{\phi}_{r\theta}) \Big\}\notag\\
& = \frac{ab - c^{2}}{2|\gbar|}\Big\{(-a_{,t} + 2aC^{\theta}_{t\theta} + 2cC^{\phi}_{t\theta}) - Q(-a_{,r} + 2aC^{\theta}_{r\theta} + 2cC^{\phi}_{r\theta}) \Big\}\label{equ:omega t theta theta}
\end{align}

\begin{align}
\omega_{\theta \phi}^{t} & = \frac{1}{2}\gbar^{ti}(\gbar_{\phi i, \theta} + \gbar_{\theta i, \phi} - \gbar_{\theta \phi, i} + C_{i \theta \phi} + C_{i \phi \theta} - C_{\theta \phi i})\notag\\
& = \frac{1}{2}\gbar^{ti}(\gbar_{\phi i, \theta} + \gbar_{\theta i, \phi} - \gbar_{\theta \phi, i} + \gbar_{\phi j}C^{j}_{i \theta} + \gbar_{\theta j}C^{j}_{i \phi} - \cancelto{0}{\gbar_{ij}C^{j}_{\theta \phi})}\notag\\
& = \frac{1}{2}\Big\{\gbar^{tt}(\cancelto{0}{\gbar_{\phi t, \theta}} + \cancelto{0}{\gbar_{\theta t, \phi}} - \gbar_{\theta \phi, t} + \gbar_{\phi \theta}C^{\theta}_{t\theta} + \gbar_{\phi \phi}C^{\phi}_{t\theta} + \gbar_{\theta \theta}C^{\theta}_{t\phi} + \gbar_{\theta \phi}C^{\phi}_{t\phi})\notag\\
& + \gbar^{tr}\cancelto{0}{\gbar_{\phi r, \theta}} + \cancelto{0}{\gbar_{\theta r, \phi}}  - \gbar_{\theta \phi, r} + \gbar_{\phi \theta}C^{\theta}_{r\theta} + \gbar_{\phi \phi}C^{\phi}_{r\theta} + \gbar_{\theta \theta}C^{\theta}_{r\phi} + \gbar_{\theta \phi}C^{\phi}_{r\phi})   \Big\}\notag\\
& = \frac{ab - c^{2}}{2|\gbar|}\Big\{(-c_{,t} + cC^{\theta}_{t\theta} + bC^{\phi}_{t\theta} + aC^{\theta}_{t \phi} + cC^{\phi}_{t\phi}) - Q(-c_{, r} + cC^{\theta}_{r\theta} + bC^{\phi}_{r\theta} + aC^{\theta}_{r\phi} \notag\\
& + cC^{\phi}_{r\phi})  \Big\}\label{equ:omega t theta phi}
\end{align}

\begin{align}
\omega_{\phi \theta}^{t} & = \frac{1}{2}\gbar^{ti}(\gbar_{\theta i, \phi} + \gbar_{\phi i, \theta} - \gbar_{\phi \theta, i} + C_{i\phi \theta} + C_{i\theta \phi} - C_{\phi \theta i})\notag\\
& = \frac{1}{2}\gbar^{ti}(\gbar_{\theta i, \phi} + \gbar_{\phi i, \theta} - \gbar_{\phi \theta, i} + \gbar_{\theta j}C^{j}_{i\phi} + \gbar_{\phi j}C^{j}_{i\theta} - \cancelto{0}{\gbar_{ij}C^{j}_{\phi \theta})}\notag\\
& = \omega_{\theta \phi}^{t}\label{equ:omega t phi theta} 
\end{align}

\begin{align}
\omega_{\phi \phi}^{t} & = \frac{1}{2}\gbar^{ti}(2\gbar_{\phi i, \phi} - \gbar_{\phi \phi, i} + 2C_{i\phi \phi} - C_{\phi \phi i})\notag\\
& = \frac{1}{2}\gbar^{ti}(2\gbar_{\phi i, \phi} - \gbar_{\phi \phi, i} + 2\gbar_{\phi j}C^{j}_{i\phi} - \cancelto{0}{\gbar_{ij}C^{j}_{\phi \phi}})\notag\\
& = \frac{1}{2}\Big\{\gbar^{tt}(\cancelto{0}{2\gbar_{\phi t, \phi}} - \gbar_{\phi \phi, t} + 2\gbar_{\phi \theta}C^{\theta}_{t \phi} + 2\gbar_{\phi \phi}C^{\phi}_{t\phi}) + \gbar^{tr}(\cancelto{0}{2\gbar_{\phi r, \phi}} - \gbar_{\phi \phi, r} + 2\gbar_{\phi \theta}C^{\theta}_{r\phi} \notag\\
& + 2\gbar_{\phi \phi}C^{\phi}_{r\phi})\Big\}\notag\\
& = \frac{ab - c^{2}}{2|\gbar|}\Big\{(-b_{,t} + 2cC^{\theta}_{t\phi} + 2bC^{\phi}_{t\phi}) - Q(-b_{,r} + 2cC^{\theta}_{r\phi} + 2bC^{\phi}_{r\phi})   \Big\}\label{equ:omega t phi phi}
\end{align}

\underline{Step 5}:
Now plug (\ref{equ:omega t theta theta}), (\ref{equ:omega t theta phi}), (\ref{equ:omega t phi theta}) and (\ref{equ:omega t phi phi}) back into (\ref{equ:mean curvature vector in steering metric step 1}), we get:
\begin{align}
& \<\vec{H}_{\Sigma}, e_{n}\> = - (1 + Q^{2})^{1/2}\frac{1}{ab - c^{2}}(b\omega_{\theta \theta}^{t} - c\omega_{\theta \phi}^{t} - c\omega_{\phi \theta}^{t} + a\omega_{\phi \phi}^{t}) \notag\\
& = - (1 + Q^{2})^{1/2}\frac{1}{ab - c^{2}}(b\omega_{\theta \theta}^{t} - 2c\omega_{\theta \phi}^{t} + a\omega_{\phi \phi}^{t})\tag{By (\ref{equ:omega t phi theta})} \notag\\
& = - (1 + Q^{2})^{1/2}\frac{1}{ab - c^{2}}\frac{ab - c^{2}}{2|\gbar|}\Big\{(-a_{,t}b + 2abC^{\theta}_{t\theta} + \cancel{2bcC^{\phi}_{t\theta}} + 2cc_{,t} - 2c^{2}C^{\theta}_{t \theta} - \cancel{2bcC^{\phi}_{t\theta}}\notag\\
& - \bcancel{2acC^{\theta}_{t\phi}} - 2c^{2}C^{\phi}_{t\phi} - ab_{,t} + \bcancel{2acC^{\theta}_{t\phi}} + 2abC^{\phi}_{t\phi}) -Q (-a_{,r}b + 2abC^{\theta}_{r\theta} + \cancel{2bcC^{\phi}_{r\theta}} \notag\\
& + 2cc_{,r} - 2c^{2}C^{\theta}_{r\theta} - \cancel{2bcC^{\phi}_{r \theta}} - \bcancel{2acC^{\theta}_{r\phi}} - 2c^{2}C^{\phi}_{r\phi} -ab_{,r} + \bcancel{2acC^{\theta}_{r\phi}} + 2abC^{\phi}_{r\phi})\Big\}\notag\\
& = -\frac{(1 + Q^{2})^{1/2}}{2|\gbar|}\Big\{ -e_{t}(ab - c^{2}) + 2(ab - c^{2})(C^{\theta}_{t\theta} + C^{\phi}_{t\phi})\notag\\
& -Q\left(-e_{r}(ab - c^{2}) + 2(ab - c^{2})\cancelto{0}{(C^{\theta}_{r\theta} + C^{\phi}_{r\phi})} \right)\Big\}\notag\\
& = -\frac{(1 + Q^{2})^{1/2}}{2|\gbar|}\Big\{e_{r}(ab - c^{2}) Q -e_{t}(ab - c^{2}) + 2(ab - c^{2})(C^{\theta}_{t\theta} + C^{\phi}_{t\phi})   \Big\}\label{equ:equ:mean curvature vector in steering metric step 2}
\end{align}
since $C^{\theta}_{r \theta} = C^{\phi}_{r\phi} = 0$ by Equation (\ref{equ:zero commutator relation}).
Therefore $\vec{H}_{\Sigma}$ is parallel to $e_{r}$ if and only if 
\begin{equation}\label{equ:condition on steering parameter}
e_{r}(ab - c^{2}) Q - e_{t}(ab - c^{2}) + 2(ab - c^{2})(C^{\theta}_{t\theta} + C^{\phi}_{t\phi}) = 0.
\end{equation}
Notice that Equation (\ref{equ:condition on steering parameter}) is zeroth order in $Q$. Since $\Sigma$ is area expanding, $e_{r}(ab - c^{2}) > 0$, and hence we get a unique solution
\begin{equation}\label{equ:unique solution of the steer parameter assuming steerable condition}
Q = \frac{e_{t}(ab - c^{2}) - 2(ab - c^{2})(C^{\theta}_{t\theta} + C^{\phi}_{t\phi})}{e_{r}(ab - c^{2})}.
\end{equation}
\end{proof}

\begin{lemma}\label{lemma:mean curvature scalars and variations of area form}
$e_{r}(ab - c^{2}) = 0$ if and only if  $H_{e_{r}} = -\<\vec{H}_{\Sigma}, e_{r}\> = 0$, i.e. $\Sigma$ is a minimal surface in $M$, where $H_{e_{r}}$ is the mean curvature scalar of $\Sigma$ in the direction of $e_{r}$.
\end{lemma}

\begin{proof}
For the first claim, note that
\begin{equation}
e_{r}(dA_{\Sigma}) = e_{r}(\sqrt{ab - c^{2}}\beta_{\theta}\beta_{\phi}) = \frac{1}{2\sqrt{ab - c^{2}}}e_{r}(ab - c^{2})\beta_{\theta}\beta_{\phi} = \frac{1}{2(ab - c^{2})}e_{r}(ab - c^{2})dA_{\Sigma}.
\end{equation}
On the other hand by the first variation of area formula (\ref{equ:first variation of area}), 
\begin{equation}
e_{r}(dA_{\Sigma}) = -\<\vec{H}_{\Sigma}, e_{r}\>dA_{\Sigma} = H_{e_{r}}dA_{\Sigma}.
\end{equation}
Combing the two equations we get:
\begin{equation}\label{equ:radial mean curvature scalar and radial variation of area form}
2H_{e_{r}}\cdot (ab - c^{2}) = e_{r}(ab - c^{2}).
\end{equation}
Therefore $e_{r}(ab - c^{2}) = 0$ if and only if $H_{e_{r}} = 0$. 
\end{proof}

An application of Theorem \ref{thm:mean curvature vector is tangential to a spacelike hypersurface in steered spacetime metric} is to generate more examples of inverse mean curvature vector flow solutions. Let $(N^{4}, \gbar)$ be a spacetime and $(M, g)$ a spacelike hypersurface with induced metric $g$. Suppose $\{\Sigma_{s}\}$ is a solution to the smooth inverse mean curvature flow in $M$. Let $\vec{H}_{s}$ and $dA_{s}$ be the mean curvature vector field and area form of $\Sigma_{s}$, respectively. By the first variation formula, 
\begin{equation}\label{equ:first variation of area form under IMCF revisit}
\frac{d}{ds}(dA_{s}) = dA_{s}.
\end{equation}
Theorem \ref{thm:mean curvature vector is tangential to a spacelike hypersurface in steered spacetime metric} allows us to steer the metric smoothly along inverse mean curvature flow to keep $\vec{H}_{s}$ tangential to $M$. Note that (\ref{equ:first variation of area form under IMCF revisit}) still holds after the steering since we are not changing the metric on $M$. Therefore by Proposition \ref{prop:spacelike IMCVF equals tangential mean curvature in IMCF}, in the steered metric $\{\Sigma_{s}\}$ is a solution to smooth inverse mean curvature vector flow equation.

One could generalize this technique to study weak solutions to inverse mean curvature vector flow equation (defined in (\ref{def:weak solution to IMCVF})) using solutions to Huisken-Ilmanen inverse mean curvature flow (with jumps) in a hypersurface, but we will not give a rigorous treatment here.

\section{Generalizations}
\label{sec:generalizations of imcvf}

Given a surface $(\Sigma, g_{\Sigma})$ inside a spacetime $(N^{4}, \gbar)$. Let $\vec{I}$ be the inverse mean curvature vector. A more general flow than inverse mean curvature vector flow is to flow out $\Sigma$ in the following direction:
\begin{equation}\label{equ:uniformly area expanding flow with constant beta}
\vec{\xi} := \vec{I} + \beta\vec{I}^{\perp},
\end{equation}
where the $\perp$ operation is a linear isomorphism on the normal bundle defined in Section \ref{subsec:normal bundle geometry}, and $\beta$ is a constant on each flow surface (hence is only a function of the flow parameter) such that $-1\leq \beta \leq 1$.

Therefore inverse mean curvature vector flow corresponds to the case where $\beta = 0$. The procedure for constructing spacetimes with inverse mean curvature vector flow solutions can be generalized to constructing spacetimes in which this more general flow exists. The idea is to construct a spacetime metric $\gbar$ that admit a coordinate chart $\{t, r, \theta, \phi\}$, such that conditions (\ref{def:first conditions of IMCVF coordinate chart}), (\ref{def:second conditions of IMCVF coordinate chart}), (\ref{def:third conditions of IMCVF coordinate chart}) and the fourth condition:
\begin{equation}\label{equ:radial flow direction in uniformly area expanding flow with constant beta}
\vec{\xi}\,\,\text{is parallel to}\,\, \ppr.
\end{equation}
are satisfied.
\chapter{Uniformly Area Expanding Straight Out Flows}
\label{chap:straightoutflow}

In this chapter we study \emph{uniformly area expanding straight out flows}, or simply \emph{straight out flows} of spacelike surfaces in a spacetime. Consider a spacetime $(N^{4}, \gbar)$ and a closed embedded surfaces $\Sigma$ with the induced metric $g_{\Sigma}$. A straight out direction of $\Sigma$, first studied by M. Mars, E. Malec and Simon \cite{MalecMarsetal2002}, is a normal vector field that has ``minimal variations'' along $\Sigma$. Such a normal vector field is a minimizer of a natural functional defined on the normal bundle. The Hawking mass is also monotonically non-decreasing under smooth straight out flows assuming the spacetime satisfies the dominant energy condition. A condition for a spacetime to admit straight out flow coordinate charts is derived in this chapter. Complete understanding of such spacetimes is still work in progress.

\section{Background and Straight Out Flow Coordinate Chart}
\label{sec:background and straight out flow coordinate chart}
Let $(N^{4}, \gbar)$ be a time-oriented spacetime and $(\Sigma^{2}, g_{\Sigma})$ be a closed embedded surface. The normal bundle of $\Sigma$, $N\Sigma$, has an induced metric of signature $(-, +)$. On each fiber $N_{p}\Sigma$, the nonzero vectors get partitioned into four quadrants: outward-spacelike, inward-spacelike, future-timelike and past-timelike. Let $U^{+}N (\Sigma)$ denote the subbundle of $N\Sigma$ that consists of outward-spacelike normal vector fields of unit length. Given a smooth section $\nu$ of $U^{+}N(\Sigma)$, its associated connection one-form $\alpha_{\nu}$ on $\Sigma$ is defined by:
\begin{equation}
\label{def:connection one-form revisit}
\alpha_{\nu}(X): = \<\naperp_{X}v, v^{\perp}\>_{\gbar}, \quad \forall\, X\in \Gamma(T\Sigma).
\end{equation}
where $\perp$ is a fiberwise linear isomorphism defined in Section \ref{sec:codim two surface geometry}. Notice that $\{\nu, \nu^{\perp}\}$ forms an orthonormal frame of $N\Sigma$.

Given another smooth section $\overline{\nu}$ of $U^{+}N(\Sigma)$, there exists a constant $\theta > 0$ such that
\begin{align}
& \overline{\nu} = \cosh{\theta}\nu + \sinh{\theta}\nu^{\perp}\notag\\
& \overline{\nu}^{\perp} = \sinh{\theta}\nu + \cosh{\theta}\nu^{\perp}\label{equ:change of frame corresponds to hyperbolic rotation}
\end{align}
 
Geometrically $\{\overline{\nu}, \overline{\nu}^{\perp}\}$ can be viewed as hyperbolic rotation of $\{\nu, \nu^{\perp}\}$ by angle $\theta$. The associated connection one-forms are related by:
\begin{lemma}[\cite{BrayJauregui2013}]
\label{lemma:relationships of hyperbolically rotated connection one forms}
Let $\alpha_{\overline{\nu}}$ be the associated connection one-form of $\overline{\nu}$, then $\alpha_{\overline{\nu}}$ is related to $\alpha_{\nu}$ by
\begin{equation}\label{equ:change of connection one form when we change the frame}
\alpha_{\overline{\nu}} = \alpha_{v} - d\theta.
\end{equation}
\end{lemma}

Let $E$ be an energy functional on $U^{+}N(\Sigma)$ such that:
\begin{equation}\label{def:energy functional on normal bundle}
E(\nu): = \int_{\Sigma}||\na^{\perp}\nu||^{2} \, dA_{\Sigma}, \quad \forall\, \nu \in \Gamma\big(U^{+}N(\Sigma) \big).
\end{equation}

H. Bray and J. Jauregui proved the following proposition:
 
\begin{proposition}[\cite{BrayJauregui2013}]\label{prop:straight out direction is minimizer of energy functional}
$\nu \in \Gamma(N\Sigma)$, outward spacelike and is of unit length, is a minimizer of $E$ if and only if $\divg_{\Sigma}(\alpha_{\nu}) = 0$. 
\end{proposition}

Intuitively, a minimizer of $E$ is a normal vector field with ``minimal variations'' along $\Sigma$. Bray and Jauregui also proved the existence of such minimizers:

\begin{proposition}[\cite{BrayJauregui2013}]\label{prop:existence of divergence free sections}
Minimizers of $E$ exist, and any two such minimizers differ by a hyperbolic rotation as in (\ref{equ:change of frame corresponds to hyperbolic rotation}) of a constant angle $\theta$.
\end{proposition}

\begin{proof}
We present the proof in \cite{BrayJauregui2013}. Fix a section $\nu$ of $U^{+}N(\Sigma)$. Given any other section $\overline{\nu}$, $\overline{\nu}$ and $\nu$ are related by a hyperbolic rotation as in (\ref{equ:change of frame corresponds to hyperbolic rotation}). Let $\alpha_{\nu}$ and $\alpha_{\overline{\nu}}$ be the associated connection one-forms of $\nu$ and $\overline{\nu}$, respectively. By Lemma \ref{lemma:relationships of hyperbolically rotated connection one forms}, $\alpha_{\overline{\nu}}$ is a divergence free if and only if 
\begin{equation}
0 = \divg_{\Sigma}(\alpha_{\overline{\nu}}) = -d^{*}\alpha_{\overline{\nu}} = -d^{*}\alpha_{\nu} + d^{*}d\theta = \divg_{\Sigma}(\alpha_{\nu}) + \Delta_{\Sigma}\theta,
\end{equation}
where $d^{*}$ is the $L^{2}$-adjoint of $d$ on $\Sigma$. The above equation is equivalent to:
\begin{equation}\label{equ:elliptic PDF for existence of straight out normal}
\Delta_{\Sigma}\theta = \divg_{\Sigma}(\alpha_{\nu}).
\end{equation}
This is a Poisson equation. By the divergence theorem, $\displaystyle\int_{\Sigma}\divg_{\Sigma}(\alpha_{\nu})\, dA_{\Sigma} = 0$. Therefore Equation (\ref{equ:elliptic PDF for existence of straight out normal}) has smooth solutions. Moreover, any two solutions differ by an element in the kernel of $\Delta_{\Sigma}$, which consists of constant functions on $\Sigma$.
\end{proof}

A minimizer of $E$ can be viewed as a \emph{straight out} direction of the surface $\Sigma$, as it tries to level the surface as much as possible. This motivates the following definition:

\begin{mydef}[Uniformly Area Expanding Straight Out Flow]\label{def:straight out flow}
Given a surface $\Sigma$ in a time-oriented spacetime $(N^{4},\gbar)$. Let $\vec{T}$ be a global timelike vector field. A smooth straight out flow of $\Sigma$ in $N$ is a smooth normal variation $F: \Sigma \times[0, T] \longrightarrow N$ such that
\begin{equation}\label{def:straight out flow equation}
\frac{\p}{\p s}F(x, s) = \beta\nu_{s}(x), \quad \divg_{\Sigma_{s}}(\alpha_{\nu_{s}}) = 0, \,\, \int_{\Sigma_{s}}\<\vec{T}, \nu_{s}\> \, dA_{s} = 0,
\end{equation}
where $\beta$ is chosen such that $\frac{d}{ds}(dA_{s}) = dA_{s}$; and $s\in [0, T]$, $\Sigma_{s}: = F(\Sigma, s)$, and $\nu_{s}$ is a outward spacelike unit normal vector field along $\Sigma_{s}$ with zero divergence.
\end{mydef}

On each surface of $\Sigma_{s}$ in Definition \ref{def:straight out flow}, the outward-spacelike vector $\nu_{s}$ with zero divergence is only unique up to a hyperbolic rotation with angle $\theta$ (see Proposition \ref{prop:existence of divergence free sections}). The additional condition that $\displaystyle\int_{\Sigma_{s}}\<\vec{T}, \nu_{s}\>\, dA_{s} = 0$ defines each $\nu_{s}$ uniquely. This condition is due to J. Jauregui. 

E. Malec, M. Mars and W. Simon have studied such flows in \cite{MalecMarsetal2002}, and they have shown that given a smooth solution $\{\Sigma_{t}\}$ to the straight out flow equation (\ref{def:straight out flow equation}), 
\[\frac{d}{dt}m_{H}(\Sigma_{s}) \geq 0,\]
i.e. that Hawking mass is also monotonically non-decreasing (also see \cite{BrayJauregui2014}).

The goal in this chapter is to construct spacetimes that admit smooth solutions to straight out flow equation. The idea is similar to the construction in Chapter \ref{chap:imcvf}. Given $U: = (\RR^{3}\backslash B_{1})\times \RR$, where $B_{1}$ is the unit ball in $\RR^{3}$, we seek a spacetime metric $\gbar$ on $U$ that admit a coordinate chart $\{t, r, \theta, \phi\}$, such that $\gbar$ satisfies the following four conditions:

\begin{align}
& (1) \quad  \<\ppr, \pptheta\> = 0;\label{def:first conditions of straight out coordinate chart}\\
& (2) \quad  \<\ppr, \ppphi\> = 0;\label{def:second conditions of straight out coordinate chart}\\
& (3) \quad \text{Area form of coordinate sphere $S_{t, r}$ satisfies}\quad dA_{t, r} = r^{2}\sin{\theta}d\theta d\phi; \label{def:third conditions of straight out coordinate chart}\\
& (4) \quad  e_{r} := \frac{1}{u}\ppr\,\, \text{is straight out, i.e.}\,\, \divg_{g_{S}} \alpha_{e_{r}} = 0;\label{def:fourth conditions of straight out coordinate chart}
\end{align}
where $\alpha_{e_{r}}$ is the connection one-form associated with $e_{r}$.

The following definition is analogous to Definition \ref{def:IMCVF coordinates}:
\begin{mydef}[Straight Out Flow Coordinate Chart]\label{def:straight out flow coordinate chart}
If a spacetime $(N^{4}, \gbar)$ admits a coordinate chart $\{t, r, \theta, \phi\}$ such that the four conditions (\ref{def:first conditions of straight out coordinate chart}), (\ref{def:second conditions of straight out coordinate chart}), (\ref{def:third conditions of straight out coordinate chart}) and (\ref{def:fourth conditions of straight out coordinate chart}) are satisfied, then $\{t, r, \theta, \phi\}$ is called a straight out coordinate chart, and $N$ is called a spacetime that admits a straight out coordinate chart. 
\end{mydef}

\section{Construction of Spacetimes That Admit Uniformly Area Expanding Straight Out Flow Coordinate Charts (work in progress)}
\label{sec:construction of spacetimes that admit straight out flow coordinate charts}

Let $U = (\RR^{3}\backslash B_{1})\times \RR$. Given a coordinate chart $\{t, r, \theta, \phi\}$, define a spacetime metric $\gbar$ to be
\begin{equation}\label{equ:local representation of gbar in straight out flow coordinates step 1}
\overline{g}: = \bordermatrix{~ & t & r & \theta & \phi \cr
                  t & -v^{2} & d & e & f\cr
                  r & d & u^{2} & 0 & 0\cr
                  \theta & e & 0 & a & c\cr
                  \phi & f & 0 & c & b\cr}
\end{equation}
for smooth functions $a,b,c,d,e,f,u,v$ on $U$.  $\gbar$ satisfies conditions (\ref{def:first conditions of straight out coordinate chart}) and (\ref{def:second conditions of straight out coordinate chart}). Choosing two of the three variables $a, b$ and $c$ such that $ab - c^{2} = r^{4}\sin^{2}{\theta}$ satisfies condition (\ref{def:third conditions of straight out coordinate chart}). We show that the fourth condition (\ref{def:fourth conditions of straight out coordinate chart}) is a second order elliptic PDE in the metric component $d$.

The normal bundle of the coordinate sphere $S_{t, r}$ has an orthonormal frame $\{e_{r}, e_{n}\}$, where $e_{r} = \frac{1}{u}\ppr$, and $e_{n} = \frac{\n}{||\n||_{\gbar}}$. The normal vector $\n$ is the same as in Equation (\ref{equ:normal vector field n}). We compute the condition on $\gbar$ such that $e_{r}$ is straight out on each $S_{t, r}$.
 
\emph{Notation}: For the rest of this chapter, the induced metric on $S_{t, r}$ is denoted as $g_{S}$; $\nabar$ is the Levi-Civita connection with respect to $\gbar$; $|\gbar|$ and $|g_{S}|$ are determinants of $\gbar$ and $g_{S}$ respectively; and the subscript $e_{r}$ is omitted from $\alpha_{e_{r}}$ whenever there is no confusion. 

\subsection{Computation of The Connection One Form and Its Divergence}
\label{subsec:computation of the connection one form and its divergence}
In terms of the frame $\{d\theta, d\phi\}$ for $T^{*}S_{t, r}$, write $\alpha = \alpha_{\theta}d\theta + \alpha_{\phi}d\phi$, where
\begin{align*}
\alpha_{\theta}: & = \alpha\left(\pptheta\right) = \<\nabar^{\perp}_{\pptheta}e_{r}, e_{r}^{\perp}\> \tag{By definition \ref{def:connection one-form revisit}}\\
& =  \<\nabar^{\perp}_{\pptheta}e_{r}, e_{n}\> = \left\<\nabar^{\perp}_{\pptheta}\left(\frac{1}{u}\ppr\right)\,,\, \frac{\n}{||\n||}\right\>\\
& = \left\<(1/u)_{,\theta}\ppr + \frac{1}{u}\nabar^{\perp}_{\pptheta}\ppr\,,\, \frac{\n}{|\n|}\right\> \\
& = \frac{1}{u||\n||}\left\<\nabar^{\perp}_{\pptheta}\ppr, \n\right\> =  \frac{1}{u||\n||}\left\<\nabar_{\pptheta}\ppr, \n\right\> \tag{$\n$ is normal}\\
& = \frac{1}{u||\n||}\Gammabar^{t}_{\theta r}\left\<\ppt, \n\right\> \tag{$\n$ is perpendicular to $\ppr, \pptheta, \ppphi$}\\
& = \frac{1}{u||\n||}\Gammabar^{t}_{\theta r} \cdot \frac{|\gbar|}{u^{2}|g_{S}|}\tag{By (\ref{equ:norm squred of n})}\\
& = \frac{1}{u} \left(\frac{u^{2}|g_{S}|}{-|\gbar|} \right)^{1/2}\Gammabar^{t}_{\theta r} \cdot \frac{|\gbar|}{u^{2}|g_{S}|} \tag{By (\ref{equ:norm of n})}\\
& = -\frac{1}{u^{2}}\Gammabar^{t}_{\theta r} \left(\frac{-|\gbar|}{|g_{S}|} \right)^{1/2} \tag{$|\gbar| = -(|\gbar|^{2})^{1/2}$, since $|\gbar| < 0$}
\end{align*}
Therefore
\begin{equation}\label{theta component of alpha in straight out flow}
\alpha_{\theta} =  -\frac{1}{u^{2}}\Gammabar^{t}_{\theta r} \left(\frac{-|\gbar|}{|g_{S}|} \right)^{1/2}.
\end{equation}
By a similar computation, we have:
\begin{equation}\label{phi component of alpha in straight out flow}
\alpha_{\phi} :=  \alpha\left(\ppphi\right) = -\frac{1}{u^{2}}\Gammabar^{t}_{\phi r} \left(\frac{-|\gbar|}{|g_{S}|} \right)^{1/2}.
\end{equation}
Therefore:
\begin{equation}\label{equ:connection one form in straight out flow}
\alpha = -\frac{1}{u^{2}}\Gammabar^{t}_{\theta r} \left(\frac{-|\gbar|}{|g_{S}|} \right)^{1/2} \cdot d\theta  -\frac{1}{u^{2}}\Gammabar^{t}_{\phi r} \left(\frac{-|\gbar|}{|g_{S}|} \right)^{1/2}\cdot d\phi.
\end{equation}

The divergence of $\alpha$ is defined to be the divergence of its dual vector field, which is $\beta = \beta^{\theta}\pptheta + \beta^{\phi}\ppphi$, where 
\begin{equation}\label{dual vector field of connection one form in straight out flow}
\beta^{\theta} = g_{S}^{\theta \theta}\alpha_{\theta} + g_{S}^{\theta \phi}\alpha_{\phi},\quad \beta^{\phi} = g_{S}^{\phi \theta}\alpha_{\theta} + g_{S}^{\phi \phi}\alpha_{\phi}
\end{equation}

Recall that 
\begin{mydef}\label{def:divergence of vector field}
For a vector field $X = \sum_{i = 1}^{n}X^{i}\ppxi$ in local coordinates, its divergence with respect to a metric $g$ is given by:
\[\divg_{g}X: = \frac{1}{\sqrt{|g|}}\sum_{i = 1}^{n}\ppxi\left(X^{i} \sqrt{|g|} \right).\]
\end{mydef}

Therefore the divergence of $\alpha$ along $S_{t, r}$ is
\begin{align*} 
\divg_{g_{S}}\alpha: & = \divg_{g_{S}}\beta = \frac{1}{\sqrt{|g_{S}|}}\sum_{i = 1}^{2}\ppxi\left(\beta^{i}\sqrt{|g_{S}|}\right)\tag{$\ppxi\in \{\pptheta, \ppphi\}$}\\
& = \frac{1}{\sqrt{|g_{S}|}}\left[\pptheta \left(\beta^{\theta} \sqrt{|g_{S}|}\right) + \ppphi \left(\beta^{\phi} \sqrt{|g_{S}|} \right) \right]\\
& = \frac{1}{r^{2}\sin{\theta}}\left[\pptheta \left(\beta^{\theta} r^{2}\sin{\theta}\right) + \ppphi \left(\beta^{\phi} r^{2}\sin{\theta} \right) \right]\\
& = \frac{1}{r^{2}\sin{\theta}}\left(\beta^{\theta}_{,\theta} \cdot r^{2}\sin{\theta} + \beta^{\theta} \cdot r^{2}\cos{\theta} + \beta^{\phi}_{,\phi} \cdot r^{2}\sin{\theta} + 0 \right)
\end{align*}
that is:
\begin{equation}\label{equ:divergence of connection one form in straight out flow}
\divg_{g_{S}}\alpha = \beta^{\theta}_{,\theta} + \beta^{\phi}_{,\phi} + \beta^{\theta}\cot{\theta}.
\end{equation}

\subsection{Straight Out Flow Coordinate Chart: Big Picture}
\label{subsec:straight out flow equation big picture}

We expand out Equation (\ref{equ:divergence of connection one form in straight out flow}), separating higher order derivatives from lower order derivatives.

\begin{align}
\divg_{g_{S}}\alpha & = \beta^{\theta}_{,\theta} + \beta^{\phi}_{,\phi} + \beta^{\theta}\cot{\theta}\notag\\
& = (g_{S}^{\theta \theta} \alpha_{\theta} + g_{S}^{\theta \phi} \alpha_{\phi})_{,\theta} + (g_{S}^{\phi \theta} \alpha_{\theta} + g_{S}^{\phi \phi} \alpha_{\phi})_{,\phi} + (g_{S}^{\theta \theta} \alpha_{\theta} + g_{S}^{\theta \phi} \alpha_{\phi}) \cot{\theta}\notag \\
& = \Red(g_{S}^{\theta \theta} \alpha_{\theta, \theta} + g_{S}^{\theta \phi} \alpha_{\phi, \theta} + g_{S}^{\phi \theta}\alpha_{\theta, \phi} + g_{S}^{\phi \phi}\alpha_{\phi, \phi})\Black \tag{\Red 2nd derivative in $\gbar$\Black}\notag \\ 
& + \CornflowerBlue (g_{S, \theta}^{\theta \theta}\alpha_{\theta} + g_{S,\theta}^{\theta \phi}\alpha_{\phi} + g_{S, \phi}^{\phi \theta}\alpha_{\theta} + g_{S, \phi}^{\phi \phi} \alpha_{\phi}) + (g_{S}^{\theta \theta} \alpha_{\theta} + g_{S}^{\theta \phi} \alpha_{\phi}) \cot{\theta} \Black\tag{\CornflowerBlue 1st derivative in $\gbar$\Black}\notag \\
& =: \Red(B)\Black + \CornflowerBlue (l1)\Black \label{equ:lower order term in straight out flow 1}
\end{align}

From (\ref{theta component of alpha in straight out flow}) and (\ref{phi component of alpha in straight out flow}), $\alpha_{\theta}$ and $\alpha_{\phi}$ are first order in the derivative of the spacetime metric $\gbar$. Therefore the divergence equation (\ref{equ:divergence of connection one form in straight out flow}) is second order in the derivative of $\gbar$. In (\ref{equ:lower order term in straight out flow 1}) above, $(B)$ denotes the second derivative terms, and $(l1)$ denotes the first derivative terms. 

To compute $\Red(B)\Black$, notice that:
\begin{align*}
g_{S}^{\theta \theta} \alpha_{\theta, \theta} & = g_{S}^{\theta \theta} \left(-\frac{1}{u^{2}}\Gammabar^{t}_{\theta r} \left(\frac{-|\gbar|}{|g_{S}|} \right)^{1/2} \right)_{,\theta} \tag{plug in (\ref{theta component of alpha in straight out flow}) for $\alpha_{\theta}$}\\
& = g_{S}^{\theta \theta}\Bigg(\CornflowerBlue\left(-\frac{1}{u^{2}} \right)_{,\theta} \Gammabar^{t}_{\theta r} \left(\frac{-|\gbar|}{|g_{S}|} \right)^{1/2} \Black + \left(-\frac{1}{u^{2}}\right)\Gammabar^{t}_{\theta r,\theta} \left(\frac{-|\gbar|}{|g_{S}|} \right)^{1/2} \\
& + \CornflowerBlue\left(-\frac{1}{u^{2}}\right)\Gammabar^{t}_{\theta r} \left(\frac{-|\gbar|}{|g_{S}|} \right)^{1/2}_{,\theta}\Black \Bigg)\\
& =: -\frac{1}{u^{2}}\left(\frac{-|\gbar|}{|g_{S}|} \right)^{1/2}g_{S}^{\theta \theta} \cdot \Gammabar^{t}_{\theta r,\theta} + \CornflowerBlue \text{lower order derivatives in $\gbar$}\Black.
\end{align*}

Similarly for the other three terms in $(B)$, the highest order derivative terms are the terms containing derivatives of the Christoffel symbols. Therefore
\begin{equation}\label{equ:lower order term in straight out flow 2}
\Red(B)\Black = \Red(C) \Black + \CornflowerBlue(l2)\Black,
\end{equation}
where 
\begin{align}
\Red (C) \Black &: =  -\frac{1}{u^{2}}\left(\frac{-|\gbar|}{|g_{S}|} \right)^{1/2} \left(g_{S}^{\theta \theta} \Gammabar^{t}_{\theta r, \theta} + g_{S}^{\theta \phi} \Gammabar^{t}_{\phi r, \theta} + g_{S}^{\phi \theta} \Gammabar^{t}_{\theta r, \phi} + g_{S}^{\phi \phi}\Gammabar^{t}_{\phi r, \phi} \right)  \notag\\
& = -\frac{1}{u^{2}}\left(\frac{-|\gbar|}{|g_{S}|} \right)^{1/2} \cdot \frac{1}{|g_{S}|} \left(b \Gammabar^{t}_{\theta r, \theta} -c \Gammabar^{t}_{\phi r, \theta} -c \Gammabar^{t}_{\theta r, \phi} + a\Gammabar^{t}_{\phi r, \phi}  \right) \notag \\
& =  -\frac{1}{u^{2}}\frac{(-|\gbar|)^{1/2}}{|g_{S}|^{3/2}}\cdot \left(b \Gammabar^{t}_{\theta r, \theta} -c \Gammabar^{t}_{\phi r, \theta} -c \Gammabar^{t}_{\theta r, \phi} + a\Gammabar^{t}_{\phi r, \phi}  \right)\notag\\
& =:  -\frac{1}{u^{2}}\frac{(-|\gbar|)^{1/2}}{|g_{S}|^{3/2}} \cdot \Red (D) \Black \label{equ:highest order derivative term in straight out flow equation}
\end{align}

and 
\begin{align}
\CornflowerBlue (l2)\Black & := g_{S}^{\theta \theta}\left(\left(-\frac{1}{u^{2}} \right)_{,\theta}  \left(\frac{-|\gbar|}{|g_{S}|} \right)^{1/2} + \left(-\frac{1}{u^{2}}\right)\left(\frac{-|\gbar|}{|g_{S}|} \right)^{1/2}_{,\theta} \right)\Gammabar^{t}_{\theta r}\notag\\
& + g_{S}^{\theta \phi} \left( \left(-\frac{1}{u^{2}}\right)_{,\theta} \left(\frac{-|\gbar|}{|g_{S}|} \right)^{1/2} +  \left(-\frac{1}{u^{2}}\right) \left(\frac{-|\gbar|}{|g_{S}|} \right)^{1/2}_{,\theta} \right)\Gammabar^{t}_{\phi r}\notag\\
&  + g_{S}^{\phi \theta}\left(\left(-\frac{1}{u^{2}} \right)_{,\phi}  \left(\frac{-|\gbar|}{|g_{S}|} \right)^{1/2} + \left(-\frac{1}{u^{2}}\right)\left(\frac{-|\gbar|}{|g_{S}|} \right)^{1/2}_{,\phi} \right)\Gammabar^{t}_{\theta r}\notag\\
&  + g_{S}^{\phi \phi}\left(\left(-\frac{1}{u^{2}}\right)_{,\phi} \left(\frac{-|\gbar|}{|g_{S}|} \right)^{1/2} +  \left(-\frac{1}{u^{2}}\right) \left(\frac{-|\gbar|}{|g_{S}|} \right)^{1/2}_{,\phi} \right)\Gammabar^{t}_{\phi r}\notag\\
& = \CornflowerBlue\left(\left(-\frac{1}{u^{2}} \right)_{,\theta}  \left(\frac{-|\gbar|}{|g_{S}|} \right)^{1/2} + \left(-\frac{1}{u^{2}}\right)\left(\frac{-|\gbar|}{|g_{S}|} \right)^{1/2}_{,\theta} \right)\cdot (g_{S}^{\theta \theta}\Gammabar^{t}_{\theta r} + g_{S}^{\theta \phi}\Gammabar^{t}_{\phi r}) \Black\notag\\
& + \CornflowerBlue\left(\left(-\frac{1}{u^{2}}\right)_{,\phi} \left(\frac{-|\gbar|}{|g_{S}|} \right)^{1/2} +  \left(-\frac{1}{u^{2}}\right) \left(\frac{-|\gbar|}{|g_{S}|} \right)^{1/2}_{,\phi} \right)\cdot (g_{S}^{\phi \theta}\Gammabar^{t}_{\theta r} + g_{S}^{\phi \phi}\Gammabar^{t}_{\phi r})\Black
\end{align}

Notice that all the second derivatives of $\gbar$ in the divergence equation (\ref{equ:lower order term in straight out flow 1}) lie in $(D)$. Note that:
\begin{align*}
\Gammabar^{t}_{\theta r} & = \frac{1}{2}\Big(\gbar^{tt}(\gbar_{\theta t, r} + \gbar_{r t, \theta} - \gbar_{\theta r, t}) + \gbar^{tr}(\gbar_{\theta r, r} + \gbar_{rr, \theta} - g_{\theta r, r}) + \gbar^{t\theta}(\gbar_{\theta \theta, r} + \gbar_{r\theta, \theta}- \gbar_{\theta r, \theta})\\
& + \gbar^{t\phi}(\gbar_{\theta \phi, r} + \gbar_{r \phi, \theta} - g_{\theta r, \phi})\Big)\\
& = \frac{1}{2}\Big(\gbar^{tt}(\gbar_{\theta t, r} + \gbar_{rt, \theta}) + \gbar^{tr}\gbar_{rr, \theta} + \gbar^{t\theta}\gbar_{\theta \theta, r} + \gbar^{t\phi}\gbar_{\theta \phi, r} \Big)\\
& = \frac{1}{2|\gbar|}\Big(u^{2}(ab - c^{2})(e_{,r} + d_{,\theta}) -d(ab - c^{2})2uu_{,\theta} + u^{2}(cf - be)a_{,r} + u^{2}(ce - af)c_{,r} \Big)
\end{align*}

\begin{align*}
\Gammabar^{t}_{\phi r} & = \frac{1}{2}\Big(\gbar^{tt}(\gbar_{\phi t, r} + \gbar_{r t, \phi} - g_{\phi r, t}) + \gbar^{tr}(\gbar_{\phi r, r} + \gbar_{rr, \phi} - \gbar_{\phi r, r}) + \gbar^{t\theta}(\gbar_{\phi \theta, r} + \gbar_{r \theta, \phi} - \gbar_{\phi r, \theta}) \\
& + \gbar^{t\phi}(\gbar_{\phi \phi, r} + \gbar_{r \phi, \phi} - \gbar_{\phi r, \phi})  \Big)\\
& = \frac{1}{2}\Big(\gbar^{tt}(\gbar_{\phi t, r} + \gbar_{rt, \phi}) + \gbar^{tr}\gbar_{rr, \phi} + \gbar^{t\theta} \gbar_{\phi \theta, r} + \gbar^{t\phi}\gbar_{\phi \phi, r} \Big)\\
& = \frac{1}{2|\gbar|}\Big(u^{2}(ab -  c^{2})(f_{,r} + d_{,\phi}) -d(ab - c^{2})2uu_{,\phi} + u^{2}(cf - be)c_{,r} + u^{2}(ce - af)b_{,r} \Big)
\end{align*}

We compute $(D)$ as follows:
\begin{align}
(D) & = \frac{b}{2}[\gbar^{tt}(\gbar_{\theta t, r} + \gbar_{rt, \theta}) + \gbar^{tr}\gbar_{rr, \theta} + \gbar^{t\theta}\gbar_{\theta \theta, r} + \gbar^{t\phi}\gbar_{\theta \phi, r} ]_{,\theta} \notag\\
& - \frac{c}{2}[\gbar^{tt}(\gbar_{\phi t, r} + \gbar_{rt, \phi}) + \gbar^{tr}\gbar_{rr, \phi} + \gbar^{t\theta} \gbar_{\phi \theta, r} + \gbar^{t\phi}\gbar_{\phi \phi, r} ]_{,\theta}\notag\\
& - \frac{c}{2}[\gbar^{tt}(\gbar_{\theta t, r} + \gbar_{rt, \theta}) + \gbar^{tr}\gbar_{rr, \theta} + \gbar^{t\theta}\gbar_{\theta \theta, r} + \gbar^{t\phi}\gbar_{\theta \phi, r} ]_{,\phi}\notag\\
& + \frac{a}{2}[\gbar^{tt}(\gbar_{\phi t, r} + \gbar_{rt, \phi}) + \gbar^{tr}\gbar_{rr, \phi} + \gbar^{t\theta} \gbar_{\phi \theta, r} + \gbar^{t\phi}\gbar_{\phi \phi, r} ]_{,\phi}\notag\\
& = \Red\frac{1}{2}\gbar^{tt} (b\gbar_{\theta t, r\theta} + b\gbar_{rt, \theta \theta} -c\gbar_{\phi t, r\theta} -c\gbar_{rt, \phi \theta} -c \gbar_{\theta t, r\phi} - c\gbar_{rt, \theta \phi} + a\gbar_{\phi t, r\phi} + a\gbar_{rt, \phi \phi})\Black\notag\\
& + \Green\frac{1}{2}\gbar^{tr}(b\gbar_{rr, \theta \theta} - c\gbar_{rr, \phi \theta} - c\gbar_{rr, \theta \phi } + a\gbar_{rr, \phi \phi})\Black + \BlueViolet\frac{1}{2}\gbar^{t\theta}(b\gbar_{\theta \theta, r\theta} - c\gbar_{\phi \theta, r\theta} - c\gbar_{\theta \theta, r\phi} \Black\notag\\
& + a \gbar_{\phi \theta, r\phi}) + \YellowOrange\frac{1}{2}\gbar^{t\phi} (b\gbar_{\theta \phi, r\theta} - c\gbar_{\phi \phi, r\theta} -c\gbar_{\theta \phi, r\phi} +a\gbar_{\phi \phi, r\phi} )\Black + \CornflowerBlue (E) \Black \notag\\
& =: \Red(I)\Black + \Green(II)\Black + \BlueViolet(III)\Black + \YellowOrange(IV)\Black + \CornflowerBlue (E) \Black \label{equ:lower order term in straight out flow 3}
\end{align}

where $\CornflowerBlue (E) \Black$ consists of lower (first) order derivatives in $\gbar$:
\begin{align}
\CornflowerBlue (E) \Black & = \frac{b}{2}[ \gbar^{tt}_{,\theta}\cdot (\gbar_{\theta t, r} + \gbar_{rt, \theta}) + \gbar^{tr}_{,\theta}\cdot\gbar_{rr, \theta} + \gbar^{t\theta}_{,\theta}\cdot\gbar_{\theta \theta, r} + \gbar^{t\phi}_{,\theta}\cdot\gbar_{\theta \phi, r} ] \notag\\
& - \frac{c}{2}[\gbar^{tt}_{,\theta}\cdot(\gbar_{\phi t, r} + \gbar_{rt, \phi}) + \gbar^{tr}_{,\theta}\cdot\gbar_{rr, \phi} + \gbar^{t\theta}_{,\theta}\cdot \gbar_{\phi \theta, r} + \gbar^{t\phi}_{,\theta}\cdot\gbar_{\phi \phi, r} ]\notag\\
& - \frac{c}{2}[\gbar^{tt}_{,\phi}\cdot(\gbar_{\theta t, r} + \gbar_{rt, \theta}) + \gbar^{tr}_{,\phi}\cdot\gbar_{rr, \theta} + \gbar^{t\theta}_{,\phi}\cdot\gbar_{\theta \theta, r} + \gbar^{t\phi}_{,\phi}\cdot\gbar_{\theta \phi, r}]\notag\\
& + \frac{a}{2}[\gbar^{tt}_{,\phi}\cdot(\gbar_{\phi t, r} + \gbar_{rt, \phi}) + \gbar^{tr}_{,\phi}\cdot\gbar_{rr, \phi} + \gbar^{t\theta}_{,\phi} \cdot\gbar_{\phi \theta, r} + \gbar^{t\phi}_{,\phi}\cdot\gbar_{\phi \phi, r}]\notag\\
& = \frac{1}{2}\gbar^{tt}_{,\theta}(be_{,r} + bd_{,\theta} - cf_{,r} - cd_{,\phi})+ \frac{1}{2}\gbar^{tt}_{,\phi} (- ce_{,r} - cd_{,\theta} + af_{,r} + ad_{,\phi})\notag\\
& + \frac{1}{2}\gbar^{tr}_{,\theta}(2buu_{,\theta} - 2cuu_{,\phi}) + \frac{1}{2}\gbar^{tr}_{,\phi}(-2cuu_{,\theta} + 2auu_{,\phi}) + \frac{1}{2}\gbar^{t\theta}_{,\theta}(ba_{,r} - cc_{,r})\notag\\ 
& + \frac{1}{2}\gbar^{t\theta}_{,\phi}(-ca_{,r} + ac_{,r}) + \frac{1}{2}\gbar^{t\phi}_{,\theta}(bc_{,r} - cb_{,r}) + \frac{1}{2}\gbar^{t\phi}_{,\phi}(-cc_{,r} + ab_{,r})\label{equ:E term in the straight out flow}\tag{see (\ref{equ:lower order term in straight out flow 3})}
\end{align}

Therefore the fourth condition (\ref{def:fourth conditions of straight out coordinate chart}) becomes:
\begin{align}
0 = \divg_{g_{S}}\alpha^{\nu} & = \beta^{\theta}_{,\theta} + \beta^{\phi}_{,\phi} + \beta^{\theta}\cot{\theta} = \Red (B) \Black + \CornflowerBlue (l1) \Black \tag{see (\ref{equ:lower order term in straight out flow 1})}\notag\\
& = \Red (C) \Black + \CornflowerBlue (l2) \Black + \CornflowerBlue (l1) \Black \tag{see (\ref{equ:lower order term in straight out flow 2})} \notag\\
& =  -\frac{1}{u^{2}}\frac{(-|\gbar|)^{1/2}}{|g_{S}|^{3/2}} \cdot \Red (D) \Black + \CornflowerBlue (l2) \Black + \CornflowerBlue (l1) \Black  \tag{see (\ref{equ:highest order derivative term in straight out flow equation})} \notag\\
& =  -\frac{1}{u^{2}}\frac{(-|\gbar|)^{1/2}}{|g_{S}|^{3/2}} \cdot \Big(\Red(I)\Black + \Green(II)\Black + \BlueViolet(III)\Black + \YellowOrange(IV)\Black + \CornflowerBlue (E) \Black\Big) + \CornflowerBlue (l2) \Black + \CornflowerBlue (l1) \Black\label{equ:fourth condition in straight out flow broken down}
\end{align}

The highest derivatives of $\gbar$ (second derivatives) lie in $\Red(I)\Black + \Green(II)\Black + \BlueViolet(III)\Black + \YellowOrange(IV)$.

\subsection{Straight Out Flow Coordinate Chart: Complete Form}
\label{subsec:straight out flow equation complete form}

In this subsection, (\ref{equ:fourth condition in straight out flow broken down}) is computed explicitly from left to right.

\emph{\Red(I)\Black}:
\begin{align}
\Red(I)\Black & = \frac{u^{2}(ab - c^{2})}{2|\gbar|}(be_{,r \theta} +  bd_{,\theta \theta} - cf_{,r \theta} - cd_{,\phi \theta} - ce_{, r \phi} - cd_{,\theta \phi} + af_{,r\phi} + ad_{,\phi \phi} )\notag\\
& = \frac{u^{2}}{2}\frac{|g_{S}|}{|\gbar|} (\Red bd_{,\theta \theta} - 2cd_{,\theta \phi} + ad_{,\phi \phi} \Black + be_{,r \theta} - cf_{,r \theta} - ce_{, r \phi} + af_{,r\phi}). \label{equ:second order derivative in d in straight out flow equation}
\end{align}

Therefore
\begin{align}
& -\frac{1}{u^{2}}\frac{(-|\gbar|)^{1/2}}{|g_{S}|^{3/2}}\cdot\Red(I)\Black = -\frac{1}{\cancel{u^{2}}}\frac{(-|\gbar|)^{1/2}}{|g_{S}|^{3/2}}\cdot \frac{\cancel{u^{2}}}{2}\frac{|g_{S}|}{|\gbar|}\cdot \notag\\
& (\Red bd_{,\theta \theta} - 2cd_{,\theta \phi} + ad_{,\phi \phi} \Black + be_{,r \theta} - cf_{,r \theta} - ce_{, r \phi} + af_{,r\phi} ) \notag\\
& = \frac{1}{2}\left(\frac{-1}{|g_{S}||\gbar|} \right)^{1/2}\Big((\Red bd_{,\theta \theta} - 2cd_{,\theta \phi} + ad_{,\phi \phi}\Black) + (be_{,r \theta} - cf_{,r \theta} - ce_{, r \phi} + af_{,r\phi}) \Big).\label{equ:scaled I term}
\end{align}

\emph{\Green (II) \Black}:
\begin{align}
\Green (II) \Black & = \frac{- d (ab - c^{2})}{2|\gbar|}(b(u^{2})_{,\theta \theta} - 2c(u^{2})_{,\theta \phi} + a(u^{2})_{,\phi \phi}) \notag \\
& =  \Green\frac{- d|g_{S}|}{2|\gbar|}(b(u^{2})_{,\theta \theta} - 2c(u^{2})_{,\theta \phi} + a(u^{2})_{,\phi \phi}) \Black.
\end{align}
Therefore
\begin{align}
-\frac{1}{u^{2}}\frac{(-|\gbar|)^{1/2}}{|g_{S}|^{3/2}} \cdot\Green(II)\Black & = -\frac{1}{u^{2}}\frac{(-|\gbar|)^{1/2}}{|g_{S}|^{3/2}} \cdot  \Green\frac{- d|g_{S}|}{2|\gbar|}(b(u^{2})_{,\theta \theta} - 2c(u^{2})_{,\theta \phi} + a(u^{2})_{,\phi \phi}) \Black\notag\\
& = -\frac{1}{2}\left(\frac{-1}{|g_{S}||\gbar|} \right)^{1/2}\frac{d}{u^{2}}\Green(b(u^{2})_{,\theta \theta} - 2c(u^{2})_{,\theta \phi} + a(u^{2})_{,\phi \phi}) \Black.\label{equ:scaled II term}
\end{align}

\emph{$\BlueViolet (III) \Black$}:
\begin{align}
\BlueViolet (III) \Black & = \BlueViolet\frac{u^{2}(cf - be)}{2|\gbar|} (ba_{,r \theta} - cc_{,r \theta} - ca_{,r \phi} + ac_{, r\phi})\Black.
\end{align}
Therefore
\begin{align}
-\frac{1}{u^{2}}\frac{(-|\gbar|)^{1/2}}{|g_{S}|^{3/2}} \cdot \BlueViolet(III)\Black & = -\frac{1}{\cancel{u^{2}}}\frac{(-|\gbar|)^{1/2}}{|g_{S}|^{3/2}} \cdot\BlueViolet\frac{\cancel{u^{2}}(cf - be)}{2|\gbar|} (ba_{,r \theta} - cc_{,r \theta} - ca_{,r \phi} + ac_{, r\phi})\Black\notag\\
& = \frac{cf - be}{2}\frac{1}{|g_{S}|^{3/2}}\frac{1}{(-|\gbar|)^{1/2}}\BlueViolet(ba_{,r \theta} - cc_{,r \theta} - ca_{,r \phi} + ac_{, r\phi})\Black\notag\\
& = \frac{1}{2}\left(\frac{-1}{|g_{S}||\gbar|} \right)^{1/2}\frac{cf - be}{ab - c^{2}}\BlueViolet(ba_{,r \theta} - cc_{,r \theta} - ca_{,r \phi} + ac_{, r\phi})\Black. \label{equ:scaled III term}
\end{align}

\emph{$\YellowOrange (IV)$}:
\begin{align}
\YellowOrange (IV) \Black & = \YellowOrange\frac{u^{2}(ce - af)}{2|\gbar|} (bc_{,r\theta} - cb_{,r\theta} - cc_{,r\phi} + ab_{,r\phi})\Black.
\end{align}
Therefore
\begin{align}
-\frac{1}{u^{2}}\frac{(-|\gbar|)^{1/2}}{|g_{S}|^{3/2}} \cdot \YellowOrange(IV)\Black & =  -\frac{1}{u^{2}}\frac{(-|\gbar|)^{1/2}}{|g_{S}|^{3/2}} \cdot \YellowOrange\frac{u^{2}(ce - af)}{2|\gbar|} (bc_{,r\theta} - cb_{,r\theta} - cc_{,r\phi} + ab_{,r\phi})\Black\notag\\
& = \frac{ce - af}{2}\frac{1}{|g_{S}|^{3/2}}\frac{1}{(-|\gbar|)^{1/2}} \YellowOrange(bc_{,r\theta} - cb_{,r\theta} - cc_{,r\phi} + ab_{,r\phi})\Black\notag\\
& = \frac{1}{2}\left(\frac{-1}{|g_{S}||\gbar|}\right)^{1/2}\frac{ce - af}{ab - c^{2}}\YellowOrange(bc_{,r\theta} - cb_{,r\theta} - cc_{,r\phi} + ab_{,r\phi})\Black. \label{equ:scaled IV term}
\end{align}

We now compute $\CornflowerBlue (E)\Black$. First we compute the derivative of the inverse metric:
\begin{itemize}
\item[(1)]
\begin{align*}
\gbar^{tt}_{,\theta} & = \left(\frac{u^{2}|g_{S}|}{|\gbar|}\right)_{,\theta} = \frac{(2uu_{,\theta}|g_{S}| + u^{2}|g_{S}|_{,\theta})|\gbar| - u^{2}|g_{S}||\gbar|_{,\theta}}{|\gbar|^{2}}\\
& = \frac{2uu_{,\theta}|g_{S}| + 2u^{2}|g_{S}|\cot{\theta}}{|\gbar|} - \frac{u^{2}|g_{S}||\gbar|_{,\theta}}{|\gbar|^{2}}\tag{By (\ref{equ:derivative of det of gs})}\\
& = \frac{u^{2}|g_{S}|}{|\gbar|}\left(2\frac{u_{,\theta}}{u} + 2\cot{\theta} - \frac{|\gbar|_{,\theta}}{|\gbar|}\right). 
\end{align*}

Similarly
\begin{align*}
\gbar^{tt}_{,\phi} & = \left(\frac{u^{2}|g_{S}|}{|\gbar|}\right)_{,\phi} = \frac{(2uu_{,\phi}|g_{S}| + u^{2}\cancelto{0}{|g_{S}|_{,\phi}})|\gbar| - u^{2}|g_{S}||\gbar|_{,\phi}}{|\gbar|^{2}}\tag{$|g_{S}| = r^{4}\sin^{2}{\theta}$}\\
& = \frac{2uu_{,\phi}|g_{S}|}{|\gbar|} - \frac{u^{2}|g_{S}||\gbar|_{,\phi}}{|\gbar|^{2}}\\
& = \frac{u^{2}|g_{S}|}{|\gbar|}\left(2\frac{u_{,\phi}}{u} - \frac{|\gbar|_{,\phi}}{|\gbar|}\right). 
\end{align*}

\item[(2)]
\begin{align*}
\gbar^{tr}_{,\theta} & = \left(\frac{-d|g_{S}|}{|\gbar|} \right)_{,\theta} = \frac{(-d_{,\theta}|g_{S}| - d|g_{S}|_{,\theta})|\gbar| - (-d|g_{S}|)|\gbar|_{,\theta}}{|\gbar|^{2}}\\
& = \frac{-d_{,\theta}|g_{S}| - 2d|g_{S}|\cot{\theta}}{|\gbar|} + \frac{d|g_{S}||\gbar|_{,\theta}}{|\gbar|^{2}} \tag{By (\ref{equ:derivative of det of gs})}\\
& = \frac{-d|g_{S}|}{|\gbar|}\left(\frac{d_{,\theta}}{d} + 2\cot{\theta} - \frac{|\gbar|_{,\theta}}{|\gbar|} \right).
\end{align*}

Similarly
\begin{align*}
\gbar^{tr}_{,\phi} & = \left(\frac{-d|g_{S}|}{|\gbar|} \right)_{,\phi} = \frac{(-d_{,\phi}|g_{S}| - d\cancelto{0}{|g_{S}|_{,\phi}})|\gbar| - (-d|g_{S}|)|\gbar|_{,\phi}}{|\gbar|^{2}}\\
& = \frac{-d_{,\phi}|g_{S}|}{|\gbar|} + \frac{d|g_{S}||\gbar|_{,\phi}}{|\gbar|^{2}}\\
& = \frac{-d|g_{S}|}{|\gbar|}\left(\frac{d_{,\phi}}{d} - \frac{|\gbar|_{,\phi}}{|\gbar|} \right).
\end{align*}

\item[(3)]
\begin{align*}
\gbar^{t\theta}_{,\theta} & = \left(\frac{u^{2}(cf - be)}{|\gbar|} \right)_{,\theta} = \frac{\big(2uu_{,\theta}(cf - be) + u^{2}(cf -be)_{,\theta}\big)|\gbar| - u^{2}(cf - be)|\gbar|_{,\theta}}{|\gbar|^{2}}\\
& = \frac{2uu_{,\theta}(cf - be) + u^{2}(cf -be)_{,\theta}}{|\gbar|} - \frac{u^{2}(cf - be)|\gbar|_{,\theta}}{|\gbar|^{2}}\\
& = \frac{u^{2}(cf - be)}{|\gbar|}\left(2\frac{u_{,\theta}}{u} + \frac{(cf -be)_{,\theta}}{cf - be} - \frac{|\gbar|_{,\theta}}{|\gbar|} \right).
\end{align*}

Similarly
\begin{align*}
\gbar^{t\theta}_{,\phi} & = \left(\frac{u^{2}(cf - be)}{|\gbar|} \right)_{,\phi} = \frac{\big(2uu_{,\phi}(cf - be) + u^{2}(cf -be)_{,\phi}  \big) |\gbar| - u^{2}(cf - be)|\gbar|_{,\phi} }{|\gbar|^{2}}\\
& = \frac{2uu_{,\phi}(cf - be) + u^{2}(cf -be)_{,\phi}}{|\gbar|} - \frac{u^{2}(cf - be)|\gbar|_{,\phi}}{|\gbar|^{2}}\\
& = \frac{u^{2}(cf - be)}{|\gbar|}\left(2\frac{u_{,\phi}}{u} + \frac{(cf -be)_{,\phi}}{cf - be} -\frac{|\gbar|_{,\phi}}{|\gbar|} \right).
\end{align*}

\item[(4)]
\begin{align*}
\gbar^{t\phi}_{,\theta} & = \left(\frac{u^{2}(ce - af)}{|\gbar|} \right)_{,\theta} = \frac{\big(2uu_{,\theta}(ce - af) + u^{2}(ce - af)_{,\theta} \big)|\gbar| - u^{2}(ce - af)|\gbar|_{,\theta}}{|\gbar|^{2}}\\
& = \frac{2uu_{,\theta}(ce - af) + u^{2}(ce - af)_{,\theta}}{|\gbar|} - \frac{u^{2}(ce - af)|\gbar|_{,\theta}}{|\gbar|^{2}}\\
& = \frac{u^{2}(ce - af)}{|\gbar|}\left(2\frac{u_{,\theta}}{u} + \frac{(ce - af)_{,\theta}}{ce - af} - \frac{|\gbar|_{,\theta}}{|\gbar|} \right).
\end{align*}

Similarly
\begin{align*}
\gbar^{t\phi}_{,\phi} & = \left(\frac{u^{2}(ce - af)}{|\gbar|}\right)_{,\phi} = \frac{\big(2uu_{,\phi}(ce - af) + u^{2}(ce - af)_{,\phi} \big)|\gbar| - u^{2}(ce - af)|\gbar|_{,\phi}}{|\gbar|^{2}}\\
& = \frac{2uu_{,\phi}(ce - af) + u^{2}(ce - af)_{,\phi}}{|\gbar|} - \frac{u^{2}(ce - af)|\gbar|_{,\phi}}{|\gbar|^{2}}\\
& = \frac{u^{2}(ce - af)}{|\gbar|}\left(2\frac{u_{,\phi}}{u} + \frac{(ce - af)_{,\phi}}{ce - af} - \frac{|\gbar|_{,\phi}}{|\gbar|}\right).
\end{align*}

\end{itemize}

Therefore 
\begin{align*}
\CornflowerBlue (E) \Black & = \frac{1}{2}\frac{u^{2}|g_{S}|}{|\gbar|}\left(2\frac{u_{,\theta}}{u} + 2\cot{\theta} - \frac{|\gbar|_{,\theta}}{|\gbar|}\right)(be_{,r} + bd_{,\theta} - cf_{,r} - cd_{,\phi})\\
& + \frac{1}{2}\frac{u^{2}|g_{S}|}{|\gbar|}\left(2\frac{u_{,\phi}}{u} - \frac{|\gbar|_{,\phi}}{|\gbar|}\right) (- ce_{,r} - cd_{,\theta} + af_{,r} + ad_{,\phi})\\
& + \frac{-du|g_{S}|}{|\gbar|}\left(\frac{d_{,\theta}}{d} + 2\cot{\theta} - \frac{|\gbar|_{,\theta}}{|\gbar|} \right)(bu_{,\theta} - cu_{,\phi}) + \frac{-du|g_{S}|}{|\gbar|}\left(\frac{d_{,\phi}}{d} - \frac{|\gbar|_{,\phi}}{|\gbar|} \right)(-cu_{,\theta} \\
& + au_{,\phi}) + \frac{1}{2}\frac{u^{2}(cf - be)}{|\gbar|}\left(2\frac{u_{,\theta}}{u} + \frac{(cf -be)_{,\theta}}{cf - be} - \frac{|\gbar|_{,\theta}}{|\gbar|} \right)(ba_{,r} - cc_{,r}) \\
& + \frac{1}{2}\frac{u^{2}(cf - be)}{|\gbar|}\left(2\frac{u_{,\phi}}{u} + \frac{(cf-be)_{,\phi}}{cf - be} -\frac{|\gbar|_{,\phi}}{|\gbar|} \right)(-ca_{,r} + ac_{,r}) \\
& + \frac{1}{2}\frac{u^{2}(ce - af)}{|\gbar|}\left(2\frac{u_{,\theta}}{u} + \frac{(ce - af)_{,\theta}}{ce - af} - \frac{|\gbar|_{,\theta}}{|\gbar|} \right)(bc_{,r} - cb_{,r}) \\
& + \frac{1}{2}\frac{u^{2}(ce - af)}{|\gbar|}\left(2\frac{u_{,\phi}}{u} + \frac{(ce - af)_{,\phi}}{ce - af} - \frac{|\gbar|_{,\phi}}{|\gbar|}\right)(-cc_{,r} + ab_{,r}).
\end{align*}

Therefore
\begin{align*}
& -\frac{1}{u^{2}}\frac{(-|\gbar|)^{1/2}}{|g_{S}|^{3/2}}\CornflowerBlue (E) \Black = \frac{1}{2}\left(\frac{-1}{|g_{S}||\gbar|} \right)^{1/2}\left[\left(2\frac{u_{,\theta}}{u} + 2\cot{\theta} - \frac{|\gbar|_{,\theta}}{|\gbar|}\right)(be_{,r} + bd_{,\theta} - cf_{,r} - cd_{,\phi})\right.\notag\\
& \left. + \left(2\frac{u_{,\phi}}{u} - \frac{|\gbar|_{,\phi}}{|\gbar|}\right) (- ce_{,r} - cd_{,\theta} + af_{,r} + ad_{,\phi})\right]\notag\\
& -\left(\frac{-1}{|g_{S}||\gbar|}\right)^{1/2}\frac{d}{u}\left[\left(\frac{d_{,\theta}}{d} + 2\cot{\theta} - \frac{|\gbar|_{,\theta}}{|\gbar|} \right)(bu_{,\theta} - cu_{,\phi}) + \left(\frac{d_{,\phi}}{d} - \frac{|\gbar|_{,\phi}}{|\gbar|} \right)(-cu_{,\theta} + au_{,\phi})\right] \notag\\
& +\frac{1}{2}\left(\frac{-1}{|g_{S}||\gbar|}\right)^{1/2}\frac{cf - be}{ab - c^{2}}\left[\left(2\frac{u_{,\theta}}{u} + \frac{(cf -be)_{,\theta}}{cf - be} - \frac{|\gbar|_{,\theta}}{|\gbar|} \right)(ba_{,r} - cc_{,r})\right. \notag\\
& \left. + \left(2\frac{u_{,\phi}}{u} + \frac{(cf-be)_{,\phi}}{cf - be} -\frac{|\gbar|_{,\phi}}{|\gbar|} \right)(-ca_{,r} + ac_{,r})\right] \notag\\
& +\frac{1}{2}\left(\frac{-1}{|g_{S}||\gbar|}\right)^{1/2}\frac{ce - af}{ab - c^{2}}\left[\left(2\frac{u_{,\theta}}{u} + \frac{(ce - af)_{,\theta}}{ce - af} - \frac{|\gbar|_{,\theta}}{|\gbar|} \right)(bc_{,r} - cb_{,r})\right. \notag\\
& \left.+ \left(2\frac{u_{,\phi}}{u} + \frac{(ce - af)_{,\phi}}{ce - af} - \frac{|\gbar|_{,\phi}}{|\gbar|}\right)(-cc_{,r} + ab_{,r})\right].\notag
\end{align*}

That is:
\begin{align}
& -\frac{1}{u^{2}}\frac{(-|\gbar|)^{1/2}}{|g_{S}|^{3/2}}\CornflowerBlue (E) \Black = \frac{1}{2}\left(\frac{-1}{|g_{S}||\gbar|} \right)^{1/2}\Red\left\{\Black\left(2\frac{u_{,\theta}}{u} + 2\cot{\theta} - \frac{|\gbar|_{,\theta}}{|\gbar|}\right)(be_{,r} + bd_{,\theta} - cf_{,r} - cd_{,\phi})\right.\notag\\
& + \left(2\frac{u_{,\phi}}{u} - \frac{|\gbar|_{,\phi}}{|\gbar|}\right) (- ce_{,r} - cd_{,\theta} + af_{,r} + ad_{,\phi}) - \frac{2d}{u}\Green\left[\Black\left(\frac{d_{,\theta}}{d} + 2\cot{\theta} - \frac{|\gbar|_{,\theta}}{|\gbar|} \right)(bu_{,\theta} - cu_{,\phi})\right. \notag\\
& \left.+ \left(\frac{d_{,\phi}}{d} - \frac{|\gbar|_{,\phi}}{|\gbar|} \right)(-cu_{,\theta} + au_{,\phi})\Green\right]\Black +\frac{cf - be}{ab - c^{2}}\Green\left[\Black\left(2\frac{u_{,\theta}}{u} + \frac{(cf -be)_{,\theta}}{cf - be} - \frac{|\gbar|_{,\theta}}{|\gbar|} \right)(ba_{,r} - cc_{,r})\right. \notag\\
& \left. + \left(2\frac{u_{,\phi}}{u} + \frac{(cf-be)_{,\phi}}{cf - be} -\frac{|\gbar|_{,\phi}}{|\gbar|} \right)(-ca_{,r} + ac_{,r})\Green\right]\Black + \frac{ce - af}{ab - c^{2}}\Green\left[\Black\left(2\frac{u_{,\theta}}{u} + \frac{(ce - af)_{,\theta}}{ce - af} - \frac{|\gbar|_{,\theta}}{|\gbar|} \right)\cdot\right. \notag\\
& \left.\left. (bc_{,r} - cb_{,r}) + \left(2\frac{u_{,\phi}}{u} + \frac{(ce - af)_{,\phi}}{ce - af} - \frac{|\gbar|_{,\phi}}{|\gbar|}\right)(-cc_{,r} + ab_{,r})\Green\right]\Black\Red\right\}.\Black\label{equ:scaled E term}
\end{align}

Next we compute $\CornflowerBlue (l2) \Black$. To do that, we first have
\begin{align*}
& \left(-\frac{1}{u^{2}}\right)_{,\theta} \left(\frac{-|\gbar|}{|g_{S}|} \right)^{1/2} +  \left(-\frac{1}{u^{2}}\right) \left(\frac{-|\gbar|}{|g_{S}|} \right)^{1/2}_{,\theta}\\
& = \frac{2}{u^{3}}u_{,\theta} \left(\frac{-|\gbar|}{|g_{S}|} \right)^{1/2} - \frac{1}{u^{2}}\frac{1}{2}\left(\frac{-|\gbar|}{|g_{S}|} \right)^{-1/2}\left(\frac{-|\gbar|_{,\theta}|g_{S}| + |\gbar||g_{S}|_{,\theta}}{|g_{S}|^{2}}\right)\\
& = \frac{2}{u^{3}}u_{,\theta} \left(\frac{-|\gbar|}{|g_{S}|} \right)^{1/2} - \frac{1}{2u^{2}}\left(\frac{|g_{S}|}{-|\gbar|}\right)^{1/2}\left(\frac{-|\gbar|_{,\theta}}{|g_{S}|} + \frac{2|\gbar||g_{S}|\cot{\theta}}{|g_{S}|^{2}} \right)\tag{By (\ref{equ:derivative of det of gs})}\\
& =  \frac{2}{u^{3}}u_{,\theta} \left(\frac{-|\gbar|}{|g_{S}|} \right)^{1/2} - \frac{1}{2u^{2}}\left(\frac{-1}{|g_{S}||\gbar|}\right)^{1/2}\left(-|\gbar|_{,\theta} + 2|\gbar|\cot{\theta}\right)\\
& = \frac{1}{u^{2}}\left(\frac{-|\gbar|}{|g_{S}|}\right)^{1/2}\left(\frac{2u_{,\theta}}{u} - \frac{1}{2}\left(\frac{-1}{|\gbar|} \right)(-|\gbar|_{,\theta} + 2|\gbar|\cot{\theta}) \right)\\
& = \frac{1}{u^{2}}\left(\frac{-|\gbar|}{|g_{S}|}\right)^{1/2}\left(\frac{2u_{,\theta}}{u} + \cot{\theta} - \frac{|\gbar|_{,\theta}}{2|\gbar|}\right).
\end{align*}

Similarly we have:
\begin{align*}
& \left(-\frac{1}{u^{2}}\right)_{,\phi} \left(\frac{-|\gbar|}{|g_{S}|} \right)^{1/2} +  \left(-\frac{1}{u^{2}}\right) \left(\frac{-|\gbar|}{|g_{S}|} \right)^{1/2}_{,\phi}\\
& = \frac{2}{u^{3}}u_{,\phi} \left(\frac{-|\gbar|}{|g_{S}|} \right)^{1/2} - \frac{1}{u^{2}}\frac{1}{2}\left(\frac{-|\gbar|}{|g_{S}|} \right)^{-1/2}\left(\frac{-|\gbar|_{,\phi}|g_{S}| + |\gbar|\cancelto{0}{|g_{S}|_{,\phi}}}{|g_{S}|^{2}}\right)\\
& = \frac{2}{u^{3}}u_{,\phi} \left(\frac{-|\gbar|}{|g_{S}|} \right)^{1/2} - \frac{1}{2u^{2}}\left(\frac{-|g_{S}|}{|\gbar|} \right)^{1/2}\left(\frac{-|\gbar|_{,\phi}}{|g_{S}|}\right)\\
& = \frac{2}{u^{3}}u_{,\phi} \left(\frac{-|\gbar|}{|g_{S}|} \right)^{1/2} - \frac{1}{2u^{2}}\left(\frac{-1}{|g_{S}||\gbar|} \right)^{1/2}(-|\gbar|_{,\phi})\\
& = \frac{1}{u^{2}}\left(\frac{-|\gbar|}{|g_{S}|} \right)^{1/2}\left(\frac{2u_{,\phi}}{u} - \frac{1}{2}\left(\frac{-1}{|\gbar|} \right)(-|\gbar|_{,\phi})  \right)\\
& = \frac{1}{u^{2}}\left(\frac{-|\gbar|}{|g_{S}|} \right)^{1/2}\left(\frac{2u_{,\phi}}{u} - \frac{|\gbar|_{,\phi}}{2|\gbar|}  \right).
\end{align*}

Then we need to compute:
\begin{align*}
& g_{S}^{\theta \theta}\Gammabar^{t}_{\theta r} + g_{S}^{\theta \phi}\Gammabar^{t}_{\phi r} = \frac{1}{|g_{S}|}(b\Gammabar^{t}_{\theta r} - c\Gammabar^{t}_{\phi r})\\
& = \frac{1}{2}\frac{1}{|g_{S}||\gbar|}\Big(b \Red\big(\Black u^{2}(ab - c^{2})(e_{,r} + d_{,\theta}) -d(ab - c^{2})2uu_{,\theta} + u^{2}(cf - be)a_{,r} + u^{2}(ce - af)c_{,r} \Red\big)\Black\\
& - c\Red\big(\Black u^{2}(ab -  c^{2})(f_{,r} + d_{,\phi}) -d(ab - c^{2})2uu_{,\phi} + u^{2}(cf - be)c_{,r} + u^{2}(ce - af)b_{,r} \Red\big)\Black\Big)\\
& = \frac{1}{2}\frac{1}{|g_{S}||\gbar|}\Big(u^{2}|g_{S}|(be_{,r} + bd_{,\theta} - cf_{,r} - cd_{,\phi}) -d|g_{S}|(2buu_{,\theta} - 2cuu_{,\phi}) \\
& + u^{2}(cf - be)(ba_{,r} - cc_{,r}) + u^{2}(ce -af)(bc_{,r} - cb_{,r})  \Big)\\
& = \frac{u^{2}}{2}\frac{1}{|\gbar|}\Big((be_{,r} + bd_{,\theta} - cf_{,r} - cd_{,\phi}) -\frac{2d}{u}(bu_{,\theta} - cu_{,\phi}) + \frac{cf - be}{ab - c^{2}}(ba_{,r} - cc_{,r}) \\
& + \frac{ce -af}{ab - c^{2}}(bc_{,r} - cb_{,r})  \Big).
\end{align*}

Similarly:
\begin{align*}
& g_{S}^{\phi \theta}\Gammabar^{t}_{\theta r} + g_{S}^{\phi \phi}\Gammabar^{t}_{\phi r} = \frac{1}{|g_{S}|}(-c\Gammabar^{t}_{\theta r} + a\Gammabar^{t}_{\phi r})\\
& = \frac{1}{2}\frac{1}{|g_{S}||\gbar|}\Big(-c\Red\big(\Black u^{2}(ab - c^{2})(e_{,r} + d_{,\theta}) -d(ab - c^{2})2uu_{,\theta} + u^{2}(cf - be)a_{,r} + u^{2}(ce - af) \Black\\
& \cdot c_{,r} \Red\big) + a \Red\big(\Black  u^{2}(ab -  c^{2})(f_{,r} + d_{,\phi}) -d(ab - c^{2})2uu_{,\phi} + u^{2}(cf - be)c_{,r} + u^{2}(ce - af)b_{,r}  \Red\big)\Black  \Big)\\
& = \frac{1}{2}\frac{1}{|g_{S}||\gbar|}\Big(u^{2}|g_{S}|(-ce_{,r} - cd_{,\theta} + af_{,r} + ad_{,\phi}) - d|g_{S}|(-2cuu_{,\theta} + 2auu_{,\phi})\\
& + u^{2}(cf - be)(-ca_{,r} + ac_{,r}) + u^{2}(ce - af)(-cc_{,r} + ab_{,r}) \Big)\\
& = \frac{u^{2}}{2}\frac{1}{|\gbar|}\Big((-ce_{,r} - cd_{,\theta} + af_{,r} + ad_{,\phi}) - \frac{2d}{u}(-cu_{,\theta} + au_{,\phi})+ \frac{cf - be}{ab - c^{2}}(-ca_{,r} + ac_{,r})\\
& + \frac{ce - af}{ab - c^{2}}(-cc_{,r} + ab_{,r}) \Big).
\end{align*}

Putting them all together, we get:
\begin{align*}
\CornflowerBlue (l2) \Black & = \frac{1}{\cancel{u^{2}}}\left(\frac{-|\gbar|}{|g_{S}|}\right)^{1/2}\left(\frac{2u_{,\theta}}{u} + \cot{\theta} - \frac{|\gbar|_{,\theta}}{2|\gbar|}\right)\cdot \frac{\cancel{u^{2}}}{2}\frac{1}{|\gbar|}\Big((be_{,r} + bd_{,\theta} - cf_{,r} - cd_{,\phi})\\
&  -\frac{2d}{u}(bu_{,\theta} - cu_{,\phi}) + \frac{cf - be}{ab - c^{2}}(ba_{,r} - cc_{,r}) + \frac{ce -af}{ab - c^{2}}(bc_{,r} - cb_{,r})  \Big)\\
& + \frac{1}{\bcancel{u^{2}}}\left(\frac{-|\gbar|}{|g_{S}|} \right)^{1/2}\left(\frac{2u_{,\phi}}{u} - \frac{|\gbar|_{,\phi}}{2|\gbar|}  \right)\cdot \frac{\bcancel{u^{2}}}{2}\frac{1}{|\gbar|}\Big((-ce_{,r} - cd_{,\theta} + af_{,r} + ad_{,\phi})\\
& - \frac{2d}{u}(-cu_{,\theta} + au_{,\phi})+ \frac{cf - be}{ab - c^{2}}(-ca_{,r} + ac_{,r}) + \frac{ce - af}{ab - c^{2}}(-cc_{,r} + ab_{,r}) \Big).
\end{align*}

Simplifying $\CornflowerBlue (l2) \Black$, we get:
\begin{align}
\CornflowerBlue (l2) \Black & = -\frac{1}{2}\left(\frac{-1}{|g_{S}||\gbar|}\right)^{1/2}\left\{\left(\frac{2u_{,\theta}}{u} + \cot{\theta} - \frac{|\gbar|_{,\theta}}{2|\gbar|}\right)\cdot \Big((be_{,r} + bd_{,\theta} - cf_{,r} - cd_{,\phi})\right. \notag\\
& - \frac{2d}{u}(bu_{,\theta} - cu_{,\phi}) + \frac{cf - be}{ab - c^{2}}(ba_{,r} - cc_{,r}) + \frac{ce -af}{ab - c^{2}}(bc_{,r} - cb_{,r})  \Big)\notag\\    
& + \left(\frac{2u_{,\phi}}{u} - \frac{|\gbar|_{,\phi}}{2|\gbar|}  \right)\cdot \Big((-ce_{,r} - cd_{,\theta} + af_{,r} + ad_{,\phi}) - \frac{2d}{u}(-cu_{,\theta} + au_{,\phi})\notag\\
& + \left. \frac{cf - be}{ab - c^{2}}(-ca_{,r} + ac_{,r}) + \frac{ce - af}{ab - c^{2}}(-cc_{,r} + ab_{,r}) \Big)\right\}. \label{equ:l2 term simplified}
\end{align}

Next we compute $\CornflowerBlue (l1)\Black$:
\begin{align*}
& \CornflowerBlue (l1)\Black = (g_{S, \theta}^{\theta \theta}\alpha_{\theta} + g_{S,\theta}^{\theta \phi}\alpha_{\phi} + g_{S, \phi}^{\phi \theta}\alpha_{\theta} + g_{S, \phi}^{\phi \phi} \alpha_{\phi}) + (g_{S}^{\theta \theta} \alpha_{\theta} + g_{S}^{\theta \phi} \alpha_{\phi}) \cot{\theta} \\
& = (g_{S,\theta}^{\theta \theta} + g_{S,\phi}^{\phi \theta} + g_{S}^{\theta \theta}\cot{\theta})\alpha_{\theta} + (g_{S,\theta}^{\theta \phi} + g_{S,\phi}^{\phi \phi} + g_{S}^{\theta \phi}\cot{\theta})\alpha_{\phi}\\
& =  -\frac{1}{u^{2}} \left(\frac{-|\gbar|}{|g_{S}|} \right)^{1/2} (g_{S,\theta}^{\theta \theta} + g_{S,\phi}^{\phi \theta} + g_{S}^{\theta \theta}\cot{\theta})\Gammabar^{t}_{\theta r} \\
& + \left(-\frac{1}{u^{2}}\right) \left(\frac{-|\gbar|}{|g_{S}|} \right)^{1/2} (g_{S,\theta}^{\theta \phi} + g_{S,\phi}^{\phi \phi} + g_{S}^{\theta \phi}\cot{\theta})\Gammabar^{t}_{\phi r} \\
& =  -\frac{1}{u^{2}} \left(\frac{-|\gbar|}{|g_{S}|} \right)^{1/2} \left((g_{S,\theta}^{\theta \theta} + g_{S,\phi}^{\phi \theta} + g_{S}^{\theta \theta}\cot{\theta})\Gammabar^{t}_{\theta r} + (g_{S,\theta}^{\theta \phi} + g_{S,\phi}^{\phi \phi} + g_{S}^{\theta \phi}\cot{\theta})\Gammabar^{t}_{\phi r}\right)\\
& =:  -\frac{1}{u^{2}} \left(\frac{-|\gbar|}{|g_{S}|} \right)^{1/2} \Big((l1a) + (l1b)\Big),
\end{align*}

where $(l1a): = (g_{S,\theta}^{\theta \theta} + g_{S,\phi}^{\phi \theta} + g_{S}^{\theta \theta}\cot{\theta})\Gammabar^{t}_{\theta r}$, and $(l1b): = (g_{S,\theta}^{\theta \phi} + g_{S,\phi}^{\phi \phi} + g_{S}^{\theta \phi}\cot{\theta})\Gammabar^{t}_{\phi r}$.

Note that
\begin{align}
(l1a) & =  (g_{S,\theta}^{\theta \theta} + g_{S,\phi}^{\phi \theta} + g_{S}^{\theta \theta}\cot{\theta}) \frac{1}{2|\gbar|}\Big(u^{2}(ab - c^{2})(e_{,r} + d_{,\theta}) -d(ab - c^{2})2uu_{,\theta}\notag\\
& + u^{2}(cf - be)a_{,r} + u^{2}(ce - af)c_{,r} \Big)\notag\\
& = \Peach \frac{u^{2}}{2}\frac{|g_{S}|}{|\gbar|}(g_{S,\theta}^{\theta \theta} + g_{S,\phi}^{\phi \theta} + g_{S}^{\theta \theta}\cot{\theta})d_{,\theta}\Black + \frac{1}{2|\gbar|}\Green\Big(\Black u^{2}|g_{S}|(g_{S,\theta}^{\theta \theta}+ g_{S,\phi}^{\phi \theta} + g_{S}^{\theta \theta}\cot{\theta})e_{,r}\notag \\
& - 2d|g_{S}|uu_{,\theta}(g_{S,\theta}^{\theta \theta} + g_{S,\phi}^{\phi \theta} + g_{S}^{\theta \theta}\cot{\theta}) + u^{2}(cf - be)(g_{S,\theta}^{\theta \theta} + g_{S,\phi}^{\phi \theta} + g_{S}^{\theta \theta}\cot{\theta})a_{,r} \notag\\
& + u^{2}(ce - af)(g_{S,\theta}^{\theta \theta} + g_{S,\phi}^{\phi \theta} + g_{S}^{\theta \theta}\cot{\theta})c_{,r} \Green\Big).\Black \notag
\end{align}

Similarly:
\begin{align}
(l1b) & = (g_{S,\theta}^{\theta \phi} + g_{S,\phi}^{\phi \phi} + g_{S}^{\theta \phi}\cot{\theta})  \frac{1}{2|\gbar|}\Big(u^{2}(ab -  c^{2})(f_{,r} + d_{,\phi}) -d(ab - c^{2})2uu_{,\phi} \notag\\
& + u^{2}(cf - be)c_{,r} + u^{2}(ce - af)b_{,r} \Big) \notag\\
 & = \Peach \frac{u^{2}}{2}\frac{|g_{S}|}{|\gbar|}(g_{S,\theta}^{\theta \phi} + g_{S,\phi}^{\phi \phi} + g_{S}^{\theta \phi}\cot{\theta})d_{,\phi}  \Black + \frac{1}{2|\gbar|}\Green\Big(\Black u^{2}|g_{S}|(g_{S,\theta}^{\theta \phi} + g_{S,\phi}^{\phi \phi} + g_{S}^{\theta \phi}\cot{\theta})f_{,r}\notag\\
&  - 2d|g_{S}|uu_{,\phi} (g_{S,\theta}^{\theta \phi} + g_{S,\phi}^{\phi \phi}+ g_{S}^{\theta \phi}\cot{\theta}) + u^{2}(cf - be)(g_{S,\theta}^{\theta \phi} + g_{S,\phi}^{\phi \phi} + g_{S}^{\theta \phi}\cot{\theta})c_{,r} \notag\\
& + u^{2}(ce - af)(g_{S,\theta}^{\theta \phi} + g_{S,\phi}^{\phi \phi} + g_{S}^{\theta \phi}\cot{\theta})b_{,r}\Green\Big).\Black\notag
\end{align}

Plug them back into $\CornflowerBlue(l1)\Black$, we get:
\begin{align}
& \CornflowerBlue(l1) \Black =  -\frac{1}{u^{2}} \left(\frac{-|\gbar|}{|g_{S}|} \right)^{1/2} \Peach \frac{u^{2}}{2}\frac{|g_{S}|}{|\gbar|}\Big((g_{S,\theta}^{\theta \theta} + g_{S,\phi}^{\phi \theta} + g_{S}^{\theta \theta}\cot{\theta})d_{,\theta}\notag\\
&  +  (g_{S,\theta}^{\theta \phi} + g_{S,\phi}^{\phi \phi} + g_{S}^{\theta \phi}\cot{\theta})d_{,\phi}\Big) \Black\notag\\
& + \left(-\frac{1}{u^{2}}\right) \left(\frac{-|\gbar|}{|g_{S}|} \right)^{1/2} \bigg\{\frac{u^{2}}{2}\frac{|g_{S}|}{|\gbar|}\Big((g_{S,\theta}^{\theta \theta}+ g_{S,\phi}^{\phi \theta} + g_{S}^{\theta \theta}\cot{\theta})e_{,r} \notag\\
& + (g_{S,\theta}^{\theta \phi} + g_{S,\phi}^{\phi \phi} + g_{S}^{\theta \phi}\cot{\theta})f_{,r}  \Big)\notag\\
& + \frac{(-2du)}{2}\frac{|g_{S}|}{|\gbar|}\Big((g_{S,\theta}^{\theta \theta} + g_{S,\phi}^{\phi \theta} + g_{S}^{\theta \theta}\cot{\theta})u_{,\theta} +  (g_{S,\theta}^{\theta \phi} + g_{S,\phi}^{\phi \phi}+ g_{S}^{\theta \phi}\cot{\theta})u_{,\phi} \Big)\notag\\
& +\frac{u^{2}(cf - be)}{2|\gbar|}\Big((g_{S,\theta}^{\theta \theta} + g_{S,\phi}^{\phi \theta} + g_{S}^{\theta \theta}\cot{\theta})a_{,r}  + (g_{S,\theta}^{\theta \phi} + g_{S,\phi}^{\phi \phi} + g_{S}^{\theta \phi}\cot{\theta})c_{,r} \Big)\notag\\
& +\frac{u^{2}(ce - af)}{2|\gbar|}\Big((g_{S,\theta}^{\theta \theta} + g_{S,\phi}^{\phi \theta} + g_{S}^{\theta \theta}\cot{\theta})c_{,r} + (g_{S,\theta}^{\theta \phi} + g_{S,\phi}^{\phi \phi} + g_{S}^{\theta \phi}\cot{\theta})b_{,r} \Big)\bigg\}.\label{equ:l1 term step 1}
\end{align}
And
\begin{align*}
(\ref{equ:l1 term step 1}) & =  \Peach\frac{1}{2}\left(\frac{|g_{S}|}{-|\gbar|}\right)^{1/2} \Big((g_{S,\theta}^{\theta \theta} + g_{S,\phi}^{\phi \theta} + g_{S}^{\theta \theta}\cot{\theta})d_{,\theta} +  (g_{S,\theta}^{\theta \phi} + g_{S,\phi}^{\phi \phi} + g_{S}^{\theta \phi}\cot{\theta})d_{,\phi}\Big) \Black\\
& + \frac{1}{2}\left(\frac{|g_{S}|}{-|\gbar|}\right)^{1/2}\bigg\{\Big((g_{S,\theta}^{\theta \theta}+ g_{S,\phi}^{\phi \theta} + g_{S}^{\theta \theta}\cot{\theta})e_{,r} + (g_{S,\theta}^{\theta \phi} + g_{S,\phi}^{\phi \phi} + g_{S}^{\theta \phi}\cot{\theta})f_{,r}  \Big)\\
& -\frac{2d}{u} \Big((g_{S,\theta}^{\theta \theta} + g_{S,\phi}^{\phi \theta} + g_{S}^{\theta \theta}\cot{\theta})u_{,\theta} +(g_{S,\theta}^{\theta \phi} + g_{S,\phi}^{\phi \phi}+ g_{S}^{\theta \phi}\cot{\theta}) u_{,\phi}  \Big)\\
& + \frac{cf - be}{ab -  c^{2}}\Big((g_{S,\theta}^{\theta \theta} + g_{S,\phi}^{\phi \theta} + g_{S}^{\theta \theta}\cot{\theta})a_{,r}  + (g_{S,\theta}^{\theta \phi} + g_{S,\phi}^{\phi \phi} + g_{S}^{\theta \phi}\cot{\theta})c_{,r} \Big)\\
& + \frac{ce - af}{ab - c^{2}}\Big((g_{S,\theta}^{\theta \theta} + g_{S,\phi}^{\phi \theta} + g_{S}^{\theta \theta}\cot{\theta})c_{,r} + (g_{S,\theta}^{\theta \phi} + g_{S,\phi}^{\phi \phi} + g_{S}^{\theta \phi}\cot{\theta})b_{,r} \Big)\bigg\}.
\end{align*}

It can be computed that 
\[g_{S,\theta}^{\theta \theta} + g_{S,\phi}^{\phi \theta} + g_{S}^{\theta \theta}\cot{\theta} = \frac{1}{|g_{S}|}(b_{,\theta} - b\cot{\theta} - c_{,\phi});\]
and
\[g_{S,\theta}^{\theta \phi} + g_{S,\phi}^{\phi \phi} + g_{S}^{\theta \phi}\cot{\theta} = \frac{1}{|g_{S}|}(-c_{,\theta} + c\cot{\theta} + a_{,\phi}).\]

Using the above, we can simplify $\CornflowerBlue(l1) \Black$ as follows:
\begin{align}
\CornflowerBlue(l1) \Black & = \Peach\frac{1}{2}\left(\frac{-1}{|g_{S}||\gbar|}\right)^{1/2}\Big((b_{,\theta} - b\cot{\theta} - c_{,\phi})d_{,\theta} +  (-c_{,\theta} + c\cot{\theta} + a_{,\phi})d_{,\phi})\Big) \Black\notag\\
& + \frac{1}{2}\left(\frac{-1}{|g_{S}||\gbar|}\right)^{1/2}\Red\left[\Black(b_{,\theta} - b\cot{\theta} - c_{,\phi})\left(e_{,r} - 2\frac{u_{,\theta}}{u}d + \frac{cf - be}{ab - c^{2}}a_{,r} +\frac{ce - af}{ab - c^{2}}c_{,r}\right) \right.\notag\\
& \left.+ (-c_{,\theta} + c\cot{\theta} + a_{,\phi})\left(f_{,r} - 2\frac{u_{,\phi}}{u}d + \frac{cf - be}{ab - c^{2}}c_{,r} + \frac{ce - af}{ab - c^{2}}b_{,r} \right)  \Red \right].\Black\label{equ:l1 term simplified}
\end{align}

Plug (\ref{equ:scaled I term}), (\ref{equ:scaled II term}), (\ref{equ:scaled III term}), (\ref{equ:scaled IV term}), (\ref{equ:scaled E term}), (\ref{equ:l2 term simplified}) and (\ref{equ:l1 term simplified}) back into the divergence free equation (\ref{equ:fourth condition in straight out flow broken down}), we get: 

\begin{align}
& 2(-|g_{S}||\gbar|)^{1/2}\cdot\divg_{g_{S}}\alpha^{\nu} =  -2(-|g_{S}||\gbar|)^{1/2}\frac{1}{u^{2}}\frac{(-|\gbar|)^{1/2}}{|g_{S}|^{3/2}} \cdot \Big(\Red(I)\Black + \Green(II)\Black + \BlueViolet(III)\Black + \YellowOrange(IV)\Black \notag\\
& + \CornflowerBlue (E) \Black\Big) + \CornflowerBlue (l2) \Black + \CornflowerBlue (l1) \Black\tag{By (\ref{equ:fourth condition in straight out flow broken down})}\notag\\
& = \underline{(\Red bd_{,\theta \theta} - 2cd_{,\theta \phi} + ad_{,\phi \phi}\Black)} + (be_{,r \theta} - cf_{,r \theta} - ce_{, r \phi} + af_{,r\phi}) -\frac{d}{u^{2}}\Green(b(u^{2})_{,\theta \theta} - 2c(u^{2})_{,\theta \phi}  \Black\notag\\
& + a(u^{2})_{,\phi \phi}) + \frac{cf - be}{ab - c^{2}}\BlueViolet(ba_{,r \theta} - cc_{,r \theta} - ca_{,r \phi} + ac_{, r\phi})\Black +\frac{ce - af}{ab - c^{2}}\YellowOrange(bc_{,r\theta} - cb_{,r\theta} - cc_{,r\phi} + ab_{,r\phi})\Black\notag\\
& + \left(2\frac{u_{,\theta}}{u} + 2\cot{\theta} - \frac{|\gbar|_{,\theta}}{|\gbar|}\right)\RedViolet(be_{,r} + bd_{,\theta} - cf_{,r} - cd_{,\phi})\Black\notag\\
& + \left(2\frac{u_{,\phi}}{u} - \frac{|\gbar|_{,\phi}}{|\gbar|}\right) \RoyalPurple(- ce_{,r} - cd_{,\theta} + af_{,r} + ad_{,\phi})\Black - \frac{2d}{u}\Green\left[\Black\left(\frac{d_{,\theta}}{d} + 2\cot{\theta} - \frac{|\gbar|_{,\theta}}{|\gbar|} \right)(bu_{,\theta} - cu_{,\phi})\right. \notag\\
& \left.+ \left(\frac{d_{,\phi}}{d} - \frac{|\gbar|_{,\phi}}{|\gbar|} \right)(-cu_{,\theta} + au_{,\phi})\Green\right]\Black +\frac{cf - be}{ab - c^{2}}\Green\left[\Black\left(2\frac{u_{,\theta}}{u} + \frac{(cf -be)_{,\theta}}{cf - be} - \frac{|\gbar|_{,\theta}}{|\gbar|} \right)(ba_{,r} - cc_{,r})\right. \notag\\
& \left. + \left(2\frac{u_{,\phi}}{u} + \frac{(cf-be)_{,\phi}}{cf - be} -\frac{|\gbar|_{,\phi}}{|\gbar|} \right)(-ca_{,r} + ac_{,r})\Green\right]\Black + \frac{ce - af}{ab - c^{2}}\Green\left[\Black\left(2\frac{u_{,\theta}}{u} + \frac{(ce - af)_{,\theta}}{ce - af} - \frac{|\gbar|_{,\theta}}{|\gbar|} \right)\cdot\right. \notag\\
& \left. (bc_{,r} - cb_{,r}) + \left(2\frac{u_{,\phi}}{u} + \frac{(ce - af)_{,\phi}}{ce - af} - \frac{|\gbar|_{,\phi}}{|\gbar|}\right)(-cc_{,r} + ab_{,r})\Green\right]\Black\notag\\
& - \Green\left(\Black\frac{2u_{,\theta}}{u} + \cot{\theta} - \frac{|\gbar|_{,\theta}}{2|\gbar|}\Green\right)\Black\cdot \Green\Big(\Black\RedViolet(be_{,r} + bd_{,\theta} - cf_{,r} - cd_{,\phi})\Black - \frac{2d}{u}(bu_{,\theta} - cu_{,\phi})  \notag\\
& + \frac{cf - be}{ab - c^{2}}(ba_{,r} - cc_{,r}) + \frac{ce -af}{ab - c^{2}}(bc_{,r} - cb_{,r})  \Green\Big)\Black - \Green\left(\Black\frac{2u_{,\phi}}{u} - \frac{|\gbar|_{,\phi}}{2|\gbar|}  \Green\right)\Black \cdot \Green\Big(\Black\RoyalPurple(-ce_{,r} - cd_{,\theta}\Black \notag\\    
& \RoyalPurple + af_{,r} + ad_{,\phi})\Black - \frac{2d}{u}(-cu_{,\theta} + au_{,\phi}) + \frac{cf - be}{ab - c^{2}}(-ca_{,r} + ac_{,r}) + \frac{ce - af}{ab - c^{2}}(-cc_{,r} + ab_{,r}) \Green\Big)\Black\notag\\
& + \underline{\Peach(b_{,\theta} - b\cot{\theta} - c_{,\phi})d_{,\theta} +  (-c_{,\theta} + c\cot{\theta} + a_{,\phi})d_{,\phi})\Black}\notag\\
& + (b_{,\theta} - b\cot{\theta} - c_{,\phi})\left(e_{,r} - 2\frac{u_{,\theta}}{u}d + \frac{cf - be}{ab - c^{2}}a_{,r} +\frac{ce - af}{ab - c^{2}}c_{,r}\right) \notag\\
& + (-c_{,\theta} + c\cot{\theta} + a_{,\phi})\left(f_{,r} - 2\frac{u_{,\phi}}{u}d + \frac{cf - be}{ab - c^{2}}c_{,r} + \frac{ce - af}{ab - c^{2}}b_{,r} \right). \label{equ:scaled divergence equation step 1}
\end{align}

Notice that the terms underlined, up to a $|g_{S}|$ factor, equal the Laplacian of $d$ along $S_{t,r}$ (see equation (\ref{equ:laplacian of d on the surface str}) in Appendix \ref{appen:laplacian of coordinate sphere}). We can simplify (\ref{equ:scaled divergence equation step 1}) as:
\begin{align}
& (\ref{equ:scaled divergence equation step 1})  = \Red|g_{S}|\Delta_{g_{S}}d\Black + (be_{,r \theta} - cf_{,r \theta} - ce_{, r \phi} + af_{,r\phi})\notag\\
& -\frac{d}{u^{2}}\Green(b(u^{2})_{,\theta \theta} - 2c(u^{2})_{,\theta \phi} + a(u^{2})_{,\phi \phi}) \Black + \frac{cf - be}{ab - c^{2}}\BlueViolet(ba_{,r \theta} - cc_{,r \theta} - ca_{,r \phi} + ac_{, r\phi})\Black\notag\\
& +\frac{ce - af}{ab - c^{2}}\YellowOrange(bc_{,r\theta} - cb_{,r\theta} - cc_{,r\phi} + ab_{,r\phi})\Black + \left(\cot{\theta} - \frac{|\gbar|_{,\theta}}{2|\gbar|}\right)\RedViolet(be_{,r} + bd_{,\theta} - cf_{,r} - cd_{,\phi})\Black\notag\\
& - \frac{|\gbar|_{,\phi}}{2|\gbar|}\RoyalPurple(-ce_{,r} - cd_{,\theta} + af_{,r} + ad_{,\phi})\Black\notag\\
& - \frac{2d}{u}\Green\left[\Black\left(\frac{d_{,\theta}}{d} - \frac{2u_{,\theta}}{u} + \cot{\theta} - \frac{|\gbar|_{,\theta}}{2|\gbar|} \right)(bu_{,\theta} - cu_{,\phi}) + \left(\frac{d_{,\phi}}{d} - \frac{2u_{,\phi}}{u} - \frac{|\gbar|_{,\phi}}{2|\gbar|} \right)(-cu_{,\theta} + au_{,\phi})\Green\right]\Black  \notag\\
& +\frac{cf - be}{ab - c^{2}}\Green\bigg[\Black\left(\frac{(cf -be)_{,\theta}}{cf - be} - \frac{|\gbar|_{,\theta}}{2|\gbar|} - \cot{\theta}\right)(ba_{,r} - cc_{,r}) \Black\notag\\
& + \left(\frac{(cf-be)_{,\phi}}{cf - be} -\frac{|\gbar|_{,\phi}}{2|\gbar|} \right)(-ca_{,r} + ac_{,r})\Green\bigg] \Black \notag \\
& + \frac{ce - af}{ab - c^{2}}\Green\bigg[\Black\left(\frac{(ce - af)_{,\theta}}{ce - af} - \frac{|\gbar|_{,\theta}}{2|\gbar|} - \cot{\theta} \right)(bc_{,r} - cb_{,r}) \Black \notag\\
& + \left(\frac{(ce - af)_{,\phi}}{ce - af} - \frac{|\gbar|_{,\phi}}{2|\gbar|}\right)(-cc_{,r} + ab_{,r})\Green\bigg]\Black \notag\\
& + (b_{,\theta} - b\cot{\theta} - c_{,\phi})\left(e_{,r} - \frac{2d}{u}u_{,\theta} + \frac{cf - be}{ab - c^{2}}a_{,r} +\frac{ce - af}{ab - c^{2}}c_{,r}\right) \notag\\
& + (-c_{,\theta} + c\cot{\theta} + a_{,\phi})\left(f_{,r} - \frac{2d}{u}u_{,\phi} + \frac{cf - be}{ab - c^{2}}c_{,r} + \frac{ce - af}{ab - c^{2}}b_{,r} \right).\label{equ:scaled straight out equation before simplification 1}
\end{align}

\begin{mydef}\label{def:second order derivative terms in straight out flow equation}
We set $L_{2}(\gbar) = L_{2}(d,e,f,u,a,b,c)$ to be the second derivative terms in the straight out flow equation, that is:
\begin{align}
& L_{2}(\gbar) := \Red|g_{S}|\Delta_{g_{S}}d\Black + (be_{,r \theta} - cf_{,r \theta} - ce_{, r \phi} + af_{,r\phi}) -\frac{d}{u^{2}}\Green(b(u^{2})_{,\theta \theta} - 2c(u^{2})_{,\theta \phi} + a(u^{2})_{,\phi \phi}) \Black\notag\\
& + \frac{cf - be}{ab - c^{2}}\BlueViolet(ba_{,r \theta} - cc_{,r \theta} - ca_{,r \phi} + ac_{, r\phi})\Black +\frac{ce - af}{ab - c^{2}}\YellowOrange(bc_{,r\theta} - cb_{,r\theta} - cc_{,r\phi} + ab_{,r\phi}).\Black\label{equ:second order derivative terms in straight out flow equation}
\end{align}
\end{mydef}

With this defintion, we further simplify (\ref{equ:scaled straight out equation before simplification 1}) as follows:

\begin{align}
& (\ref{equ:scaled straight out equation before simplification 1}) = L_{2}(\gbar) + \cot{\theta}\RedViolet(\cancel{be_{,r}} + bd_{,\theta} - \bcancel{cf_{,r}} - cd_{,\phi})\Black -  \frac{|\gbar|_{,\theta}}{2|\gbar|}\RedViolet(be_{,r} + bd_{,\theta} - cf_{,r} - cd_{,\phi})\Black\notag\\
& - \frac{|\gbar|_{,\phi}}{2|\gbar|}\RoyalPurple(-ce_{,r} - cd_{,\theta} + af_{,r} + ad_{,\phi})\Black - \frac{2d}{u}\Green\Big[\Black\left(\frac{d_{,\theta}}{d} - \frac{2u_{,\theta}}{u} + \cot{\theta} - \frac{|\gbar|_{,\theta}}{2|\gbar|} \right)(bu_{,\theta} - cu_{,\phi})\notag\\
& + \left(\frac{d_{,\phi}}{d} - \frac{2u_{,\phi}}{u} - \frac{|\gbar|_{,\phi}}{2|\gbar|} \right)(-cu_{,\theta} + au_{,\phi})\Green\Big]\Black  \notag\\
& +\frac{cf - be}{ab - c^{2}}\Green\Big[\Black\left(\frac{(cf -be)_{,\theta}}{cf - be} - \frac{|\gbar|_{,\theta}}{2|\gbar|} - \cot{\theta}\right)(ba_{,r} - cc_{,r}) + \left(\frac{(cf-be)_{,\phi}}{cf - be} -\frac{|\gbar|_{,\phi}}{2|\gbar|} \right)\cdot\notag\\
& (-ca_{,r} + ac_{,r})\Green\Big]\Black\notag\\
& + \frac{ce - af}{ab - c^{2}}\Green\Big[\Black\left(\frac{(ce - af)_{,\theta}}{ce - af} - \frac{|\gbar|_{,\theta}}{2|\gbar|} - \cot{\theta} \right)(bc_{,r} - cb_{,r}) + \left(\frac{(ce - af)_{,\phi}}{ce - af} - \frac{|\gbar|_{,\phi}}{2|\gbar|}\right)\cdot \notag\\
& (-cc_{,r} + ab_{,r})\Green\Big]\Black \notag\\
& + b_{,\theta}e_{,r} - \frac{2d}{u}b_{,\theta}u_{,\theta} + \frac{cf - be}{ab - c^{2}}b_{,\theta}a_{,r} + \frac{ce - af}{ab - c^{2}}b_{,\theta}c_{,r} - \cancel{be_{,r}\cot{\theta}} \underline{+ \frac{2d}{u}bu_{,\theta}\cot{\theta}} \notag\\
& - \frac{cf - be}{ab - c^{2}}ba_{,r}\cot{\theta} - \frac{ce - af}{ab - c^{2}}bc_{,r}\cot{\theta} - c_{,\phi}e_{,r} + \frac{2d}{u}c_{,\phi}u_{,\theta} - \frac{cf - be}{ab - c^{2}}c_{,\phi}a_{,r} \notag\\
& - \frac{ce - af}{ab - c^{2}}c_{,\phi}c_{,r} - c_{,\theta}f_{,r} + \frac{2d}{u}c_{,\theta}u_{,\phi} - \frac{cf - be}{ab - c^{2}}c_{,\theta}c_{,r} - \frac{ce - af}{ab - c^{2}}c_{,\theta}b_{,r} + \bcancel{cf_{,r}\cot{\theta}}\notag\\
&  \underline{- \frac{2d}{u}cu_{,\phi}\cot{\theta}} + \frac{cf - be}{ab - c^{2}}cc_{,r}\cot{\theta} + \frac{ce - af}{ab - c^{2}}cb_{,r}\cot{\theta} + a_{,\phi}f_{,r} - \frac{2d}{u}a_{,\phi}u_{,\phi} \notag\\
&+ \frac{cf - be}{ab - c^{2}}a_{,\phi}c_{,r} + \frac{ce - af}{ab - c^{2}}a_{,\phi}b_{,r}.   \label{equ:scaled straight out equation simplification step 1}
\end{align}

\begin{align}
& (\ref{equ:scaled straight out equation simplification step 1}) = L_{2}(\gbar) + \cot{\theta}\RedViolet(bd_{,\theta} - cd_{,\phi})\Black -  \frac{|\gbar|_{,\theta}}{2|\gbar|}\RedViolet(be_{,r} + bd_{,\theta} - cf_{,r} - cd_{,\phi})\Black \notag\\
& - \frac{|\gbar|_{,\phi}}{2|\gbar|}\RoyalPurple(-ce_{,r} - cd_{,\theta} + af_{,r} + ad_{,\phi})\Black\notag\\
& - \frac{2d}{u}\Green\Big[\Black\left(\frac{d_{,\theta}}{d} - \frac{2u_{,\theta}}{u} + \cot{\theta} - \frac{|\gbar|_{,\theta}}{2|\gbar|} \right)(bu_{,\theta} - cu_{,\phi}) + \left(\frac{d_{,\phi}}{d} - \frac{2u_{,\phi}}{u} - \frac{|\gbar|_{,\phi}}{2|\gbar|} \right)\cdot \notag\\
& (-cu_{,\theta} + au_{,\phi})\Green\Big]\Black  \notag\\
& +\frac{cf - be}{ab - c^{2}}\Green\Big[\Black\left(\frac{(cf -be)_{,\theta}}{cf - be} - \frac{|\gbar|_{,\theta}}{2|\gbar|} - \cot{\theta}\right)(ba_{,r} - cc_{,r}) + \left(\frac{(cf-be)_{,\phi}}{cf - be} -\frac{|\gbar|_{,\phi}}{2|\gbar|} \right)\cdot\notag\\
& (-ca_{,r} + ac_{,r})\Green\Big]\Black\notag\\
& + \frac{ce - af}{ab - c^{2}}\Green\Big[\Black\left(\frac{(ce - af)_{,\theta}}{ce - af} - \frac{|\gbar|_{,\theta}}{2|\gbar|} - \cot{\theta} \right)(bc_{,r} - cb_{,r}) + \left(\frac{(ce - af)_{,\phi}}{ce - af} - \frac{|\gbar|_{,\phi}}{2|\gbar|}\right)\cdot \notag\\
& (-cc_{,r} + ab_{,r})\Green\Big]\Black \notag\\
& + (b_{,\theta}e_{,r} - c_{,\phi}e_{,r} - c_{,\theta}f_{,r} + a_{,\phi}f_{,r}) - \frac{2d}{u}(b_{,\theta}u_{,\theta} - c_{,\phi}u_{,\theta} - c_{,\theta}c_{\phi} + a_{,\phi}u_{,\phi})\notag\\
& + \frac{cf - be}{ab - c^{2}}(b_{,\theta}a_{,r} - a_{,r}c_{,\phi} - c_{,r}c_{,\theta} + a_{,\phi}c_{,r}) + \frac{ce - af}{ab - c^{2}}(b_{,\theta}c_{,r} - c_{,r}c_{,\phi} - b_{,r}c_{,\theta} + a_{,\phi}b_{,r})\notag\\
& + \underline{\frac{2d}{u}\cot{\theta}(bu_{,\theta} - cu_{,\phi})} - \frac{cf - be}{ab - c^{2}}\cot{\theta} (ba_{,r}  - cc_{,r}) - \frac{ce - af}{ab - c^{2}}\cot{\theta} (bc_{,r} - cb_{,r})\label{equ:scaled straight out equation before simplification 2}
\end{align}

Finally the above can be simplified to:
\begin{align}
& L_{2}(\gbar) + \cot{\theta}\RedViolet(bd_{,\theta} - cd_{,\phi})\Black -  \frac{|\gbar|_{,\theta}}{2|\gbar|}\RedViolet(be_{,r} + bd_{,\theta} - cf_{,r} - cd_{,\phi})\Black - \frac{|\gbar|_{,\phi}}{2|\gbar|}\RoyalPurple(-ce_{,r} - cd_{,\theta} + af_{,r} \Black\notag\\
& + ad_{,\phi}) \notag\\
& - \frac{2d}{u}\Green\left[\Black\left(\frac{d_{,\theta}}{d} - \frac{2u_{,\theta}}{u} - \frac{|\gbar|_{,\theta}}{2|\gbar|} \right)(bu_{,\theta} - cu_{,\phi}) + \left(\frac{d_{,\phi}}{d} - \frac{2u_{,\phi}}{u} - \frac{|\gbar|_{,\phi}}{2|\gbar|} \right)(-cu_{,\theta} + au_{,\phi})\Green\right]\Black  \notag\\
& +\frac{cf - be}{ab - c^{2}}\Green\bigg[\Black\left(\frac{(cf -be)_{,\theta}}{cf - be} - \frac{|\gbar|_{,\theta}}{2|\gbar|} - 2\cot{\theta}\right)(ba_{,r} - cc_{,r}) + \left(\frac{(cf-be)_{,\phi}}{cf - be} -\frac{|\gbar|_{,\phi}}{2|\gbar|} \right)\cdot \notag\\
& (-ca_{,r} + ac_{,r})\Green\bigg]\Black\notag\\
& + \frac{ce - af}{ab - c^{2}}\Green\bigg[\Black\left(\frac{(ce - af)_{,\theta}}{ce - af} - \frac{|\gbar|_{,\theta}}{2|\gbar|} - 2\cot{\theta} \right)(bc_{,r} - cb_{,r}) + \left(\frac{(ce - af)_{,\phi}}{ce - af} - \frac{|\gbar|_{,\phi}}{2|\gbar|}\right)\notag\\
& (-cc_{,r} + ab_{,r})\Green\bigg]\Black \notag\\
& + (b_{,\theta}e_{,r} - c_{,\phi}e_{,r} - c_{,\theta}f_{,r} + a_{,\phi}f_{,r}) - \frac{2d}{u}(b_{,\theta}u_{,\theta} - c_{,\phi}u_{,\theta} - c_{,\theta}c_{\phi} + a_{,\phi}u_{,\phi})\notag\\
& + \frac{cf - be}{ab - c^{2}}(b_{,\theta}a_{,r} - a_{,r}c_{,\phi} - c_{,r}c_{,\theta} + a_{,\phi}c_{,r}) + \frac{ce - af}{ab - c^{2}}(b_{,\theta}c_{,r} - c_{,r}c_{,\phi} - b_{,r}c_{,\theta} + a_{,\phi}b_{,r}).\label{equ:scaled straight out equation simplified}\\
& =:\Red|g_{S}|\Delta_{g_{S}}d\Black + F(d, d').\notag
\end{align}

Therefore we have the following characterization of the fourth condition (\ref{def:fourth conditions of straight out coordinate chart}):
\begin{proposition}
The fourth condition (\ref{def:fourth conditions of straight out coordinate chart}) is a second order elliptic PDE in $d$:
\begin{equation}\label{equ:straight out flow equals second order elliptic in d}
\Red \Delta_{g_{S}}d \Black + \CornflowerBlue G(d, d') \Black = 0,
\end{equation}
where $G := \frac{F}{|g_{S}|} = \frac{F}{r^{4}\sin^{2}{\theta}}$. 
\end{proposition}

If Equation $(\ref{equ:straight out flow equals second order elliptic in d})$ is solvable, then many spacetimes that admit straight out flow coordinate chart exist, even beyond spherically symmetric examples. A necessary condition for (\ref{equ:straight out flow equals second order elliptic in d}) to be solvable if
\begin{equation}\label{equ:obstruction for poisson equation}
\int_{S_{t, r}}G(d, d')\,\,dA_{t, r} = 0.
\end{equation}
We conjecture that this is always the case. The verification of (\ref{equ:obstruction for poisson equation}) is still work in progress.
\chapter{Conclusions and Open Problems}
\label{chap:futurework}

We have constructed many examples of non-spherically symmetric, non-static spacetimes that admit smooth global solutions to inverse mean curvature vector flow. Prior to our work, such solutions were only known in spherically symmetric and static spacetimes. Our work seems to suggest that spacetimes that admit inverse mean curvature vector flow solutions might exist generically. However, this more general problem is still open:

\begin{problem}
Given an arbitrary spacetime. Can we always find a ``right'' initial surface such that inverse mean curvature vector flow starting with this surface exists for all time?
\end{problem}

Going to the big picture of relating local and global notions of mass, it is still unknown that:

\begin{problem}
Given a spacetime that is sufficiently asymptotically flat (e.g. 

Schwarzschild outside a compact set). Does the Hawking of inverse mean curvature vector flow surfaces approach the total mass of the spacetime?
\end{problem}

A natural next step is to consider the following problem:
\begin{problem}
Given a spherically symmetric spacetime that admits an inverse mean curvature vector flow coordinate chart. Consider a perturbation of the spacetime metric. Does the perturbed metric admit an inverse mean curvature vector flow coordinate chart as well?
\end{problem}

We conjecture that this is always the case for some Minkowski spacetimes:
\begin{conjecture}
\label{conj:perturbation of minkowski spacetime admits IMCVF coordinate chart}
Given Minkowski space with inverse mean curvature vector flow coordinate chart that can be smoothly extended to the boundary, consider a perturbation of the spacetime metric. The resulting spacetime still admits inverse mean curvature vector flow solutions (in a single spacelike hypersurface) that exist for all time. 
\end{conjecture}

\appendix
\chapter{Geometric Calculations}
\label{chap:geometric calculations}

\section{Ricci, Scalar and Einstein Curvature of Spherically Symmetric Spacetime}\label{appen:curvatures in spherically symmetric spacetime}
In this subsection, we are going to compute the Ricci, Scalar and Einstein curvature of the spherically symmetric space time $(N^{4}, \gbar)$, with coordinates $(t, r, \theta, \phi)$, such that $g$ has the local coordinate representation as in $(\ref{equ:coordinate repren of spherically symmetric metric})$:
\begin{equation}\label{spherically symmetric metrics satisfying radial mean curvature vector condition again again}
\gbar = \bordermatrix{~ & t & r & \theta & \phi \cr
                  t & -v^{2}(t, r) & 0 & 0 & 0\cr
                  r & 0 & u^{2}(t, r) & 0 & 0\cr
                  \theta & 0 & 0 &  r^{2} & 0\cr
                  \phi & 0 & 0 & 0 & r^{2}\sin^{2}{\theta}\cr}
\end{equation}
Note that $\{\p_{t}, \p_{r}, \p_{\theta}, \p_{\phi}\}$ form a local frame of the tangent bundle. We assume that the connection on $(N^{4}, \gbar)$ is the Levi-Civita connection. The Einstein summation convention will be used, and colons will always denote the coordinate chart derivative whereas semicolons will denote covariant derivatives.In the following, we will use Latin letters $i, j, k, l$ and so on to be indices taking values in $\{t, r, \theta, \phi\}$.  Recall the following formulas from Riemannian geometry:
\begin{align}
& \Gamma^{k}_{ij} = \frac{1}{2}g^{kl}(g_{jl, i} + g_{il, j} - g_{ij, l}) \tag{Christoffel symbols}\label{c symbols}\\
& \Ric_{ij} = \Ric(\p_{i}, \p_{j}) = \Gamma^{k}_{ij, k} - \Gamma^{k}_{ik, j} + \Gamma^{k}_{km}\Gamma^{m}_{ij} - \Gamma^{k}_{jm}\Gamma^{m}_{ik} \tag{Ricci curvature} \label{Ricci curvature}\\
& \R = \text{tr}_{g}\Ric \tag{scalar curvature} \label{scalar curvature}
\end{align}
We will write the Christoffel symbols in four matrices $\Gamma^{t}, \Gamma^{r},\Gamma^{\theta}$, and $\Gamma^{\phi}$ as in Alan Parry's survey paper [\cite{Parry2012}]. Now we compute them in sequence, using the fact that the metric $g$ is diagonal in our coordinates. 
\begin{enumerate}
\item The computation of $\Gamma^{t}$:
\begin{itemize}
\item $\Gamma^{t}_{tt} = \frac{1}{2}g^{tt}(g_{tt, t}) = \frac{1}{2}\frac{1}{v^{2}} 2 vv_{, t} = \frac{v_{, t}}{v}$.
\item $\Gamma^{t}_{tr} = \Gamma^{t}_{rt} = \frac{1}{2}g^{tt}(g_{rt, t} + g_{tt, r}  - g_{tr, t}) = \frac{1}{2}\frac{1}{v^{2}}(2vv_{, r}) = \frac{v_{, r}}{v}$.
\item $\Gamma^{t}_{t\theta} = \Gamma^{t}_{\theta t} = \frac{1}{2}g^{tt}(g_{\theta t, t} + g_{tt, \theta} - g_{t\theta, t}) = 0$.
\item $\Gamma^{t}_{t\phi} = \Gamma^{t}_{\phi t} = \frac{1}{2}g^{tt}(g_{\phi t, t} + g_{tt, \phi} - g_{t\phi, t}) = 0$.
\item $\Gamma^{t}_{rr} = \frac{1}{2}g^{tt}(2g_{rt, r} - g_{rr, t}) = \frac{1}{2}\frac{1}{v^{2}}(2u u_{, t}) = \frac{u u_{, t}}{v^{2}}$.
\item $\Gamma^{t}_{r\theta} = \Gamma^{t}_{\theta r} = \frac{1}{2}g^{tt}(g_{\theta t, r} + g_{rt, \theta} - g_{r\theta, t}) = 0$.
\item $\Gamma^{t}_{r\phi} = \Gamma^{t}_{\phi r} = \frac{1}{2}g^{tt}(g_{\phi t, r} + g_{rt, \phi} - g_{r\phi, t}) = 0$.
\item $\Gamma^{t}_{\theta \theta} = \frac{1}{2}g^{tt}(2g_{\theta t, \theta} - g_{\theta \theta, t}) = 0$.
\item $\Gamma^{t}_{\theta \phi} = \Gamma^{t}_{\phi \theta} = \frac{1}{2}g^{tt}(g_{\phi t, \theta} + g_{\theta t, \phi} - g_{\theta \phi, t}) = 0$.
\item $\Gamma^{t}_{\phi \phi} = \frac{1}{2}g^{tt}(2g_{\phi t, \phi} - g_{\phi \phi, t}) = 0$.
\end{itemize}
Therefore we have
\begin{equation}\label{$Gamma^{t}$ matrix}
\Gamma^{t} = \bordermatrix{~ & t & r & \theta & \phi \cr
                  t & \frac{v_{, t}}{v} & \frac{v_{, r}}{v} & 0 & 0\cr
                  r & \frac{v_{, r}}{v} & \frac{u u_{, t}}{v^{2}} & 0 & 0\cr
                  \theta & 0 & 0 &  0 & 0\cr
                  \phi & 0 & 0 & 0 & 0 \cr}
\end{equation}

\item The computation of $\Gamma^{r}$:
\begin{itemize}
\item $\Gamma^{r}_{tt} = \frac{1}{2}g^{rr}(2g_{tr, t} - g_{tt, r}) = \frac{1}{2}\frac{1}{u^{2}}(2vv_{, r}) = \frac{vv_{, r}}{u^{2}}$.
\item $\Gamma^{r}_{tr} = \Gamma^{r}_{rt} = \frac{1}{2}g^{rr}(g_{rr, t} + g_{tr, r} - g_{tr, r}) = \frac{1}{2}\frac{1}{u^{2}}(2u u_{,t}) = \frac{u_{, t}}{u}$.
\item $\Gamma^{r}_{t\theta} = \Gamma^{r}_{\theta t} = \frac{1}{2}g^{rr}(g_{\theta r, t} + g_{tr, \theta} - g_{t \theta, r}) = 0$.
\item $\Gamma^{r}_{t\phi} = \Gamma^{r}_{\phi t} = \frac{1}{2}g^{rr}(g_{\phi r, t} + g_{tr, \phi} - g_{t\phi, r}) = 0$.
\item $\Gamma^{r}_{rr} = \frac{1}{2}g^{rr}g_{rr, r} = \frac{1}{2}\frac{1}{u^{2}}(2u u_{, r}) = \frac{u_{, r}}{u}$.
\item $\Gamma^{r}_{r \theta} = \Gamma^{r}_{\theta r} = \frac{1}{2}g^{rr}(g_{\theta r, r} + g_{rr, \theta} - g_{r\theta, r}) = 0$.
\item $\Gamma^{r}_{r\phi} = \Gamma^{r}_{\phi r} = \frac{1}{2}g^{rr}(g_{\phi r, r} + g_{rr, \phi} - g_{r\phi, r}) = 0$.
\item $\Gamma^{r}_{\theta \theta} = \frac{1}{2}(2g_{\theta r, \theta} - g_{\theta \theta, r}) = \frac{1}{2}\frac{1}{u^{2}} (-2r) = - \frac{r}{u^{2}}$.
\item $\Gamma^{r}_{\theta \phi} = \Gamma^{r}_{\phi \theta} = \frac{1}{2}g^{rr}(g_{\phi r, \theta} + g_{\theta r, \phi} - g_{\theta \phi, r}) = 0$.
\item $\Gamma^{r}_{\phi \phi} = \frac{1}{2}g^{rr}(2g_{\phi r, \phi} - g_{\phi \phi, r}) = \frac{1}{2}\frac{1}{u^{2}}(-2r\sin^{2}{\theta}) = -\frac{r\sin^{2}{\theta}}{u^{2}}$.
\end{itemize}
Therefore we have
\begin{equation}\label{$Gamma^{r}$ matrix}
\Gamma^{r} = \bordermatrix{~ & t & r & \theta & \phi \cr
                  t & \frac{vv_{, r}}{u^{2}} & \frac{u_{, t}}{u} & 0 & 0\cr
                  r & \frac{u_{, t}}{u} & \frac{u_{, r}}{u} & 0 & 0\cr
                  \theta & 0 & 0 &  -\frac{r}{u^{2}} & 0\cr
                  \phi & 0 & 0 & 0 & -\frac{r\sin^{2}{\theta}}{u^{2}} \cr}
\end{equation}

\item The computation of $\Gamma^{\theta}$:
\begin{itemize}
\item $\Gamma^{\theta}_{tt} = \frac{1}{2}g^{\theta \theta}(2g_{t\theta, t} - g_{tt, \theta}) = 0$.
\item $\Gamma^{\theta}_{tr} = \Gamma^{\theta}_{rt} = \frac{1}{2}g^{\theta \theta}(g_{t\theta, r} + g_{r\theta, t} - g_{tt, \theta}) = 0$.
\item $\Gamma^{\theta}_{t\theta} = \Gamma^{\theta}_{\theta t} = \frac{1}{2}g^{\theta \theta}(g_{\theta \theta, t} + g_{t\theta, \theta} - g_{t \theta, \theta}) = 0$.
\item $\Gamma^{\theta}_{t\phi} = \Gamma^{\theta}_{\phi t} = \frac{1}{2}g^{\theta \theta}(g_{\phi \theta, t} + g_{t\theta,\phi} - g_{t\phi, \theta}) = 0$.
\item $\Gamma^{\theta}_{rr} = \frac{1}{2}g^{\theta \theta}(2g_{r\theta, r} - g_{rr,\theta}) = 0$.
\item $\Gamma^{\theta}_{r \theta} = \Gamma^{\theta}_{\theta r} = \frac{1}{2}g^{\theta \theta}(g_{\theta \theta, r} + g_{r\theta, \theta} - g_{r\theta, \theta}) = \frac{1}{2}\frac{1}{r^{2}}2r = \frac{1}{r}$.
\item $\Gamma^{\theta}_{r\phi} = \Gamma^{\theta}_{\phi r} = \frac{1}{2}g^{\theta \theta}(g_{\phi \theta, r} + g_{r\theta, \phi} - g_{r\phi, \theta}) = 0$.
\item $\Gamma^{\theta}_{\theta \theta} = \frac{1}{2}g^{\theta \theta}g_{\theta \theta, \theta} = 0$.
\item $\Gamma^{\theta}_{\theta \phi} = \Gamma^{\theta}_{\phi \theta} = \frac{1}{2}g^{\theta \theta}(g_{\phi \theta, \theta} + g_{\theta \theta, \phi} - g_{\theta \phi, \theta}) = 0$.
\item $\Gamma^{\theta}_{\phi \phi} = \frac{1}{2}g^{\theta \theta}(2g_{\phi \theta, \phi} - g_{\phi \phi, \theta}) = \frac{1}{2}\frac{1}{r^{2}}(-2r^{2}\sin{\theta}\cos{\theta}) = -sin{\theta}\cos{\theta}$.
\end{itemize}

Therefore we have
\begin{equation}\label{$Gamma^{theta}$ matrix}
\Gamma^{\theta} = \bordermatrix{~ & t & r & \theta & \phi \cr
                  t & 0 & 0 & 0 & 0\cr
                  r & 0 & 0 & \frac{1}{r} & 0\cr
                  \theta & 0 & \frac{1}{r} & 0 & 0\cr
                  \phi & 0 & 0 & 0 & -\sin{\theta}\cos{\theta} \cr}
\end{equation}

\item The computation of $\Gamma^{\phi}$:
\begin{itemize}
\item $\Gamma^{\phi}_{tt} = \frac{1}{2}g^{\phi \phi}(2g_{t\phi, t} - g_{tt, \phi}) = 0$.
\item $\Gamma^{\phi}_{tr} = \Gamma^{\phi}_{rt} = \frac{1}{2}g^{\phi \phi}(g_{r\phi, t} + g_{t\phi, r} - g_{tr, \phi}) = 0$.
\item $\Gamma^{\phi}_{t\theta} = \Gamma^{\phi}_{\theta t} = \frac{1}{2}g^{\phi \phi}(g_{\theta \phi, t} + g_{t\phi, \theta} - g_{t\theta, \phi}) = 0$.
\item $\Gamma^{\phi}_{t\phi} = \Gamma^{\phi}_{\phi t} = \frac{1}{2}g^{\phi\phi}(g_{\phi \phi, t} + g_{t\phi, \phi} - g_{t\phi, \phi}) = 0$.
\item $\Gamma^{\phi}_{rr} = \frac{1}{2}g^{\phi \phi}(2g_{r\phi, r} - g_{rr, \phi}) = 0$.
\item $\Gamma^{\phi}_{r\theta} = \Gamma^{\phi}_{\theta r} = \frac{1}{2}g^{\phi \phi} (g_{\theta \phi, r} + g_{r\phi, \theta} - g_{r\theta, \phi}) = 0$.
\item $\Gamma^{\phi}_{r\phi} = \Gamma^{\phi}_{\phi r} = \frac{1}{2}g^{\phi \phi}(g_{\phi \phi, r} + g_{r\phi, \phi} - g_{r\phi, \phi}) = \frac{1}{2}\frac{1}{r^{2}\sin^{2}{\theta}}2r\sin^{2}{\theta} = \frac{1}{r}$.
\item $\Gamma^{\phi}_{\theta \theta} = \frac{1}{2}g^{\phi \phi} (2g_{\theta \phi, \theta} - g_{\theta \theta, \phi}) = 0$.
\item $\Gamma^{\phi}_{\theta \phi} = \Gamma^{\phi}_{\phi \theta} = \frac{1}{2}g^{\phi \phi}(g_{\phi \phi, \theta} + g_{\theta \phi, \phi} - g_{\theta \phi, \phi}) = \frac{1}{2}\frac{1}{r^{2}\sin^{2}{\theta}}r^{2}2\sin{\theta}\cos{\theta} = \frac{\cos{\theta}}{\sin{\theta}}$.
\item $\Gamma^{\phi}_{\phi \phi} = \frac{1}{2}g^{\phi \phi}g_{\phi \phi, \phi} = 0$.
\end{itemize}

Therefore we have
\begin{equation}\label{$Gamma^{phi}$ matrix}
\Gamma^{\phi} = \bordermatrix{~ & t & r & \theta & \phi \cr
                  t & 0 & 0 & 0 & 0\cr
                  r & 0 & 0 & 0 & \frac{1}{r}\cr
                  \theta & 0 & 0 & 0 & \frac{\cos{\theta}}{\sin{\theta}}\cr
                  \phi & 0 & \frac{1}{r} & \frac{\cos{\theta}}{\sin{\theta}} & 0 \cr}
\end{equation}
\end{enumerate}

Now we compute each of the components of the Ricci curvature with respect to the basis $\{\p_{t}, \p_{r}, \p_{\theta}, \p_{\phi}\}$.
\begin{enumerate}
\item 
\begin{align*}
\Ric_{tt} & = \Ric(\p_{t}, \p_{t}) \\
& = \Gamma_{tt, r}^{r}  - \Gamma_{tr, t}^{r} + \Gamma_{rt}^{r}\Gamma_{tt}^{t} + \Gamma_{rr}^{r}\Gamma_{tt}^{r} - \Gamma_{tt}^{r}\Gamma_{tr}^{t} - \Gamma_{tr}^{r}\Gamma_{tr}^{r}  + \Gamma_{\theta r}^{\theta}\Gamma_{tt}^{r} + \Gamma_{\phi r}^{\phi} \Gamma_{tt}^{r}\\
& = \left(\frac{v v_{, r}}{u^{2}} \right)_{, r} - \left(\frac{u_{,t}}{u} \right)_{, t} + \frac{u_{, t}}{u}\frac{v_{, t}}{v} + \frac{u_{, r}}{u}\frac{vv_{, r}}{u^{2}} - \frac{vv_{, r}}{u^{2}}\frac{v_{, r}}{v} - \frac{u_{, t}}{u}\frac{u_{, t}}{u} + \frac{2}{r}\frac{v v_{, r}}{u^{2}}\\
& = \frac{1}{u^{2}}\left(vv_{, rr} + \frac{2}{r}vv_{, r} \right) - \frac{u_{, r} }{u^{3}} v v_{, r} +\frac{1}{u}\left(\frac{u_{, t}v_{, t}}{v} - u_{, tt} \right)
\end{align*} 

\item
\begin{align*}
\Ric_{tr} & = \Ric(\p_{t}, \p_{r}) = \Gamma_{tr, t}^{t} - \Gamma_{tt, r}^{t} + \Gamma_{tr}^{t}\Gamma_{tr}^{r} - \Gamma_{rr}^{t}\Gamma_{tt}^{r} + \Gamma_{\theta r}^{\theta}\Gamma_{tr}^{r} + \Gamma_{\phi r}^{\phi}\Gamma_{tr}^{r} \\
& = \left(\frac{v_{, r}}{v} \right)_{, t} - \left(\frac{v_{, t}}{v} \right)_{, r} + \frac{v_{, r}}{v}\frac{u_{, t}}{u} - \frac{uu_{, t}}{v^{2}} \frac{vv_{, r}}{u^{2}} + \frac{1}{r}\frac{u_{, t}}{u} + \frac{1}{r}\frac{u_{, t}}{u}\\
& =  \frac{2}{r} \frac{u_{, t}}{u}
\end{align*}

\item
\begin{align*}
\Ric_{t\theta} & = \Ric_{\theta t} = - \Gamma_{tt, \theta}^{t} - \Gamma_{tr, \theta}^{r} = 0
\end{align*}

\item
\begin{align*}
\Ric_{t\phi} & = \Ric_{\phi t} = - \Gamma_{tt, \phi}^{t} - \Gamma_{tr,\phi}^{r} = 0
\end{align*}

\item
\begin{align*}
\Ric_{rr} & = \Gamma_{rr, t}^{t} - \Gamma_{rt, r}^{t} + \Gamma_{tt}^{t} \Gamma_{rr}^{t} + \Gamma_{tr}^{t} \Gamma_{rr}^{r} - \Gamma_{rt}^{t}\Gamma_{rt}^{t} - \Gamma_{rr}^{t} \Gamma_{rt}^{r} - \Gamma_{r \theta, r}^{\theta} + \Gamma_{\theta r}^{\theta}\Gamma_{rr}^{r} - \Gamma_{r\theta}^{\theta}\Gamma_{r\theta}^{\theta}  \\
& - \Gamma_{r\phi, r}^{\phi} + \Gamma_{\phi r}^{\phi}\Gamma_{rr}^{r} - \Gamma_{r\phi}^{\phi}\Gamma_{r\phi}^{\phi}\\
& = \left(\frac{uu_{,t}}{v^{2}} \right)_{, t} - \left(\frac{v_{, r}}{v} \right)_{, r} + \frac{v_{, t}}{v} \frac{uu_{, t}}{v^{2}} + \frac{v_{, r}}{v}\frac{u_{, r}}{u} - \frac{v_{, r}}{v} \frac{v_{, r}}{v}  - \frac{uu_{, t}}{v^{2}}\frac{u_{, t}}{u} + \frac{2}{r}\frac{u_{, r}}{u}\\
& = -\frac{v_{, rr}}{v} + \frac{2}{r}\frac{u_{, r}}{u} - \frac{v_{, t}}{v^{3}}uu_{, t} + \frac{1}{v}\left(\frac{u u_{, tt}}{v}  + \frac{v_{, r}u_{, r}}{u}\right)
\end{align*}

\item
\begin{align*}
\Ric_{r\theta} & = \Ric_{\theta r} = 0
\end{align*}

\item
\begin{align*}
\Ric_{r\phi} & = \Ric_{\phi r} = 0
\end{align*}

\item
\begin{align*}
\Ric_{\theta \theta} & = \Gamma_{tr}^{t}\Gamma_{\theta \theta}^{r} + \Gamma_{\theta \theta, r}^{r} + \Gamma_{rr}^{r}\Gamma_{\theta \theta}^{r} - \Gamma_{\theta \theta}^{r}\Gamma_{\theta r}^{\theta} - \Gamma_{\theta \phi, \theta}^{\phi} + \Gamma_{\phi r}^{\phi}\Gamma_{\theta \theta}^{r} - \Gamma_{\theta \phi}^{\phi} \Gamma_{\theta \phi}^{\phi} \\
& = \frac{v_{, r}}{v} \frac{-r}{u^{2}} + \left( \frac{-r}{u^{2}} \right)_{, r} + \frac{u_{, r}}{u}\frac{-r}{u^{2}} - \frac{-r}{u^{2}}\frac{1}{r} - \left(\frac{\cos{\theta}}{\sin{\theta}} \right)_{, \theta} + \frac{1}{r}\frac{-r}{u^{2}} - \left(\frac{\cos{\theta}}{\sin{\theta}} \right)^{2}\\
& = \left(1 - \frac{1}{u^{2}} \right) + r \frac{u_{, r}}{u^{3}} - r\frac{v_{, r}}{v}\frac{1}{u^{2}}
\end{align*}

\item
\begin{align*}
\Ric_{\theta\phi} & = \Ric_{\phi \theta} = 0
\end{align*}

\item
\begin{align*}
& \Ric_{\phi \phi} = \Gamma_{tr}^{t}\Gamma_{\phi \phi}^{r} + \Gamma_{\phi \phi, r}^{r} + \Gamma_{rr}^{r}\Gamma_{\phi \phi}^{r} - \Gamma_{\phi \phi}^{r}\Gamma_{\phi r}^{\phi} + \Gamma_{\phi \phi, \theta}^{\theta} + \Gamma_{\theta r}^{\theta}\Gamma_{\phi \phi}^{r} - \Gamma_{\phi \phi}^{\theta}\Gamma_{\phi \theta}^{\phi}\\
& =  \frac{v_{, r}}{v}\frac{-r \sin^{2}{\theta}}{u^{2}} + \left(\frac{-r\sin^{2}{\theta}}{u^{2}} \right)_{, r} + \frac{u_{, r}}{u}\frac{-r\sin^{2}{\theta}}{u^{2}} + (-\sin{\theta}\cos{\theta})_{,\theta} + \sin{\theta}\cos{\theta} \frac{\cos{\theta}}{\sin{\theta}}\\
& = \left(1 - \frac{1}{u^{2}} \right)\sin^{2}{\theta} + r\sin^{2}{\theta} \frac{u_{, r}}{u^{3}} - r\sin^{2}{\theta}\frac{v_{, r}}{v} \frac{1}{u^{2}} \end{align*}

\end{enumerate}

Now we can compute the scalar curvature $R$ of this metric as follows:
\begin{align*}
R & = \text{tr}_{g}\Ric = g^{tt}\Ric_{tt} + g^{rr}\Ric_{rr} + g^{\theta \theta} \Ric_{\theta \theta} + g^{\phi \phi}\Ric_{\phi \phi} \tag{$g$ is diagonal}\\
& = \frac{-1}{v^{2}}\left[\frac{1}{u^{2}} \left(vv_{, rr}  + \frac{2}{r}vv_{, r}\right)  - \frac{u_{, r}}{u^{3}} vv_{, r} + \frac{v_{, t}}{v}\frac{u_{, t}}{u} - \frac{u_{, tt}}{u}\right] \\
& + \frac{1}{u^{2}} \left[-\frac{v_{, rr}}{v} + \frac{2}{r}\frac{u_{, r}}{u} - \frac{v_{, t}}{v^{3}}uu_{, t} + \frac{1}{v}\left(\frac{u u_{, tt}}{v}  + \frac{v_{, r}u_{, r}}{u}\right)  \right] \\
& + \frac{2}{r^{2}} \left[\left(1 - \frac{1}{u^{2}} \right)  + r\frac{u_{, r}}{u^{3}} - r\frac{v_{, r}}{v}\frac{1}{u^{2}}\right] \\
& = -\frac{1}{u^{2}} \frac{v_{, rr}}{v} - \frac{2}{r} \frac{1}{u^{2}} \frac{v_{, r}}{v} + \frac{u_{, r}}{u^{3}}\frac{v_{, r}}{v} - \frac{u_{, t}}{u}\frac{v_{, t}}{v^{3}} + \frac{u_{, tt}}{u}\frac{1}{v^{2}} - \frac{1}{u^{2}}\frac{v_{, rr}}{v} + \frac{2}{r}\frac{u_{, r}}{u^{3}} + \frac{u_{, r}}{u^{3}}\frac{v_{, r}}{v}\\
& - \frac{u_{, t}}{u}\frac{v_{, t}}{v^{3}} + \frac{u_{, tt}}{u}\frac{1}{v^{2}} + \frac{2}{r^{2}}\left(1 - \frac{1}{u^{2}} \right) + \frac{2}{r}\frac{u_{, r}}{u^{3}} - \frac{2}{r}\frac{1}{u^{2}}\frac{v_{, r}}{v} \\
& = -\frac{2}{u^{2}}\frac{v_{, rr}}{v} + 2\frac{u_{, r}}{u^{3}}\frac{v_{, r}}{v} - 2\frac{u_{, t}}{u}\frac{v_{, t}}{v^{3}} + 2\frac{u_{, tt}}{u}\frac{1}{v^{2}} + \frac{4}{r}\frac{u_{, r}}{u^{3}} - \frac{4}{r}\frac{1}{u^{2}}\frac{v_{, r}}{v} + \frac{2}{r^{2}}\left(1 - \frac{1}{u^{2}} \right)
\end{align*}
That is
\begin{equation}\label{scalar curvature explicit}
R = -\frac{2}{u^{2}}\frac{v_{, rr}}{v} + 2\frac{u_{, r}}{u^{3}}\frac{v_{, r}}{v} - 2\frac{u_{, t}}{u}\frac{v_{, t}}{v^{3}} + 2\frac{u_{, tt}}{u}\frac{1}{v^{2}} + \frac{4}{r}\frac{u_{, r}}{u^{3}} - \frac{4}{r}\frac{1}{u^{2}}\frac{v_{, r}}{v} + \frac{2}{r^{2}}\left(1 - \frac{1}{u^{2}} \right).
\end{equation}
The Einstein curvature tensor $G$ is given by
\begin{equation}\label{Einstein curvature tensor}
G = \Ric - \frac{1}{2}Rg.
\end{equation}

In our case of the spacetime metric $(\ref{spherically symmetric metrics satisfying radial mean curvature vector condition again again})$, we compute the components of $G$ with respect to the basis $\{\p_{t}, \p_{r}, \p_{\theta}, \p_{\phi}\}$:
\begin{enumerate}
\item
\begin{align}
G_{tt} & = G(\p_{t}, \p_{t}) = \Ric_{tt} - \frac{1}{2}Rg_{tt} =  \frac{1}{u^{2}}\left(vv_{, rr} + \frac{2}{r}vv_{, r} \right) - \frac{u_{, r} }{u^{3}} v v_{, r} \notag\\
&  +\frac{1}{u}\left(\frac{u_{, t}v_{, t}}{v} - u_{, tt} \right) + \frac{v^{2}}{2}R \notag \\
& = \frac{2}{r}\frac{u_{, r}}{u^{3}} v^{2} + \frac{1}{r^{2}}v^{2}\left(1 - \frac{1}{u^{2}} \right) \label{Gtt}
\end{align}

\item
\begin{align}
G_{tr} = G_{rt} = \Ric_{tr} - \frac{1}{2}Rg_{tr} = \Ric_{tr} = \frac{2}{r}\frac{u_{, t}}{u} \label{Gtr}
\end{align}

\item 
\begin{align*}
G_{t\theta}  = G_{\theta t} = \Ric_{t\theta} - \frac{1}{2}Rg_{t\theta} = \Ric_{t\theta} = 0
\end{align*}

\item 
\begin{align*}
G_{t\phi}  = G_{\phi t} = \Ric_{t\phi} - \frac{1}{2}Rg_{t\phi} = \Ric_{t\phi} = 0
\end{align*}

\item
\begin{align}
G_{rr} & = \Ric_{rr} - \frac{1}{2}Rg_{rr}  = -\frac{v_{, rr}}{v} + \frac{2}{r}\frac{u_{, r}}{u} - \frac{v_{, t}}{v^{3}}uu_{, t} + \frac{1}{v}\left(\frac{u u_{, tt}}{v}  - \frac{v_{, r}u_{, r}}{u}\right)
 - \frac{u^{2}}{2}R \notag \\
 & = \frac{2}{r}\frac{v_{, r}}{v} - \frac{u^{2}}{r^{2}} + \frac{1}{r^{2}}\label{Grr}
\end{align}

\item
\begin{align*}
G_{r\theta} & = G_{\theta r} = \Ric_{r\theta} - \frac{1}{2}Rg_{r\theta} = \Ric_{r\theta} = 0
\end{align*}

\item
\begin{align*}
G_{r\phi} & = G_{\phi r} = \Ric_{r\phi} - \frac{1}{2}Rg_{r\phi} = \Ric_{\phi} = 0
\end{align*}

\item
\begin{align}
G_{\theta \theta} & = \Ric_{\theta \theta} - \frac{1}{2}Rg_{\theta \theta} = \left(1 - \frac{1}{u^{2}} \right) + r \frac{u_{, r}}{u^{3}} - r\frac{v_{, r}}{v}\frac{1}{u^{2}} - \frac{r^{2}}{2}R \notag \\
& = \frac{r^{2}}{u^{2}}\frac{v_{, rr}}{v} - r^{2}\frac{u_{, r}}{u^{3}}\frac{v_{, r}}{v} + r^{2}\frac{u_{, t}}{u}\frac{v_{, t}}{v^{3}} - r^{2}\frac{u_{, tt}}{u}\frac{1}{v^{2}} - r\frac{u_{, r}}{u^{3}} + r\frac{1}{u^{2}}\frac{v_{, r}}{v} \label{Gthetatheta}
\end{align}

\item
\begin{align}
G_{\phi \phi} & = \Ric_{\phi \phi} - \frac{1}{2}Rg_{\phi \phi} = \sin^{2}{\theta}\cdot \Ric_{\theta \theta} - \frac{1}{2}Rr^{2}\sin^{2}{\theta} = \sin^{2}{\theta}\cdot G_{\theta \theta} \notag \\
& = \sin^{2}{\theta} \left(\frac{r^{2}}{u^{2}}\frac{v_{, rr}}{v} - r^{2}\frac{u_{, r}}{u^{3}}\frac{v_{, r}}{v} + r^{2}\frac{u_{, t}}{u}\frac{v_{, t}}{v^{3}} - r^{2}\frac{u_{, tt}}{u}\frac{1}{v^{2}} - r\frac{u_{, r}}{u^{3}} + r\frac{1}{u^{2}}\frac{v_{, r}}{v} \right) \label{Gphiphi}
\end{align}

\end{enumerate}

\section{Calculations in Inverse Mean Curvature Vector Flow Coordinates}
\subsection{Determinant of the Spacetime Metric and its Inverse in Inverse Mean Curvature Vector Flow Coordinates}
\label{subsec:determinant and inverse of spacetime metric}

Given a matrix $A$, its $(i, j)$th cofactor is $C_{ij}: = (-1)^{i + j}M_{ij}$, where $M_{ij}$ is determinant of the matrix obtained by deleting the $i$th row and the $j$th column of $A$. The adjoint matrix $\adj(A)$ is defined as $\adj(A)_{ij}: = C_{ji} = (-1)^{j+i}M_{ji}$. If $A$ is invertible, then $A^{-1} = \frac{1}{\det(A)}\adj(A)$.

Let $A$ be the matrix representation of the spacetime metric $\gbar$ as in (\ref{equ:local representation of gbar in IMCVF coordinates step 1}). 
\begin{align}
& \adj(A)_{11} = C_{11} = M_{11} = u^{2}(ab - c^{2}) \\
& \adj(A)_{12} = C_{21} = -M_{21} = -d(ab - c^{2}) \\
& \adj(A)_{13} =  C_{31} = M_{31} = u^{2}(cf - be) \\
& \adj(A)_{14} = C_{41} = -M_{41} = u^{2}(ce - af) \\
& \adj(A)_{22} = C_{22} = M_{22} = -v^{2}(ab - c^{2}) + f(ce - af) + e(cf - be) \\
& \adj(A)_{23} = C_{32} = -M_{32} = -d(cf - be) \\
& \adj(A)_{24} = C_{42} = M_{42} = -d(ce - af) \\
& \adj(A)_{33} = C_{33} = M_{33} = -u^{2}v^{2}b - u^{2}f^{2} - bd^{2} \\
& \adj(A)_{34} = C_{43} = -M_{43} = u^{2}v^{2}c + u^{2}ef + cd^{2} \\
& \adj(A)_{44} = C_{44} = M_{44} = -u^{2}v^{2}a - u^{2}e^{2} - ad^{2}
\end{align}

\begin{align}
\det(\gbar) & = A_{21}C_{21} + A_{22}C_{22} = -dM_{21} + u^{2}M_{22}\notag\\
& = -d^{2}(ab - c^{2}) + u^{2}\big[-v^{2}(ab - c^{2}) + f(ce - af) + e(cf - be)\big]\notag\\
& = (-u^{2}v^{2} - d^{2})(ab - c^{2}) + eu^{2}(cf - be) + fu^{2}(ce - af)
\end{align}

These prove Lemma \ref{lemma:inverse and det of spacetime metric in IMCVF coordinates}.

\subsection{Computation of $\<\n,\n\>$}
\label{subsec:inner product of n with n}

\begin{lemma}
\begin{equation}
\<\textbf{n}, \textbf{n}\> = \frac{\text{det}(\overline{g})}{u^{2}(ab - c^{2})} = \frac{\text{det}(\overline{g})}{u^{2}\text{det}(g_{S})} = :\frac{|\gbar|}{u^{2}|g_{S}|},
\end{equation}
where we set $|g_{S}|: = \text{det}(g_{S})$.
\end{lemma}

\begin{proof}
We simply compute that
\begin{align*}
& \<\n, \n\> = \bigg\< \frac{\p}{\p t} + \frac{-d}{u^{2}}\frac{\p}{\p r} + \frac{cf - be}{ab - c^{2}}\frac{\p}{\p \theta} + \frac{ce - af}{ab - c^{2}}\frac{\p}{\p \phi}\,,\,  \frac{\p}{\p t} + \frac{-d}{u^{2}}\frac{\p}{\p r} + \frac{cf - be}{ab - c^{2}}\frac{\p}{\p \theta}\\
& + \frac{ce - af}{ab - c^{2}}\frac{\p}{\p \phi}\bigg\> \\
& = \<\ppt, \ppt\> + \frac{-d}{u^{2}}\<\ppt, \ppr\> + \Red\frac{cf - be}{ab - c^{2}}\<\ppt, \pptheta\> + \frac{ce - af}{ab - c^{2}}\<\ppt, \ppphi\> \Black \\
& + \frac{-d}{u^{2}}\<\ppr, \ppt\> + \frac{d^{2}}{u^{4}}\<\ppr, \ppr\> + 0 + 0\\
& + \Red\frac{cf - be}{ab - c^{2}}\<\pptheta, \ppt\>\Black + 0+  \YellowOrange\left(\frac{cf - be}{ab - c^{2}}\right)^{2}\<\pptheta, \pptheta\>\Black + \Green\frac{(cf - be)(ce - af)}{(ab - c^{2})^{2}}\<\pptheta, \ppphi\>\Black\\
& + \Red\frac{ce - af}{ab - c^{2}}\<\ppphi, \ppt\>\Black + 0 + \Green\frac{(ce - af)(cf - be)}{(ab - c^{2})^{2}}\<\ppphi, \pptheta\>\Black + \YellowOrange\left( \frac{ce - af}{ab - c^{2}}\right)^{2}\<\ppphi, \ppphi\>\Black\\
& = -v^{2} + \frac{-d^{2}}{u^{2}} + \Red\frac{1}{ab - c^{2}}( 2(cf - be)e + 2(ce - af)f)\Black + \Blue \frac{-d^{2}}{u^{2}}\Black + \Blue\frac{d^{2}}{u^{2}} \Black \\
& + 2 \Green\frac{(ce - af)(cf - be)c}{(ab - c^{2})^{2}} \Black + \YellowOrange\left(\frac{cf - be}{ab - c^{2}}\right)^{2} a\Black + \YellowOrange\left( \frac{ce - af}{ab - c^{2}}\right)^{2} b\Black\\
& =: -v^{2} + \frac{-d^{2}}{u^{2}} + \Red\frac{1}{ab - c^{2}}( 2cef - 2be^{2} + 2cef - 2af^{2})\Black + (A)
\end{align*}
where $(A)$ is 
\begin{align*}
A & = 2 \Green\frac{(ce - af)(cf - be)c}{(ab - c^{2})^{2}}\Black + \YellowOrange\left(\frac{cf - be}{ab - c^{2}}\right)^{2} a\Black + \YellowOrange\left( \frac{ce - af}{ab - c^{2}}\right)^{2} b\Black\\
& = \frac{1}{(ab - c^{2})^{2}}\Big( \Green2c^{3}ef - 2bc^{2}e^{2} - 2ac^{2}f^{2} + 2abcef\Black + \YellowOrange ac^{2}f^{2} + ab^{2}e^{2} - 2abcef + bc^{2}e^{2}\Black\\
&  \YellowOrange + a^{2}bf^{2} - 2abcef\Black\Big)\\
& = \frac{1}{(ab - c^{2})^{2}}\Big(2c^{3}ef - 2abcef - bc^{2}e^{2} -  ac^{2}f^{2}+ ab^{2}e^{2} + a^{2}bf^{2} \Big)\\
& =  \frac{1}{(ab - c^{2})^{2}}\Big( 2cef(c^{2} - ab) -c^{2}(af^{2} + be^{2}) + ab(be^{2} + af^{2}) \Big)\\
& = \frac{1}{(ab - c^{2})^{2}}\Big(-2cef(ab - c^{2}) + (af^{2} + be^{2})(ab - c^{2}) \Big)\\
& = \frac{1}{ab - c^{2}}( af^{2} + be^{2} - 2cef ).
\end{align*}
Plug $(A)$ back into the above, we get:
\begin{align}
\<\n, \n\> & = -v^{2} + \frac{-d^{2}}{u^{2}} + \frac{1}{ab - c^{2}}( 2cef - 2be^{2} + 2cef - 2af^{2})\notag\\
& + \frac{1}{ab - c^{2}}( af^{2} + be^{2} - 2cef )\notag\\
& = -v^{2} + \frac{-d^{2}}{u^{2}} + \frac{1}{ab - c^{2}} (2cef - 2be^{2} + 2cef - 2af^{2} + af^{2} + be^{2} - 2cef)\notag\\
& = \frac{-u^{2}v^{2} - d^{2}}{u^{2}} + \frac{1}{ab - c^{2}} (2cef - be^{2} - af^{2}) \notag\\
& = \frac{1}{u^{2}(ab - c^{2})}\Big((-u^{2}v^{2} - d^{2})(ab - c^{2}) + u^{2}(2cef - be^{2} - af^{2})\Big)\\
& = \frac{|\gbar|}{u^{2}|g_{S}|},
\end{align}
by Equation (\ref{equ:det of spacetime metric in IMCVF coordinates}).
\end{proof}

\subsection{Laplacian along $S_{t,r}$}\label{appen:laplacian of coordinate sphere}
Note that the laplacian of $d$ on the surface $(S_{t,r}, g_{S})$ is:
\begin{align}
& \Delta_{g_{S}}d = \divg_{g_{S}}(\na_{g_{S}} d) = \frac{1}{\sqrt{|g_{S}|}}\sum_{i = 1}^{2}\ppxi\left( (\na_{g_{S}} d)^{i}\sqrt{|g_{S}|}\right) = \frac{1}{\sqrt{|g_{S}|}}\sum_{i = 1}^{2}\ppxi\left((g_{S}^{ij}d_{,j})\sqrt{|g_{S}|}\right)\notag\\
& = \frac{1}{r^{2}\sin{\theta}}\left(\pptheta\left( (g_{S}^{\theta \theta}d_{,\theta} + g_{S}^{\theta \phi}d_{,\phi}) r^{2}\sin{\theta} \right) + \ppphi\left((g_{S}^{\phi \theta}d_{,\theta} + g_{S}^{\phi \phi} d_{,\phi}) r^{2}\sin{\theta} \right) \right)\notag\\
& = \frac{1}{\Blue r^{2}\sin{\theta} \Black}\Big((g_{S,\theta}^{\theta \theta}d_{,\theta} + \Red g_{S}^{\theta \theta}d_{,\theta \theta}\Black  + g_{S, \theta}^{\theta \phi} d_{,\phi} + \Red g_{S}^{\theta \phi} d_{,\phi \theta}\Black  )\Blue r^{2}\sin{\theta}\Black + (g_{S}^{\theta \theta}d_{,\theta} + g_{S}^{\theta \phi}d_{,\phi}) \Blue r^{2}\cos{\theta}\Black \notag\\
& + (g_{S,\phi}^{\phi \theta}d_{,\theta} + \Red g_{S}^{\phi \theta}d_{,\theta \phi}\Black + g_{S,\phi}^{\phi \phi} d_{,\phi} + \Red g_{S}^{\phi \phi} d_{,\phi \phi}\Black)\Blue r^{2}\sin{\theta}\Black + 0  \Big)\notag\\
& = (\Red g_{S}^{\theta \theta}d_{,\theta \theta} + 2g_{S}^{\theta \phi} d_{,\phi \theta} +  g_{S}^{\phi \phi} d_{,\phi \phi}\Black)\tag{second order derivative in $d$} \notag\\
& + (g_{S,\theta}^{\theta \theta}d_{,\theta} + g_{S, \theta}^{\theta \phi} d_{,\phi}+ g_{S,\phi}^{\phi \theta}d_{,\theta}+ g_{S,\phi}^{\phi \phi} d_{,\phi}) +  (g_{S}^{\theta \theta}d_{,\theta} + g_{S}^{\theta \phi}d_{,\phi})\cot{\theta}
\end{align}

That is:
\begin{align} \Delta_{g_{S}}d &
= (\Red g_{S}^{\theta \theta}d_{,\theta \theta} + 2g_{S}^{\theta \phi} d_{,\phi \theta} +  g_{S}^{\phi \phi} d_{,\phi \phi}\Black) + \Peach (g_{S,\theta}^{\theta \theta} + g_{S,\phi}^{\phi \theta} + g_{S}^{\theta \theta}\cot{\theta})d_{,\theta} \Black\notag\\
& + (g_{S,\theta}^{\theta \phi} + g_{S,\phi}^{\phi \phi} + g_{S}^{\theta \phi}\cot{\theta})d_{,\phi} \Black \notag\\
& = \frac{1}{|g_{S}|}\Big((\Red bd_{,\theta \theta} - 2cd_{,\phi \theta} + ad_{,\phi \phi}\Black) + \Peach(b_{,\theta} - 2b\cot{\theta} - c_{,\phi} + b\cot{\theta})d_{,\theta}\Black\notag\\
&  + \Peach(-c_{,\theta} + 2c\cot{\theta} + a_{,\phi} - c\cot{\theta})d_{,\phi}\Black\Big) \notag\\
& = \frac{1}{|g_{S}|}\Big((\Red bd_{,\theta \theta} - 2cd_{,\phi \theta} + ad_{,\phi \phi}\Black) + \Peach(b_{,\theta} - b\cot{\theta} - c_{,\phi})d_{,\theta} + (-c_{,\theta} + c\cot{\theta} + a_{,\phi})d_{,\phi}\Black\Big) \label{equ:laplacian of d on the surface str}
\end{align}

For the second to the last equality above, we have used the following computations:
\begin{equation}\label{equ:derivative of det of gs}
|g_{S}|_{,\theta} = (r^{4}\sin^{2}{\theta})_{,\theta} = r^{4}2\sin{\theta}\cos{\theta} = 2r^{4}\sin^{2}{\theta}\cot{\theta} = 2|g_{S}|\cot{\theta}.
\end{equation}
and
\begin{equation}\label{equ:derivative of inverse of gs}
g_{S, \theta}^{\theta\theta} = \left(\frac{b}{|g_{S}|} \right)_{\theta} = \frac{b_{,\theta}|g_{S}| - b|g_{S}|_{,\theta}}{|g_{S}|^{2}} = \frac{1}{|g_{S}|}(b_{,\theta} - 2b\cot{\theta});\quad  g_{S,\phi}^{\theta \theta} = \left(\frac{b}{|g_{S}|} \right)_{\phi} = \frac{b_{,\phi}}{|g_{S}|}.
\end{equation}
Similar for the other derivatives of the inverse of $g_{S}$.

\section{First Variation of Area}\label{append:first variation of area}
Let $\Sigma^{n - 1}$ be an embedded closed (compact without boundary) hypersurface in a Riemannian manifold $(M^{n}, g, \overline{\na})$. Endow $\Sigma$ with the induced metric. We consider a variation of $\Sigma$ as follows:
\begin{equation}
F: \Sigma \times (-\delta, \delta) \longrightarrow M, \quad \delta > 0,
\end{equation}
such that for all $x\in \Sigma_{t}: = F(\Sigma, t)$, and $t \in (-\delta, \delta)$,
\begin{equation}
\frac{\p }{\p t} F(x, t) = \eta(x, t)\nu(x, t),
\end{equation}
where $\eta$ is a smooth function $\eta\in C^{\infty}(\Sigma\times (-\delta, \delta))$, and $\nu(x, t)$ is the unit outward normal vector to $\Sigma_{t}$ at $(x, t)$. Therefore the variational vector fields along each surface $\Sigma_{t}$ is $\frac{\p}{\p t} = \eta_{t}\cdot \nu_{t}$. Let $g_{t}$ be the induced metric on $\Sigma_{t}$, and let $\na^{t}$ be the associated Levi-Civita connection. Let $d\sigma_{t}$ be the corresponding $(n - 1)$-volume form on $\Sigma_{t}$, and $A_{t}$ the $(n - 1)$-volume.  Let $V_{t}$ be the $n$-volume enclosed by $\Sigma_{t}$. We shall refer to $A_{t}$ as the \emph{area} of $\Sigma_{t}$, and $V_{t}$ the \emph{volume}, in analogy to the case where $\Sigma_{t}$ are surfaces in a $3$-dimensional manifold.
Let $\II_{t}$ and $H_{t}: = \text{tr}_{g_{t}} \II$ be the second fundamental form and the mean curvature of $\Sigma$ with respect to $\nu(x, t)$ respectively. We first compute the variation of $d\sigma_{t}$. Let $\{U; x_{1}, x_{2}, \cdots, x_{n - 1}\}$ be a local coordinate chart of $\Sigma$, then $\Sigma_{t}$ can be locally parametrized as $\{x_{1}, x_{2}, \cdots, x_{n - 1}, t\}$ with each fixed $t$. Let $g_{t} = (g_{t})_{ij}dx^{i}dx^{j}$ be the local representation of the metric on $\Sigma_{t}$, $i, j  = 1, 2, \cdots, n-1$.
\begin{align}
\frac{\p}{\p t} d\sigma_{t} & = \frac{\p}{\p t}\sqrt{\text{det}(g_{t})}dx^{1}\wedge dx^{2} \cdots \wedge dx^{n-1} \notag \\
& = \frac{1}{2}\frac{1}{\sqrt{\text{det}(g_{t})}}  \text{det}(g_{t}) \cdot \text{trace} \left(g^{-1}_{t} \frac{\p}{\p t}g_{t}\right) dx^{1}\wedge dx^{2} \wedge \cdots \wedge dx^{n - 1} \label{derivative of determinant}\\
& = \frac{1}{2}\sqrt{\text{det}(g_{t})} \text{trace} \left(g^{-1}_{t} \frac{\p}{\p t}g_{t}\right) dx^{1}\wedge dx^{2} \wedge \cdots \wedge dx^{n - 1}  \notag 
\end{align}
where equation $(\ref{derivative of determinant})$ follows from the identity
\begin{equation}
\frac{d}{dt}\text{det}(A) = \text{det}(A)\text{trace}\left(A^{-1}\frac{d}{dt}A \right),
\end{equation}
for any square matrix $A$ with entries functions of $t$.
Now
\begin{align}
\frac{\p}{\p t}(g_{t})_{ij} & = \frac{\p}{\p t}\<\ppxi, \ppxj\> = \<\overlinenat \ppxi, \ppxj\> + \<\ppxi, \overlinenat \ppxj\> \notag \\
& = \<\overlinenaxi \ppt, \ppxj\> + \<\ppxi, \overlinenaxj \ppt\> \tag{$\overline{\na}$ is torsion free}\notag \\
& = \<\overlinenaxi \eta \nu, \ppxj\> + \<\ppxi, \overlinenaxj \eta \nu\>\notag \\
& = 2\eta \II_{t}\left(\ppxi, \ppxj\right).  \label{variation of metric} 
\end{align}
where the last identity follows from the fact that $\nu$ is a normal to the surfaces. Notice that Equation (\ref{variation of metric}) implies that the first derivative of the metric (along the variational vector fields) is given by the second fundamental form. Now plug (\ref{variation of metric}) into (\ref{derivative of determinant}):
\begin{equation}\label{first variation of area form}
\frac{\p}{\p t} d\sigma_{t} = \frac{1}{2}\sqrt{\text{det}(g_{t})} \trace\left(g_{t}^{-1}\cdot 2 \eta\cdot  \II_{t}(\ppxi, \ppxj) \right) dx^{1}\wedge dx^{2} \wedge \cdots \wedge dx^{n-1}  =  H_{t} \eta d\sigma_{t}
\end{equation}
Therefore
\begin{equation}\label{equ:first variation of area}
\frac{\p}{\p t}A_{t} = \int_{\Sigma_{t}} \frac{\p}{\p t} d\sigma_{t} = \int_{\Sigma_{t}}H_{t}(x)\eta(x, t)\, d\sigma_{t}(x).
\end{equation}

\section{Second Variation of Area}\label{append:second variation of area}
Now we compute the first variation of mean curvature $H_{t}$, which gives rise to the second derivative of area. Recall that $H_{t} = g_{t}^{ij}(\II_{t})_{ij}$ in local coordinates, $i, j = 1, 2, \cdots, n-1$. Thus
\begin{equation}
\frac{\p}{\p t}H_{t} =  \frac{\p}{\p t}g_{t}^{ij}(\II_{t})_{ij} + g_{t}^{ij} \frac{\p}{\p t} (\II_{t})_{ij}
\end{equation}
Since $0 = \frac{\p}{\p t} (g_{t} g_{t}^{-1})  = (\frac{\p}{\p t}g_{t}) g_{t}^{-1} + g_{t}(\frac{\p}{\p t}g_{t}^{-1})$, we have $\frac{\p}{\p t}g_{t}^{-1}  = -g_{t}^{-1} (\frac{\p}{\p t} g_{t}) g_{t}^{-1}$. Thus the first term in the above becomes
\begin{equation}
\frac{\p}{\p t}g_{t}^{ij} (\II_{t})_{ij} = - g_{t}^{ik}\left( \frac{\p}{\p t}(g_{t})_{kl} \right) g_{t}^{lj}(\II_{t})_{ij} = -g_{t}^{ik} 2\eta (\II_{t})_{kl} g_{t}^{lj} (\II_{t})_{ij} = -2\eta ||\II_{t}||^{2}
\end{equation}

We now compute the derivative of the second fundamental form.
\begin{align*}
& \frac{\p}{\p t} (\II_{t})_{ij} = \frac{\p}{\p t}\<\overlinenaxi \nu, \ppxj\>\\
& = \<\overlinenat \overlinenaxi \nu, \ppxj\> + \<\overlinenaxi \nu, \overlinenat \ppxj\>\\
& = \<\overlinenaxi \overlinenat \nu, \ppxj\> + \<(\overlinenat \overlinenaxi - \overlinenaxi \overlinenat - \overline{\na}_{[\ppxi, \ppt]}) \nu, \ppxj\> + \<\overlinenaxi \nu, \overlinenat \ppxj\>\\
& = \<\overlinenaxi \overlinenat \nu, \ppxj\> + \<R_{g}(\ppt, \ppxi)\nu, \ppxj\> + \<\overlinenaxi \nu,  \overlinenaxj \ppt\> \\
& = \<\overlinenaxi \left( -\na_{\Sigma_{t}} \eta \right), \ppxj\> + \eta \<R(\nu, \ppxi)\nu, \ppxj\> + \eta\<\overlinenaxi \nu, \overlinenaxj \nu\>
\end{align*}
where we have used two lemmas, which will be proved below: 

\begin{lemma}\label{lemma: covariant derivative of nu is tangential}
$\overlinenaxi \nu$ is tangential, i.e., $\<\overlinenaxi \nu, \nu\> = 0$.
\end{lemma}

\begin{lemma}\label{lemma: surface gradient}
$\overlinenat \nu = -\na_{\Sigma_{t}} \eta$, where $\na_{\Sigma_{t}}$ is the surface gradient on $\Sigma_{t}$.
\end{lemma}

Therefore
\begin{align*}
g_{t}^{ij} \frac{\p}{\p t}(\II_{t})_{ij} & = -\Delta_{\Sigma_{t}} \eta - \eta  g_{t}^{ij}\<R(\nu, \ppxi), \ppxj, \nu\>+ \eta ||\II_{t}||^{2}\\
& = -\Delta_{\Sigma_{t}} \eta - \eta  g^{ij}\<R(\nu, \ppxi), \ppxj, \nu\>+ \eta ||\II_{t}||^{2}\tag{ambient metric $g$ trace}\\
& =  -\Delta_{\Sigma_{t}} \eta - \eta  \Ric_{g}(\nu, \nu)+ \eta ||\II_{t}||^{2}.
\end{align*}
where we used
\begin{lemma}\label{lemma: trace of surface gradient}
$g_{t}^{ij} \<\overlinenaxi \left(\na_{\Sigma_{t}} \eta \right), \ppxj\> = \Delta_{\Sigma_{t}}\eta$.
\end{lemma}
and
\begin{lemma}\label{lemma: norm of second fundamental form}
$g_{t}^{ij} \<\overlinenaxi \nu, \overlinenaxj \nu\> = ||\II_{t}||^{2}$.
\end{lemma}

Combining above, we get
\begin{equation}\label{equ:first variation of mean curvature}
\frac{\p}{\p t}H_{t} = -\Delta_{\Sigma_{t}} \eta - \eta \left( \Ric_{g}(\nu, \nu) + ||\II_{t}||^{2}\right) =: L_{\Sigma_{t}},
\end{equation}
where $L_{\Sigma_{t}}$ is called the \emph{stability operator} of $\Sigma_{t}$. The second variation of area is then given by:
\begin{equation}\label{second variation of area}
\frac{\p^{2}}{\p t^{2}}A_{t} = \int_{\Sigma} \eta(x, t) (L_{\Sigma_{t}}\eta)(x, t)\, d\sigma_{t} (x) + H_{t}(x)\left(\frac{\p}{\p t}\eta(x, t)\right)\, d\sigma_{t}(x) + H_{t}^{2}\eta(x, t)^{2}\, d\sigma_{t}(x).
\end{equation}

Now we verity the above lemmas. 
\begin{proof}[Proof of Lemma \ref{lemma: covariant derivative of nu is tangential}]
Since $\nu$ is the unit outward normal vector field, we have
\begin{equation}
0 = \ppxi\<\nu, \nu\> = 2\<\overlinenaxi \nu, \nu\>.
\end{equation}
Thus $\overlinenaxi \nu$ is tangential. Similarly, $\overlinenat \nu$ is also tangential.
\end{proof}

\begin{proof}[Proof of Lemma \ref{lemma: surface gradient}]
Recall that the gradient of a smooth function $\eta$ along the surface $\Sigma_{t}$ is defined as
\begin{equation}
\na_{\Sigma_{t}} \eta: = \na \eta - \<\na \eta, \nu\>\nu,
\end{equation}
that is, the tangential component of the gradient with respect to the ambient metric. For any point $p\in \Sigma_{t}$, choose geodesic normal coordinates $\{U; e_{1}, e_{2}, \cdots, e_{n}\}$ around $p$ such that $e_{1}, e_{2}, \cdots, e_{n-1}$ span $T_{p}\Sigma_{t}$, and $e_{n} = \nu$. Since $\overline{\na}_{e_{i}}\nu$ is tangential, it suffices to show that $\<\overlinenat \nu, e_{i}\> (p) = \<-\na_{\Sigma_{t}} \eta, e_{i}\> (p)$, for $i = 1, 2, \cdots, n-1$. Indeed:
\begin{align*}
\<\overlinenat \nu, e_{i}\> (p) & = -\<\nu, \overlinenat e_{i}\> (p) = -\<\nu, \overline{\na}_{e_{i}} \ppt\> (p) \tag{$\overline{\na}$ is torsion free}\\
& = -  e_{i}(\eta)\<\nu,\nu\> (p) - \eta\<\nu, \overline{\na}_{e_{i}} \nu\>(p)\\
& = -e_{i}(\eta)(p) \tag{$\overline{\na}_{e_{i}}\nu$ is tangential} \\
& = -\<\sum_{j = 1}^{n-1}e_{j}(\eta)e_{j}, e_{i})\>(p) \\
& = \<-\na_{\Sigma_{t}} \eta, e_{i}\>(p).
\end{align*}
Since $p$ is arbitrary, $\overlinenat \nu = -\na_{\Sigma_{t}} \eta$, as desired.
\end{proof}

\begin{proof}[Proof of Lemma \ref{lemma: trace of surface gradient}]
First note that $g_{t}^{ij}\<\overlinenaxi (\na_{\Sigma_{t}}\eta), \ppxj\> = g_{t}^{ij}\<\na^{t}_{\ppxi} (\na_{\Sigma_{t}}\eta), \ppxj\>$. But 
\begin{equation}
g_{t}^{ij}\<\na^{t}_{\ppxi} (\na_{\Sigma_{t}}\eta), \ppxj\> = tr_{g_{t}}(\na^{t}(\na_{\Sigma_{t}}\eta) ) = \Div_{g_{t}} (\na_{\Sigma_{t}}\eta) = \Delta_{\Sigma_{t}}\eta,
\end{equation}
where $\na^{t}(\na_{\Sigma_{t}}\eta)$ is the covariant derivative of the vector field $\na_{\Sigma_{t}}\eta$, hence is a $(1, 1)$-tensor field.
\end{proof}

\begin{proof}[Proof of Lemma \ref{lemma: norm of second fundamental form}]
Define vector fields $X(i): = \overlinenaxi \nu$ and $Y(j): = \overlinenaxj \nu$. Using local coordinates, we can also write $X(i) = \displaystyle\sum_{k = 1}^{n-1}X(i)^{k}\ppxk$ and $Y(j) = \displaystyle\sum_{l = 1}^{n-1} Y(j)^{l}\ppxl$. Then 
\begin{align*}
||\II_{t}||^{2} & = g_{t}^{ij}g_{t}^{kl}(\II_{t})_{ik}(\II_{t})_{jl} = g_{t}^{ij}g_{t}^{kl}\left\<X(i), \ppxk\right\>\left\<Y(j), \ppxl\right\>\\
& = g_{t}^{ij}g_{t}^{kl}X(i)^{\alpha}(g_{t})_{\alpha k}Y(j)^{\beta}(g_{t})_{\beta l}\\
& = g_{t}^{ij}(g_{t})_{\alpha \beta}X(i)^{\alpha}Y(j)^{\beta}\\
& = g_{t}^{ij}\<X(i), Y(j)\>\\
& = g_{t}^{ij} \<\overlinenaxi \nu, \overlinenaxj \nu\>
\end{align*}
as desired. 
\end{proof}

\section{Transformation Formulae of Ricci and Scalar Metric under Conformal Change of Metrics}

Given Riemannian manifold $(M^{n}, g)$, recall that the Riemann curvature operator $R$ acts on vector fields $X, Y, Z \in \Gamma(TM)$ as follows:
\begin{equation}\label{def:Riemann curvature operator}
R(X, Y)Z = D_{X}D_{Y}Z - D_{Y}D_{X}Z - D_{[X, Y]}Z.
\end{equation}
We can thus define a $(1, 3)$ tensor field $R = R^{k}_{lij} dx^{l} \otimes dx^{i} \otimes dx^{j}\otimes \frac{\p}{\p x^{k}}$, such that 
\[R(\frac{\p}{\p x^{i}}, \frac{\p}{\p x^{j}})\frac{\p}{\p x^{l}} =: R^{k}_{lij}\frac{\p}{\p x^{k}}.\]
One can verify that 
\begin{equation}\label{equ:components of Riemann curvature operator}
R^{k}_{lij} = \frac{\p \Gamma^{k}_{jl}}{\p x^{i}} - \frac{\Gamma^{k}_{il}}{\p x^{j}} + \Gamma^{k}_{im}\Gamma^{m}_{jl} - \Gamma^{k}_{jm}\Gamma^{m}_{il}.
\end{equation}

Using the metric $g$, we can define a new $(0, 4)$ tensor field $Rm: = R_{klij}dx^{k}\otimes dx^{l}\otimes dx^{i}\otimes dx^{j}$ with components $R_{klij}: = g_{km}R^{m}_{lij}$. One verifies that
\[R_{klij} = \<R(\p_{i}, \p_{j})\p_{l}, \p_{k}\>.\]
$Rm$ is called the Riemannian curvature tensor with respect to the metric $g$.
The Ricci curvature in the direction $X\in T_{p}M$ is given by
\begin{equation}
\text{Ric}(X, X): = g^{jl}\<R(X, \frac{\p}{\p x^{j}})\frac{\p}{\p x^{l}}, X\>.
\end{equation}
Therefore the Ricci curvature tensor is given by the (2,4)-contraction of the Riemann curvature tensor, i.e., 
\begin{align}
\text{Ric}_{ij} & = \text{Ric}(\frac{\p}{\p x^{i}}, \frac{\p}{\p x^{j}}) = g^{kl}\<R(\frac{\p}{\p x^{i}}, \frac{\p}{\p x^{k} })\frac{\p}{\p x^{l}}, \frac{\p}{\p x^{j}}\> = g^{kl}R{jlik} = g^{kl}R_{ikjl} = -R^{l}_{ijl}\notag \\
& = \frac{\p\Gamma^{k}_{ij}} {\p x^{k}} - \frac{\p\Gamma^{k}_{ik}}{\p x^{j}} + \Gamma^{k}_{kl}\Gamma^{l}_{ij} - \Gamma^{k}_{jl}\Gamma^{l}_{ik}. \label{equ:Ricci curvature tensor}
\end{align}

The scalar curvature is
\begin{equation}\label{def:scalar curvature tensor}
R: = g^{ij}\text{Ric}_{ij}.
\end{equation}

\begin{proposition}\label{proposition: conformal transformation of Ricci curvature}
Let $(M^{n}, g)$ be a Riemannian manifold of dimension $n$, and $\rho > 0$ is a smooth function on $M$. Consider the new metric $\widetilde{g}: = \rho g$, i.e., $\widetilde{g}$ is a conformal change of $g$. Then the corresponding Ricci curvatures changes in the following way:
\begin{equation}\label{conformal transformation of Ricci curvature}
\color{Black}\widetilde{R}_{ij} = R_{ij} - \frac{n-2}{2} (\log{\rho})_{,ij} + \frac{n-2}{4} (\log{\rho})_{,i} (\log{\rho})_{,j} - \frac{1}{2}g_{ij}\left[\Delta_{g} \log{\rho} + \frac{n-2}{2}|\na \log{\rho}|_{g}^{2} \right]. 
\end{equation}
\end{proposition}

\begin{proof}
We will verity equation $(\ref{conformal transformation of Ricci curvature})$ at the center $p\in M$ of a geodesic normal coordinate neighborhood. Then at $p$, we have $g_{ij} = \delta_{ij}$ and $\Gamma^{k}_{ij} \equiv 0$, and therefore
\begin{equation}
g_{ij,k} = \p_{k}\<\p_{i}, \p_{j}\> = \<\Gamma^{m}_{ki}\p_{m}, \p_{j}\> + \<\p_{i}, \Gamma^{n}_{kj}\p_{n}\> = 0 \,\,\text{at} \,\, p.
\end{equation}

By $(\ref{equ:Ricci curvature tensor})$, the Ricci curvature of $\widetilde{g}$ is given by:
\begin{equation}
\widetilde{R}_{ij} = \widetilde{\Gamma}^{k}_{ij, k} - \widetilde{\Gamma}^{k}_{ik, j} + \widetilde{\Gamma}^{s}_{ij}\widetilde{\Gamma}^{k}_{sk} - \widetilde{\Gamma}^{k}_{ik}\widetilde{\Gamma}^{k}_{sj}.
\end{equation}

Thus at $p$, we have:

\begin{align}
\widetilde{R}_{ij} & = \bigg(\Gamma^{k}_{ij ,k} + \color{red} \frac{1}{2}\delta_{ik}(\log{\rho})_{,jk} + \color{black} \frac{1}{2} \delta_{jk} (\log{\rho})_{, ik} - \frac{1}{2}g_{ij ,k}g^{kl}(\log{\rho})_{,l} - \frac{1}{2}g_{ij}g^{kl}_{,k} (\log{\rho})_{,l} \notag \\
& - \frac{1}{2}g_{ij}g^{kl}(\log{\rho})_{,lk} \bigg) - \bigg(\Gamma^{k}_{ik,j} + \color{red} \frac{1}{2}\delta_{ik}(\log{\rho})_{,kj} +\color{black}  \frac{1}{2}\delta_{kk}(\log{\rho})_{,ij} - \frac{1}{2} g_{ik,j}g^{kl}(\log{\rho})_{,l}\notag\\
&  - \frac{1}{2}g_{ik}g^{kl}_{,j}(\log{\rho})_{,l} - \frac{1}{2}g_{ik}g^{kl}(\log{\rho})_{,ij} \bigg)\notag \\
& + \left(\color{blue} \Gamma^{s}_{ij} + \color{Green} \frac{1}{2}\delta_{is}(\log{\rho})_{,j} + \color{black} \frac{1}{2}\delta_{js}(\log{\rho})_{,i} - \frac{1}{2}g_{ij}g^{sl}(\log{\rho})_{,l} \right)\times\notag \\
& \left(\color{blue} \Gamma^{k}_{sk} + \color{Green} \frac{1}{2}\delta_{sk}(\log{\rho})_{,k} + \color{black} \frac{1}{2}\delta_{kk}(\log{\rho})_{,s} - \frac{1}{2}g_{sk}g^{kl}(\log{\rho})_{,l} \right)\notag \\
& - \left(\color{blue} \Gamma^{s}_{ik} +  \color{Green}\frac{1}{2}\delta_{is}(\log{\rho})_{,k} + \color{black}  \frac{1}{2}\delta_{ks}(\log{\rho})_{,i} - \frac{1}{2}g_{ik}g^{sl}(\log{\rho})_{,l} \right)\times\notag \\
& \left(\color{blue} \Gamma^{k}_{sj} + \color{Green} \frac{1}{2}\delta_{sk}(\log{\rho})_{,j} + \color{black} \frac{1}{2}\delta_{jk}(\log{\rho})_{,s} - \frac{1}{2}g_{sj}g^{kl}(\log{\rho})_{,l} \right)\notag \\
& = R_{ij} + \frac{2-n}{2}(\log{\rho})_{,ij} - \frac{1}{2}g_{ij}(\log{\rho})_{,kk}  + (\text{product terms}),
\end{align}
where the red terms cancel, green terms cancel and the blue terms vanish at $p$.
\begin{align}
& (\text{product terms}) = \frac{1}{4}\delta_{js}\delta_{kk}(\log{\rho})_{,i}(\log{\rho})_{,s} - \frac{1}{4}\delta_{ks}\delta_{jk}(\log{\rho})_{,i}(\log{\rho})_{,s} \notag\\
& + \frac{1}{4}g_{ij}g_{sk}g^{sl}g^{kl}(\log{\rho})_{,l}(\log{\rho})_{,l} - \frac{1}{4}g_{ik}g_{sj}g^{sl}g^{kl}(\log{\rho})_{,l}(\log{\rho})_{,l} \notag \\
&  - \color{red} \frac{1}{4}g_{jk}g^{kl}(\log{\rho})_{,l}(\log{\rho})_{,i} \color{black} - \frac{n}{4}g_{ij}g^{sj}(\log{\rho})_{,s}(\log{\rho})_{,l}\notag \\
 & + \color{red} g_{jk}g^{kl}(\log{\rho})_{,l}(\log{\rho})_{,i}  \color{black} + \frac{1}{4}g_{ij}g^{sl}(\log{\rho})_{,s}(\log{\rho})_{,l}\notag \\
 & = \frac{n-1}{4}(\log{\rho})_{,i}(\log{\rho})_{,j} + \frac{1-n}{4}g_{ij}(\log{\rho})_{,l}(\log{\rho})_{,l} + \frac{1}{4}g_{ij}\delta_{lk}\delta_{lk}(\log{\rho})_{,l}(\log{\rho})_{,l}\notag \\
 & - \frac{1}{4}\delta_{il}\delta_{jl}(\log{\rho})_{,l}(\log{\rho})_{,l}\notag \\
 & = \frac{n-1}{4}(\log{\rho})_{,i}(\log{\rho})_{,j} - \frac{n-1}{4}g_{ij}(\log{\rho})_{,l}(\log{\rho})_{,l} + \frac{1}{4}g_{ij}(\log{\rho})_{,k}(\log{\rho})_{,k}\notag \\
 & - \frac{1}{4}(\log{\rho})_{,i}(\log{\rho})_{,j}\notag \\
 & = \frac{n-2}{4}(\log{\rho})_{,i}(\log{\rho})_{,j} + \frac{2-n}{4}g_{ij}(\log{\rho})_{,l}(\log{\rho})_{,l}
\end{align}
Where the red terms cancel.
Therefore
\begin{align}
& \widetilde{R}_{ij} = R_{ij} + \frac{2-n}{2}(\log{\rho})_{,ij} - \frac{1}{2}g_{ij}(\log{\rho})_{,kk} + \frac{n-2}{4}(\log{\rho})_{,i}(\log{\rho})_{,j}\notag\\
& + \frac{2-n}{4}g_{ij}(\log{\rho})_{,l}(\log{\rho})_{,l}\notag \\
& = R_{ij} - \frac{n-2}{2}(\log{\rho})_{,ij} + \frac{n-2}{4}(\log{\rho})_{,i}(\log{\rho})_{,j} - \frac{1}{2}\left(\Delta_{g}((\log{\rho})) + \frac{n-2}{2}|\na (\log{\rho})|_{g}^{2} \right)
\end{align}
\end{proof}

\begin{coro}\label{coro: conformal transformation of scalar curvature}
Continue from Proposition $(\ref{proposition: conformal transformation of Ricci curvature})$, then the scalar curvature of $\widetilde{g}$ is given by:
\begin{equation}\label{conformal transformation of scalar curvature}
\color{Black}\widetilde{R} = \frac{1}{\rho}R - \frac{n-1}{\rho^{2}}\Delta_{g}\rho - \frac{(n-1)(n-6)}{4\rho^{3}}|\na \rho|_{g}^{2}. 
\end{equation}
\end{coro}

\begin{proof}
We again compute at $p$, which is the center of a geodesic normal coordinate neighborhood. The result then follows from a simple calculation:

\begin{align}
\widetilde{R}& : = \widetilde{g}^{ij}\widetilde{R}_{ij} = \frac{1}{\rho}g^{ij}\Red\bigg(\Black R_{ij} - \frac{n-2}{2}(\log{\rho})_{,ij} + \frac{n-2}{4}(\log{\rho})_{,i}(\log{\rho})_{,j}\notag \\
&  - \frac{1}{2}\Big(\Delta_{g}(\log{\rho}) + \frac{n-2}{2}|\na (\log{\rho})|_{g}^{2} \Big) \Red\bigg)\Black\notag \\
& = \frac{1}{\rho}R - \frac{n-2}{2\rho}\Delta_{g}(\log{\rho}) + \frac{n-2}{4\rho} g^{ij}(\log{\rho})_{,i}(\log{\rho})_{,j} - \frac{n}{2\rho} \left[\Delta_{g}(\log{\rho}) + \frac{n-2}{2} |\na \log{\rho}|_{g}^{2} \right]\notag \\
& = \frac{1}{\rho}R - \left(\frac{n-2}{2 \rho} + \frac{n}{2\rho} \right) \Delta_{g}(\log{\rho}) + \frac{n-2}{4\rho} |\na \log{\rho}|_{g}^{2} - \frac{n(n-2)}{4\rho}|\na \log{\rho}|_{g}^{2}\notag \\
& = \frac{1}{\rho}R - \frac{2(n-1)}{2 \rho} \Delta_{g}(\log{\rho}) + \frac{-n^{2} +3n - 2}{ 4 \rho} |\na \log{\rho}|_{g}^{2}\notag \\
& = \frac{1}{\rho}R - \frac{n-1}{\rho}\left(\frac{1}{\rho}\rho_{,i} \right)_{,i} + \frac{-n^{2} + 3n - 2}{4\rho} \left(\frac{1}{\rho^{2}}\rho_{,i}\rho_{,i} \right)\notag \\
& = \frac{1}{\rho}R - \frac{n-1}{\rho}\left(-\frac{1}{\rho^{2}}\rho_{,i}\rho_{,i} + \frac{1}{\rho}\rho_{,ii} \right) + \frac{-n^{2} + 3n -2}{4\rho^{3}} |\na \rho|_{g}^{2}\notag \\
& = \frac{1}{\rho}R - \frac{n-1}{\rho^{2}}\Delta_{g}\rho + \frac{4n - 4 -n^{2} + 3n - 2}{4\rho^{3}} |\na \rho|_{g}^{2}\notag \\
& = \frac{1}{\rho}R - \frac{n-1}{\rho^{2}}\Delta_{g}\rho - \frac{(n-1)(n-6)}{4\rho^{3}}|\na \rho|_{g}^{2},\notag 
\end{align}
as desired.
\end{proof}

We want to eliminate the gradient term, and for that purpose we need to distinguish the following two cases.

\begin{coro}\label{coro: conformal transformation of scalar curvature in dimension two}
Suppose the manifold is of dimension $n = 2$. Let $\rho: = e^{2u}$, where $u > 0$ is a smooth function on $M$, then
\begin{equation}\label{conformal transformation of scalar curvature in dimension two}
\color{BlueViolet}\widetilde{R} = e^{-2u}\left(R - 2\Delta_{g}u \right). 
\end{equation}
\end{coro}

\begin{proof}
We plug $\rho = e^{2u}$ into equation $(\ref{conformal transformation of scalar curvature})$, and again, we will be computing at $p\in M$ in a geodesic normal coordinate neighborhood:
\begin{align*}
\widetilde{R} & = \rho^{-1}R - \rho^{-2} \Delta_{g}\rho + \rho^{-3}|\na \rho|^{2}\\
& = e^{-2u}R - e^{-4u}\left(e^{2u}2u_{,i} \right)_{,i} + e^{-6u} (e^{2u}2u_{,i})^{2}\\
& = e^{-2u}R - e^{-4u}\left(\color{Red} e^{2u}2u_{,i}2u_{,i} \color{Black} + e^{2u}2u_{,ii}\right) + \color{Red} e^{-6u}4e^{4u}u_{,i}u_{,i}\\
& = e^{-2u}\left(R - 2\Delta_{g}u \right),
\end{align*}
where the red terms cancel.
\end{proof}

\begin{coro}\label{coro: conformal transformation of scalar curvature in dimension three and higher}
Now suppose manifold is of dimension $n \geq 3$, we set $\rho: = u^{\frac{4}{n-2}}$, for $u$ a positive function on $M$. Then equation $(\ref{conformal transformation of scalar curvature})$ becomes:
\begin{equation}\label{conformal transformation of scalar curvature in dimension three and higher}
\color{BlueViolet}\widetilde{R} = u^{-\frac{n+2}{n-2}}\left(Ru - \frac{4(n-1)}{n-2}\Delta_{g}u \right). 
\end{equation}
\end{coro}

\begin{proof}
The result follows from plugging $\rho = u^{\frac{4}{n-2}}$ into equation $(\ref{conformal transformation of scalar curvature})$, again we will be computing at $p$:

\begin{align}
\widetilde{R} & = u^{-\frac{4}{n-2}}R - (n-1) u^{-\frac{8}{n-2}} \Delta_{g} (u^{\frac{4}{n-2}}) - \frac{(n-1)(n-6)}{4}u^{-\frac{12}{n-2}}|\na (u^{\frac{4}{n-2}})|_{g}^{2}\notag \\
& = u^{-\frac{4}{n-2}} R - (n-1) u^{-\frac{8}{n-2}} \left(\frac{4}{n-2} u^{\frac{6-n}{n-2}} u_{,i} \right)_{,i} - \frac{(n-1)(n-6)}{4} u^{-\frac{12}{n-2}}\left(\frac{4}{n-2}\right)^{2} \notag\\
& \cdot u^{\frac{2(6-n)}{n-2}} u_{,i}u_{,i}\notag \\
& = u^{-\frac{4}{n-2}}R - (n-1)u^{-\frac{8}{n-2}} \frac{4}{n-2} \left(\frac{6-n}{n-2} u^{\frac{8-2n}{n-2}}u_{,i}u_{,i} + u^{\frac{6-n}{n-2}}u_{,ii}\right)\notag \\
& - \frac{(n-1)(n-6)}{4} \left(\frac{4}{n-2}\right)^{2} u^{\frac{-12 + 12 - 2n}{n-2}} u_{,i}u_{,i}\notag \\
& = u^{-\frac{4}{n-2}}R - (n-1) \frac{4}{n-2}u^{-\frac{n+2}{n-2}} \Delta_{g}u\notag \\
& - |\na u|_{g}^{2} \left((n-1)\frac{4}{n-2} \frac{6-n}{n-2} u^{-\frac{2n}{n-2}} + \frac{(n-1)(n-6)}{4}\left(\frac{4}{n-2} \right)^{2} u^{-\frac{2n}{n-2}}\right)\notag \\
& = u^{-\frac{4}{n-2}}R - \frac{4(n-1)}{(n-2)}u^{-\frac{n+2}{n-2}} \Delta_{g}u - |\na u|_{g}^{2} \left(\frac{-4(n-1)(n-6)}{(n-2)^{2}} + \frac{(n-1)(n-6)4^{2}}{4(n-2)^{2}}\right)\notag\\
& \cdot u^{-\frac{2n}{n-2}}\notag \\
& = u^{-\frac{4}{n-2}}R - \frac{4(n-1)}{(n-2)}u^{-\frac{n+2}{n-2}} \Delta_{g}u\notag \\
& = u^{-\frac{n+2}{n-2}}\left(Ru - \frac{4(n-1)}{n-2}\Delta_{g}u \right),
\end{align}
which is desired.
\end{proof} 

\nocite{*} 
\bibliographystyle{plain}
\normalbaselines 
\bibliography{./Bibliography/mybib} 

\biography


Hangjun Xu was born in August, 1987 in Hangzhou, China. In 2005, He went to Zhejiang Unviersity in China and obtained a Bachelor of Science degree in Mathematics and Applied Mathematics with the \emph{Chu Kochen Honors Program} certificate in June, 2009. After graduation, he went to Duke University to pursue a Ph.D. degree in Mathematics under the supervision of Professor Hubert L. Bray. His field of research has been differential geometry and geometric analysis. He developed a side interest in computational geometry, and since fall 2011, he started pursuing a Master of Science degree in Computer Science under the supervision of Professor Pankaj K. Agarwal, en route to his Ph.D. program. During his stay at Duke, he has taught $8$ undergraduate courses as an instructor, including calculus I and II, linear algebra, ordinary and partial differential equations. In 2012, he received the \emph{Graduate School Research Fellowship} for five thousand dollars.  Starting June 2014, he will be a senior software engineer of Oracle at Santa Clara, California. His dream is to become an animation and visual effects artist.

\end{document}